\documentclass{article}
\usepackage{amssymb}
\usepackage{amsthm}
\usepackage{amsmath}
\usepackage{xspace}
\usepackage{url}
\usepackage{cancel}
\usepackage{stmaryrd}
\theoremstyle{plain}
\newtheorem{theorem}[equation]{Theorem}
 \newtheorem{corollary}[equation]{Corollary}
 \newtheorem{lemma}[equation]{Lemma}
 \newtheorem{proposition}[equation]{Proposition}
 \newtheorem{definition}[equation]{Definition}
\theoremstyle{definition}
 \newtheorem{remark}{\emph{Remark}}
 \newtheorem{example}{\emph{Example}}

\newcommand{\mix}{\mathrm{Mix}}
\newcommand{\MLL}{\textbf{MLL}\xspace}

\newcommand{\MLLminus}{\textbf{MLL$^-$}\xspace}

\newcommand{\Cexists}{\mathrm{C}_\exists}
\newcommand{\CP}{\mathrm{C_P}} 
 
\newcommand{\fv}{\mathsf{free}} \newcommand{\bv}{\mathsf{bound}}

\newcommand{\finset}[1]{\{ #1 \}} 
\newcommand{\thehole}[1]{  \finset{  #1  }}
\newcommand{\set}[2]{\{\, #1 \: |\, #2 \, \} }

\renewcommand{\implies}{\Rightarrow}
\newcommand{\LK}{\textbf{LK}\xspace}

\newcommand{\C}{\textsc{C}}

\newcommand{\G}{\mathrm{\Gamma}}

\newcommand{\Cut}{\textsc{Cut}}

\newcommand{\existsR}{\exists\mathrm R}

\newcommand{\forallR}{\forall\mathrm R}

\newcommand{\turnstile}{\vdash}

\usepackage{tikz}
\usetikzlibrary{chains,arrows,decorations.pathmorphing,positioning,patterns,backgrounds,fit,scopes}
\usepackage{fullpage}
\usepackage{lineno}
\input prooftree.sty
\title{Proof nets for Herbrand's Theorem}
\author{Richard McKinley\thanks{Supported by the 
Swiss National Science Foundation grant ``Algebraic and Logical Aspects of Knowledge Processing.'',
and the ANR grant ``INFER''.}}
\begin{document}

 \titlepage
\maketitle
\abstract{This paper explores the connection between two central
  results in the proof theory of classical logic: Gentzen's
  cut-elimination for the sequent calculus and Herbrands ``fundamental
  theorem''.  Starting from Miller's expansion-tree-proofs, a highly
  structured way presentation of Herbrand's theorem, we define a
  calculus of weakening-free proof nets for (prenex) first-order
  classical logic, and give a weakly-normalizing cut-elimination
  procedure.  It is not possible to formulate the usual
  counterexamples to confluence of cut-elimination in this calculus,
  but it is nonetheless nonconfluent, lending credence to the view
  that classical logic is inherently nonconfluent.}
\section{Introduction}
\label{sec:intro}

\renewcommand{\labelitemi}{$\bullet$}

The constructive content of an intuitionistic proof of an existential
statement $\exists x.A$ is well-understood: the \emph{existence
  property} for intuitionistic logic states that a cut-free proof of
$\exists x.A$ is precisely pair of a witness $M$ and a proof of
$A[x:=M]$.  This is, of course, not true for classical logic; a famous
example is the problem ``there exists a pair of irrational numbers $a$
and $b$ such that $a^b$ is rational''.  The standard proof is to first
give $\sqrt 2$, $\sqrt 2$ as a candidate pair.  If $(\sqrt 2)^{\sqrt
  2}$ is rational we are done: if it is irrational, we abandon our
first candidates and instead pick the pair $(\sqrt 2)^{\sqrt 2}, \sqrt
2$.  This is an instance of \emph{backtracking}, and is only possible
because we admit the identity $A = \lnot \lnot A$ on propositions.

In a sense, the counterpart of the existence property in classical
logic is Herbrand's theorem.  In its simplest form, Herbrand's theorem
states that a formula of first-order logic $\exists x.A$, where $A$ is
quantifier free, is provable if and only if there exist formulae $M_1,
\dots M_n$ such that
\[ \vDash A[x:=M_1]\lor \dots \lor A[x:= M_n].\]
\noindent This simple form of Herbrand's theorem does not do the full
generality of the theorem justice, but it gives the a flavour of its
content: a classical proof of an existential does not consist of a
single witness, but a set of candidate witnesses, plus a proof that at
least one of them is an actual witness.  This is complicated by the
fact that the witnesses may interact: observe this in the proof above,
where the failure of $(\sqrt 2,\sqrt 2)$ is necessary to show that
$((\sqrt 2)^{\sqrt 2}, \sqrt 2)$ is a witness.  An example from pure
logic of this interaction is the so-called ``drinker's formula''
\[ \exists x.\forall y (A(x) \to A(y)).\] To prove this formula, we
first guess a witness $a$ (the domain of individuals should be
nonempty, to allow this).  If there is a counterexample (an individual
$b$ such that $A(b)$ does not hold), we \emph{backtrack} and instead
pick $b$ to instantiate the existential quantifier.  That we can
backtrack is expressed logically by \emph{contraction}: we can prove
$\exists x.A$ if and only if we can prove $\exists x.A \lor \exists
x.A$.

It is well known that a more general ``Herbrand's theorem'' for
\emph{prenex} formulae follows directly from Gentzen's cut-elimination
theorem, or more properly the \emph{Midsequent} theorem.  This is
usually stated in terms of permutability of inference rules, but it
can be more succinctly stated as follows:
\begin{theorem}
  The cut-free sequent system given in Figure \ref{fig:midseq} is
  complete for prenex formulae.
\end{theorem}
A proof of a prenex formula $q_1x_1.\dots q_nx_n.B$ in this
calculus yields, for each provable formula, a set of instantiated
versions of $B$ whose disjunction is a tautology (or more generally, a
quasi-tautological consequence of the relevant universal theory).  Any
proof of a prenex formula in the usual sequent calculus may be
converted to a proof in the system in Figure \ref{fig:midseq} by
permuting the quantifier rules below the propositional rules, and then
observing that any consequence of the propositional rules is a sequent
whose disjunction is a tautology.

\begin{figure}
  \centering
  \noindent\hrulefill
  \[
  \begin{prooftree}
    \vDash \bigvee P_i \justifies \turnstile P_1, \dots, P_n
  \end{prooftree}
  \]
  \[
  \begin{prooftree}
    \turnstile \mathbf{\G}, A[x:=a] \justifies \turnstile \G, \forall
    x.A \using \forallR
  \end{prooftree} \qquad \qquad
  \begin{prooftree}
    \turnstile \mathbf{\G}, A[x:= M] \justifies \turnstile
    \mathbf{\G}, \exists x.A \using \existsR
  \end{prooftree}
  \]
  \[
  \begin{prooftree}
    \turnstile\mathbf{\G}, A, \ A \justifies \turnstile\mathbf{\G}, A
    \using C
  \end{prooftree}\qquad
  \begin{prooftree}
    \turnstile\mathbf{\G} \justifies \turnstile\mathbf{\G}, B \using W
  \end{prooftree}
  \]
  \noindent\hrulefill
  \caption{A midsequent calculus $\LK_{mid}$, sound and complete for prenex
    classical logic}
  \label{fig:midseq}
\end{figure}

Herbrand's theorem is usually stated in terms of provability: a
first-order formula (typically $\Pi^0_2$) is provable in a certain
theory if and only if a an extension of that theory can be found such
that some Herbrand disjunction of the formula is provable in the
theory.  But the original theorem \cite{Herbrand30} was stated in
terms of a proof system, with an associated notion of \emph{Herbrand
  proof} \cite{buss95onherb}.  This paper examines these Herbrand
proofs from a modern perspective, answering positively the question:
do Herbrand proofs have syntactic cut-elimination?

The notion of Herbrand proof was adopted, improved and extended by
Miller~\cite{Mil87ComRep}, and it is indeed a variation on his notion
of \emph{expansion-tree proof} which we will take as our proof
objects.  Miller called expansion-trees a ``Compact Representation of
Proofs''; when moving from sequents to expansion trees, a lot of
inessential details regarding the order of application of rules is
discarded.  Another representation of proofs with this property is the
paradigm of Girard's \emph{proof nets} \cite{Gir96ProNetPar}.  In this
paper we make an explicit connection between these two previously
unconnected notions. Our proof objects, Herbrand nets, are proof
nets in the style of Girard; the cut-free nets are essentially
expansion-tree proofs.

We can view the paper from one further perspective: that of
controlling/studying the bad properties of classical sequent systems.
The major hurdle for studying the computational content of the
cut-elimination system of Gentzen is that it lacks both of the usual
``good properties'' of proof systems: it is neither confluent and nor
strongly normalizing, and because of this a proof may in general have
infinitely many normal forms, where normal means cut-free.  One might
observe that many of these normal forms differ only by inessential
details, such as the order of rule applications, yet there is no
universally accepted notion of \emph{equality} (or, better,
\emph{equivalence}) of proofs in classical logic, even in the cut-free
case.

The typical examples of bad behaviour in Gentzen's system arise by
opposing structural rules in a cut, leading to critical pairs which
cannot be resolved, and so an evident solution to this problem is to
orient these critical pairs: such an approach is taken by Curien and
Herbelin's $\bar{\lambda}\mu\tilde{\mu}$-calculus, which gives a kind
of sequent calculus for classical logic which is indeed strongly
normalizing and has natural confluent subsystems, corresponding to
call-by-name and call-by-value evaluation.  Systems like these support
the notion that proofs in classical logic have computational content
as of functional programs with control.  But Curien and Herbelin's
system is very restricted when compared to Gentzen's: a proof in this
system is assigned a term of $\bar{\lambda}\mu\tilde{\mu}$, and rule
permutations that are natural in Gentzen's system do not preserve this
annotation: thus proofs which one would naturally like to identify are
given different behaviours.

It is natural to ask if the kind of orientation found in
$\bar{\lambda}\mu\tilde{\mu}$ can be had without these structural
restrictions.  This paper begins with the observation that the
calculus in Figure \ref{fig:midseq} contains redundant rules: in
particular, the contraction rule is unnecessary on universally
quantified formulae, since the rule for the universal is invertible.
Thus, at least for quantified cut formulae, each formula has a natural
``orientation'' derived from the \emph{polarity} of the cut formula.
(where, roughly speaking, a connective is negative if it has an
invertible rule on the right-hand side of the turnstile, and a formula
negative if its main connective is).  Based on this observation, we
show an annotation of the formulae in proofs of the midsequent
calculus, with the initial idea being that two cut-free proofs are
equivalent if their endsequents receive the same annotation.  We then
consider the annotated sequents themselves as proof structures,
calling the resulting class of nets \emph{Herbrand nets}, and develop
for them a theory of correctness and sequentialization.  By harnessing
the notion of \emph{kingdom}, the smallest subnet containing a certain
formula, we give cut-reduction steps directly on Herbrand nets, and
prove cut-elimination.  We will then see an example proof which,
despite lacking the usual critical pairs, reduces
nondeterministically. We also consider the tactic of duplicating the
largest possible subnet (the \emph{empire}), which we will see may
lead to infinite reduction traces.

\subsection{Related work}\label{sec:related-work}
 Strassburger \cite{Stra09ObsPro} has adapted expansion tree
proofs to give a notion of proof net for second-order propositional
MLL.

Structures similar to those we present here are also studied in
Heijltjes (under the name ``Forest proofs'')
\cite{Hei08ProFor,Hei09ClaPro}, but from a rather different
perspective.  We will discuss in depth the differences in these two
pieces of work later: for now we simply state that our two approaches
represent two different ways to repair an intuitive but flawed idea
for cut-elimination in expansion-tree proofs.  Similar connections
between Herbrand's theorem and abstract proof objects for predicate
logic were suggested in~\cite{DJDHughes06TowHilbProb}.

\section{Preliminary definitions}

\subsection{Prenex formulae of classical first-order logic}
\label{sec:logic-prel}

\newcommand{\arity}{\operatorname{ar}}
\newcommand{\Arity}{\operatorname{Ar}}
\newcommand{\variables}{\mathbf{V}} \newcommand{\funcsymbols}{\mathbf{FS}} \newcommand{\predsymbols}{\mathbf{PS}}
A \emph{signature} $\Sigma =(\variables, \funcsymbols, \predsymbols)$
consists of
\begin{itemize}
\item $\variables$, a countable set of variable symbols.
\item $\funcsymbols$, a countable collection of function symbols,
  together with a function $\arity$ (arity) from $\funcsymbols$ to the
  natural numbers.
\item $\predsymbols$, a countable set of predicate symbols, together
  with a function $\Arity$ from $\predsymbols$ to the natural numbers.
\end{itemize}
A \emph{constant} of a signature $\Sigma$ is a function symbol 
with arity zero.
We will use metavariables $x,y,z,a,b$ to denote variable symbols, $f,g$ to
denote function symbols, and $p,q$ to denote predicate symbols.
The \emph{first-order terms} of $\Sigma$ are given by the following
grammar:
\[ M ::= x \ | \ f(M_1, \dots M_{ar(f)}).\]

\newcommand{\var}{\fv} 

\noindent Given a term $M$, the \emph{free variables of $M$} (written $\var(M)$)
are defined as follows:
\[
\var(x) = \finset{x}, \]
\[ \var(f(M_1, \dots M_n)) = \var(M_1) \cup \dots
\cup \var(M_n).\\
\]

An \emph{atomic formula} can be either positive or negative, and is a
tuple consisting of a polarity from $\finset{+,-}$, a predicate symbol
$p$ of arity $n$, and $n$ terms $M_1, \dots M_n$.  We will write an
atomic formula $(+,p, M_1, \dots, M_n)$ as $p(M_1, \dots M_n)$, and an
atomic formula $(-,q, N_1, \dots, N_n)$ as $\bar{q}(N_1,\dots N_n)$.

The quantifier-free formulae (QFFs) are generated from
the atomic formulae using the connectives $\land$ and $\lor$:
\[ P,Q := \, p(M_1, \dots M_{\Arity(p)}) \ | \ \bar{p}(M_1, \dots
M_{\Arity(p)}) \ | \ (P\lor Q) \ | \ (P \land Q) \]
Notice that we give no explicit connective for negation,
presenting instead the formulae of classical logic 
in \emph{negation normal form}.
Each formula $A$ has a dual formula $\bar{A}$ defined
by \emph{De Morgan duality}:
\[\overline{p(M_1,\dots M_n)}= \bar{p}(M_1,\dots M_n)\qquad
\overline{\bar{p}(M_1,\dots M_n)}= p(M_1,\dots M_n)\]
\[\overline{(P \lor Q)} := (\bar{P} \land
\bar{Q}),\quad \overline{P \land Q} := \bar{P} \lor \bar{Q}.\]

A formula in prenex normal form (or \emph{prenex formula} for short)
is a member of the following grammar, where $x$ ranges over the
variables in $\mathbf{V}$ and $P$ over QFFs:
\[ A ::= P \ | \ \exists x.A \ | \ \forall x.A\] 
The dual of a prenex formulae is defined, as for QFFs, using De Morgan
duality:
\[ \overline{\forall x. A} := \exists x. \bar{A}, \quad
\overline{\exists x. A} := \forall x. \bar{A}\]
We will use $qx$ to refer to an indeterminate quantifier over $x$
(q is either $\forall x$ or $\exists x$).  Given a prenex formula $A =
q_1x_1\dots q_nx_n.P$, we call $P$ the \emph{matrix} of $A$.

The bound and free variables of a prenex formula are defined as usual:
\begin{definition}
  \label{def:free-bound-in-formulae}
  Let $A$ be a formula in prenex normal form.  The set of \emph{free
    variables} $\fv(A)$ of $A$ is a set of variable symbols defined as
  follows:
  \begin{eqnarray*}
    \fv(p(M_1, \dots M_n)) =\fv(\bar{p}(M_1, \dots M_n))  & := & \var(M_1)\cup\dots\cup \var(M_n)\\
    \fv(P \land Q) = \fv (P \lor Q) & := & \fv(P) \cup \fv(Q)\\
    \fv(\forall x.A) = \fv (\exists x.A) & := & \fv(A) \setminus \finset{x}\\ 
  \end{eqnarray*}
  \noindent The set of \emph{bound variables}  
  $\bv(P)$ of a QFF $P$ is empty set.  For an arbitrary prenex
  formula $A$, $\bv(A)$ is the set of variable symbols defined as
  follows:
  \[
  \bv(\forall x.A) = \bv (\exists x.A)  := \bv(A) \cup \finset{x}\\
  \]
\end{definition}
\noindent Notice that, because of the way prenex formulae are built, for any
prenex formula $A$ we have \\$\fv(A)~\cap~\bv(A)~=~\emptyset$. We will
use the notation $A[x:=M]$ for the usual notion of capture avoiding
substitution of a first-order term $M$ for a variable $x$ in a
formula.





\subsection{Trees and terms, forests and sequents}
\label{sec:trees-terms}

\newcommand{\pr}{\mathrm{pr}}

For us, a forest will be a pair $(A, \pr)$ consisting of a set $A$ and
a partial endofunction $\pr$ (predecessor) on $A$ (the elements of $A$
on which $\pr$ is undefined being the \emph{roots} such that, for each
element $x$ of $A$, there is an $n \geq 0$ such that $\pr^n(x)$ is a
root.  Clearly, a forest with one root is a tree.  Given a $y$ such
that $\pr(x)=y$, we will say that $x$ is a \emph{successor} of $y$.  A
forest possesses a natural order structure derived from its
predecessor: $x \leq y$ if there exists $n \geq 0$ with $x = \pr^n
(y).$

The trees we deal with will be derived from subterms of terms or
subformulae of formulae.  For example, given a formula $X$, denote its
set of occurrences of subformulae $\mathcal{O}_X$.  This set has a
natural tree structure (we write $\pr_X$ for its predecessor and
$\leq_{X}$ for its order) whose root is the formula $X$, and where $Z
\leq_{X}Y$ means ``$Y$ is a subformula of $Z$''.  Similarly, if $t$ is
a term, we denote its natural tree structure by
$(\mathcal{O}_t,\pr_t)$, and the order on its subterms $\leq_t$.

We mention here our approach to multisets.  It is usual to define a
sequent as a multiset of formulae, without worrying too much about
what a multiset is, when in fact there are several notions of multiset
with differing properties.  For example, we might consider a finite
multiset of formulae to be a function (multiplicity) from the set of
formulae to the set of the natural numbers which is zero on all but a
finite set of formulae.  This is fine for many applications, but
disastrous for analysing proofs.  For example, if from a sequent $\G,
A$, we derive $\G, A, A$ by weakening and then $\G, A \lor A$ by an
application of a multiplicative $\lor$ rule, how are we to know which
disjunct arose from the weakening?  The right notion of multiset in
this context allows us to distinguish between members of
the multiset, without going so far as to order them (for further discussion,
see~\cite{Lam08EssNets}).  In this paper, we will think of sequents as
\emph{forests} of formulae: the roots of the forest are then distinct
vertices of the forest and can be distinguished.  This corresponds, of
course, to a ``right definition'' of a multiset of formulae, and so
we will write a sequent as $A_1, \dots A_n$ as usual.

\section{Herbrands Theorem and Herbrand proofs}

\newcommand{\alel}{\alpha\varepsilon}
\newcommand{\aec}{$\alel$-calculus}
\newcommand{\univax}{\mathcal{T}}

The form of Herbrand's theorem we will use is the following: let
$\Sigma$ be a signature containing at least one constant, and let
$\univax$ be a finite set of universal axioms.  The Herbrand's theorem
is the following: a prenex formula $A$ over $\Sigma$ is semantically
entailed by $\univax$ ($\univax \vDash A$) if and only if $A$ has an
\emph{Herbrand Proof} \cite{buss95onherb} over $\univax$, which is a
triple consisting of an \emph{expansion}, a \emph{prenexification},
and a \emph{witnessing substitution}.

\begin{remark}
  This restriction to prenex formulae is necessary to have a
  connection between Herbrand proofs and standard sequent proofs;
  Herbrand's theorem for general formulae does not follow directly
  from the midsequent theorem.  To prove Herbrand's theorem for
  general first-order formulae one can consider a generalized
  sequent-calculus with ``deep'' contraction; see
  \cite{McK08Herbnote}.
\end{remark}

We now define the constituents of an Herbrand proof, beginning with expansion:

\begin{definition}
  \label{def:context}
  \begin{enumerate}
  \item A \emph{context} is a prenex formula with precisely one
    occurrence of the special atomic~formula $\thehole{}$ (the
    \emph{hole}).  We write $A \thehole{}$ to denote a context.
  \item If $A\thehole{}$ is a context, $B$ a formula, we write $A\thehole{B}$ for the
    formula given by replacing the hole by $B$.
  \end{enumerate}
\end{definition}
\begin{definition}
  \label{def:expansion}
  Let $A$ be a prenex formula.  An \emph{expansion} of $A$ is defined
  as follows: 
\begin{itemize}
\item 
$A$ is an expansion of $A$; 
\item if $C \thehole{\exists
    x.B}$ is an expansion of $A$, then $C \thehole{\exists
    x.B \lor \exists x.B}$ is an expansion of $A$.
\end{itemize}
\end{definition}
Given an \emph{arbitrary} (not necessarily prenex) formula $A$ of
classical predicate logic, a \emph{prenexification} of a formula $A$
is given by ``pulling the quantifiers to the front'':
\begin{definition}
  \label{def:prenexification}
  Let $A$ be a formula of first-order classical logic.  By renaming
  bound variables, we may write $A$ such that each quantifier $q_i$ in
  $A$ binds a unique variable $x_i$ not appearing free in $A$.  A
  \emph{prenexification} of $A$ is a prenex formula $A^*$ derived from
  this form of $A$ by applications of the following rewrites (where
  $*$ is either $\land$ or $\lor$):
  \begin{eqnarray*}
    A * q x.B \to qx.(A * B) &\quad& qx.A * B \to qx.(A * B)\\
  \end{eqnarray*} 
\end{definition}
A \emph{witnessing substitution} for a closed prenex formula $A$ is a
sequence of terms which, when substituted into the matrix of $A$, make
it valid, and which respects the order of quantifiers appearing in the
prefix of $A$:

\begin{definition}
  \label{def:witnessing-substitution}
  Let $A= q_1x_1.\dots q_nx_n.B(x_1, \dots x_n)$ be a closed prenex formula.  A
  \emph{witnessing substitution} for $A$ is a finite sequence $M_1, \dots M_n$
  of first-order terms such that
  \begin{enumerate}
  \item $M_i = x_i$ if $q_i=\forall$
  \item $\fv(M_j) \subseteq \set{x_i}{q_i = \forall, i<j}$
  \item $\univax \vDash B(M_1, \dots M_n)$
  \end{enumerate}
\end{definition}

Given these components, we may now state Herbrand's theorem:

\begin{theorem}[Herbrand's theorem]
  \label{thm:Herbrand}
  Let $\Sigma$ be a signature containing at least one constant symbol,
  and let $\univax$ be a finite set of universal axioms.  Let $A$ be a
  closed prenex formula; then $\univax \vDash A$ if and only if $A$
  has an \emph{Herbrand proof} -- a tuple $(\hat{A}, A^*, \sigma)$
  such that $\hat{A}$ is an expansion of $A$, $A^*$ is a
  prenexification of $\hat{A}$, and $\sigma$ is a witnessing
  substitution for $A^*$.
\end{theorem}

One direction is easy: if a closed prenex formula $A$ has an Herbrand proof
$(\hat{A}, A^*, \sigma)$ then $\univax \vDash A$. For we have $\univax
\vDash A^*$, since $A^*$ has a witnessing substitution. Furthermore
$\univax \vDash A^*$ if and only if $\univax \vDash \hat{A}$. Lastly, $\hat{A}
\implies A$ is a classically valid implication.  We postpone the
other direction to Section~\ref{sec:cut-free-compl}.

\begin{remark}
  While Herbrand proofs provide a satisfactory abstract account of
  Herbrand's theorem, Herbrand proofs themselves are not a good
  candidate for abstract proof objects, since they lack canonicity.
  Given an Herbrand proof, we can find another with essentially the
  same content by making permutations in the quantifier prefix of the
  prenexification (such that it is still compatible with the
  witnessing substitution).  Such a permutation is the equivalent, in
  this setting, of a sequent calculus rule permutation.  Miller's
  expansion-tree proofs \cite{Mil87ComRep} provide a better notion of
  abstract proof, where a specific prenexification is replaced by a
  demonstration that such a prenexification exists: an acyclicity
  check on the dependencies on quantifiers induced by the
  substitution.  In the following section, we give a reformulation of
  expansion tree proofs and extend them to account for multiple
  conclusions and the presence of cuts.  In the presence of cuts,
  acyclicity is not enough to check correctness; instead, we treat a
  forest of expansion trees as a \emph{proof structure}, and use a
  form of proof-net correctness to identify those corresponding to
  genuine proofs.
\end{remark}

\section{$\alpha\varepsilon$-terms}
\label{sec:alph-calc}

We take, in this paper, the position that a proof net is a forest with
an additional linking structure.  This is most evident in $\MLLminus$
with or without $\mix$, where the forest structure of a net is simply
the forest given by its conclusion, and the linking structure is given
by the axiom links.  Nets for $\MLL$ can also be considered in this
light, with the linking also indicating attachments for $\bot$.  A
pleasing aspect of this approach is that, when considering proofs
which are essentially identical (in this case, by a Trimble rewiring)
the forest remains constant, and only the linking changes.

In settings where we have contraction, such as Lamarche's
essential~nets~\cite{Lam08EssNets}, or Robinson's nets for
propositional classical~logic~\cite{Rob03ProNetCla}, the forest
structure of a net is more complex. In both these settings one finds
\emph{contraction nodes}, of the form
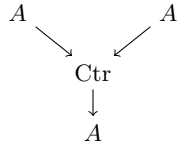
\[ 
\begin{tikzpicture}[level distance=8mm, sibling distance=20mm, edge
    from parent path= {[<-](\tikzparentnode) to (\tikzchildnode)}]]
\node {$A$}[grow=up] 
       child {node {$\mathrm{Ctr}$}
              child {node {$A$}} 
              child {node {$A$}}};
\end{tikzpicture}
\]

\noindent This is problematic from the point of view of canonical
representation of proofs: in addition to any quotienting by rewiring
weakenings (which occurs solely in the linking on the forest), we must
in addition quotient by identities between forests, such as as
those generated by the following identities of subtrees:
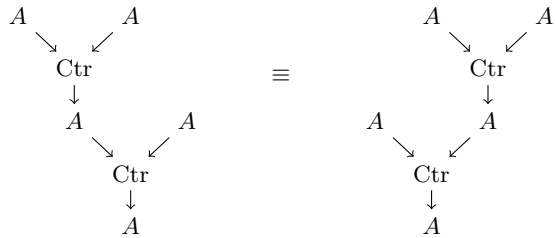
\[
\begin{tikzpicture}[level distance=7mm, sibling distance=15mm, edge
    from parent path= {[<-](\tikzparentnode) to (\tikzchildnode)}]]
\node at (0,0) {$A$}[grow=up] 
      child {node {$\mathrm{Ctr}$}
             child {node {$A$}} 
             child {node {$A$} 
                    child {node {$\mathrm{Ctr}$}
                           child {node {$A$}} 
                           child {node {$A$}}}}};
\node at (2,2) {$\equiv$};
\node at (4,0) {$A$}[grow=up] 
      child {node {$\mathrm{Ctr}$}
             child {node {$A$} 
                    child {node {$\mathrm{Ctr}$}
                           child {node {$A$}} 
                           child {node {$A$}}}}
              child {node {$A$}}};
\end{tikzpicture}
\]
and
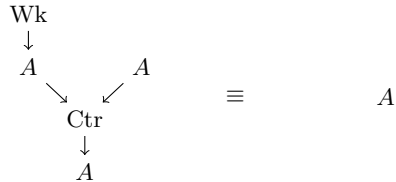
\[
 \begin{tikzpicture}[level distance=7mm, sibling distance=15mm, edge
     from parent path= {[<-](\tikzparentnode) to (\tikzchildnode)}]]
 \node at (0,0) {$A$}[grow=up] 
       child {node {$\mathrm{Ctr}$}
              child {node {$A$}} 
              child {node {$A$} 
                     child {node {$\mathrm{Wk}$}}}};
 \node at (2,1) {$\equiv$};
 \node at (4,1) {$A$};
\end{tikzpicture}
\]

\noindent Such identifications are necessary, for example, to build a
sensible category from Robinson's proof nets, see \cite{Fuhr06OrdEnr}.
The forest structure we suggest for representing contraction replaces
the usual binary contraction with an n-ary contraction. The suggestion
to use an n-ary contraction is already present in
\cite{Girard91NewCon}, along with an additional condition to enforce
canonicity: a structural rule cannot provide the premise of another
structural rule.

\subsection{$\alpha\varepsilon$ terms}
In this section we define $\alpha\varepsilon$-terms, which 
consist of the \emph{expansion-trees} (a reformulation
of Miller's expansion trees for the prenex first-order fragment
of classical logic), \emph{cuts}, and \emph{witnessing terms}.

\label{sec:alph-terms}
\newcommand{\indices}{\mathbf{I}}
\begin{definition}[$\alpha\varepsilon$ terms]
  \label{def:alel-terms}
  Let $\Sigma = (\variables, \funcsymbols, \predsymbols)$ be a
  signature, and let $\indices$ be a countable set of indices.  The
  \emph{$\alel$ terms} $t, \dots$ over $(\Sigma, \indices)$
  (consisting of the \emph{expansion trees} $p, \dots$, \emph{cuts}
  $c, \dots $, and \emph{witnessing terms} $w, \dots $) are given by
  the following grammars:
  \[ t := p \ | \ w \ | \ c\]
  \[ p := S \ | \ \alpha[a].p \ | \ ( w + \dots + w) \]
  \[ w := \varepsilon[M].p \]
  \[ c := p \bowtie p\]
  \noindent where $S$ is a nonempty finite set of indices, $M$ is a
  first-order term over the signature, $a \in \variables$, and $( w
  + \dots + w )$ denotes a finite nonempty formal sum (a member of the
  free commutative semigroup over $w$).  A \emph{non-cut term} is either an
  expansion tree or a witnessing term.
\end{definition}

The witnessing terms represent the components of (generalized)
Herbrand disjunctions.  We could of course replace the formal sums of
witnesses by nonempty finite multisets of witnesses, but this
complicates the notation a little.  The reader more comfortable with
multisets can think of $(w_1+\cdots w_n)$ as the multiset
$\finset{w_1, \dots, w_n}$, with the semigroup operation $+$ being
interpreted as disjoint multiset union. We make an explicit
distinction between the witnessing term $\varepsilon[M].t$ and the
expansion tree $(\varepsilon[M].t)$.  We will refer to a witnessing
term not in the scope of a semigroup $+$ as a \emph{naked} witness.

\begin{remark}
  \label{remark:semigroup}
  The reader might wonder why we have a commutative semigroup rather
  than commutative monoid structure on expansion trees: why are we not
  allowed to form the empty formal sum as a expansion tree (in
  multiset terms, why not also allow the empty multiset)?  This would
  amount to explicit weakening in our sequent calculus, and in the
  proof nets we will form from $\alpha\varepsilon$ terms.  Weakening
  is notoriously difficult to handle well in proof nets; as we will
  see, in this setting explicit weakening is not necessary.
\end{remark}

\subsection{Typing $\alel$-terms}
\label{sec:typing-alel-terms}
We now assign \emph{types} to these terms. The type of an expansion tree is
always a prenex formula. The witnessing terms and cuts receive special
non-logical types:

\newcommand{\witnesstype}[1]{\langle \exists #1 \rangle}

\begin{definition}
  \label{def:prenex-type}
  A \emph{type} over $\Sigma =
  (\mathcal{X}, \mathcal{F}, \mathcal{R})$ is either
  \begin{enumerate}
  \item A \emph{logical type}: a formula of classical predicate logic
    in prenex normal form, over the signature (as described in the
    preliminaries); or
  \item a \emph{non-logical} type, of which there are two
    kinds: \begin{enumerate}
     \item A \emph{witness type}, written $\witnesstype{x.A}$,
      where $\exists x.A$ is a formula in prenex normal form; or
     \item A \emph{cut type}: a pair of dual formulae of classical
      logic in prenex normal form, written $A \bowtie \bar{A}$.
    \end{enumerate}
  \end{enumerate}
  We will occasionally need to refer to a type without specifying if
  is is logical or non-logical: in that case we will use a capital
  $T$, reserving $A, B, C \dots$ for those types which are prenex
  formulae.
\end{definition}

The non-logical types are needed to type the witness and cut terms,
respectively. We use the witness types to distuish between a naked
witness, $\varepsilon[M].s$, which recieves a witness type, and the
expansion tree $(\varepsilon[M].s)$, which recieves a logical type.
Only terms of witness type can take part in an expansion, and only
terms of logical type can take part in other logical rules; in this
way, we avoid non-canonicity, preventing the premise of an expansion
to be, itself, the result of an expansion.

Each non-logical types has an underlying logical type:
\begin{definition}
  \label{def:underlying-type}
  The \emph{underlying type} of a witness type $\witnesstype{x.A}$ is
  $\exists x.A$.  The \emph{underlying type} of $A \bowtie \bar{A}$ is
  $A$.  The free/bound variables $\fv$ and $\bv$ of a witness/cut type
  are the free/bound variables of its underlying type.  We define
  substitution into witness/cut types in the obvious way
  \[ \witnesstype{x.A}[y:=M] = \witnesstype{x.A[y:= M]}\]
  \[ (A \bowtie \bar{A}) [y:=M] = A [y:=M]\bowtie \bar{A}[y:=M] \]
\end{definition}

\begin{definition}
  \label{def:typed-term}
  A typed term is a pair $t:T$ of a term $t$ and a type $T$, derivable
  in the typing system given in Figure~\ref{fig:LKH}.
\end{definition}

\begin{figure}[t]
  \noindent\hrulefill
  \centering
  \[
  \begin{prooftree}
    i_1, \dots i_n \in \indices
    \justifies \finset{i_1, \dots i_n} : P
  \end{prooftree}
  \]
  \[
  \begin{prooftree}
    t : A[x:=a] \justifies \alpha [a].t : \forall x.A
  \end{prooftree} \qquad\qquad \qquad
  \begin{prooftree}t : A[x:= M] \justifies \varepsilon[M].t: \witnesstype{
    x.A} \end{prooftree}\]\vspace{0.7em}
  \[
  \begin{prooftree} w_1 : \witnesstype{ x.A}, \ \dots, \ w_n:\witnesstype{ x.A}
    \justifies (w_1+ \dots + w_n) : \exists x.A \end{prooftree}
  \] \vspace{0.7em}
  \[
  \begin{prooftree}
    t: A \quad s: \bar{A} \justifies t \bowtie s : A\bowtie \bar{A}
  \end{prooftree}
  \]
  \noindent\hrulefill
  \caption{Typing derivations for $\alel$ terms}
  \label{fig:typing-derivations}
\end{figure}

There are some terms that cannot be typed, for simple reasons.  For
example, the term $\alpha[a].t\bowtie\alpha[b].s$ can never be
well-typed: a type for a term beginning with an $\alpha$ must be a
formula of the form $\forall x.A$, and two such formulae can never be
dual.

\begin{example}
  \label{ex:drinkers-term}
  The following is a well-typed term, which will be an important
  example for us for the rest of the paper.  Its type is the
  drinker's formula mentioned in the introduction: for that reason we
  will call it $D$, the drinker's term:
  \[
  D = (\varepsilon[\mathsf{c}].\alpha[a].\finset{1} \ + \
  \varepsilon[a].\alpha[b].\finset{1}) : \exists x. \forall y (\bar
  A(x)\lor A(y))\]
\end{example}

We can now take advantage of the fact that terms can be seen as trees:
Figure \ref{fig:extrees} gives typing tree equivalents of the
derivations in Figure~\ref{fig:typing-derivations}.  Viewed in this
way, a typed term is a tree built from the elements in Figure
\ref{fig:extrees} by matching the types of the directed
edges. Conversely, each typed term gives rise to a tree of its typed
subterms: in practice, we will annotate only the root with its type,
as the types of subterms can be inferred.

\begin{example}
  \label{ex:drinkers-tree}
  The tree corresponding to the typed drinker's term $D$ is
  
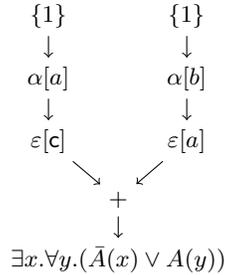
\[
  \begin{tikzpicture}
    \matrix[row sep=3mm,column sep=-10mm]{
      \node (left prop) {$\finset{1}$}; &&
      \node (right prop) {$\finset{1}$}; \\
      \node (left alpha) {$\alpha[a]$}; &&
      \node (right alpha) {$\alpha[b]$}; \\
      \node (left epsilon) {$\varepsilon[\mathsf{c}]$}; &&
      \node (right epsilon) {$\varepsilon[a]$}; \\
      & \node(plus){$+$};\\
      & \node(conclusion){$\exists x. \forall y. (\bar{A}(x) \lor A(y))$} ;\\
    }; \draw[->](left prop) to (left alpha); \draw[->](right prop) to
    (right alpha); \draw[->](left alpha) to (left epsilon);
    \draw[->](right alpha) to (right epsilon); \draw[->](left epsilon)
    to (plus); \draw[->](right epsilon) to (plus); \draw[->](plus) to
    (conclusion);
  \end{tikzpicture}
  \]
\end{example}

As mentioned above, we consider sequents
to be forests of formulae.  Continuing this convention, a proof
structure (or prenet) in this system will consist of a \emph{forest}
of typed expansion-trees and cuts, satisfying certain properties.  By
a forest of typed terms, we mean informally a multiset of typed terms,
and more formally, the following:

\begin{definition}
A \emph{typed forest} $F$ is a forest in which each root is a typed 
term, with the tree above a roots being its typing tree.
\end{definition}

\begin{figure}[t]
  \noindent\hrulefill
  \center
  \begin{tikzpicture}[level distance=8mm, sibling distance=20mm, edge
    from parent path= {[<-](\tikzparentnode) to (\tikzchildnode)}]]
    \matrix[row sep= 1cm,column sep= 1cm]{ \node
      {$\witnesstype{x.A}$}[grow=up] child {node {$\varepsilon[M]$}
        child {node {$A[x:=M]$}} }; &

      \node {$P$}[grow=up] child {node {$S$}}; &

      \node {$\forall x.A$}[grow=up] child {node {$\alpha[a]$} child
        {node {$A[x:=a]$}}
      };\\

      \node {$\exists x.A$}[grow=up] child {node (plus) {$+$} child
        {node {$\witnesstype{ x.A}$}} child {node {$\witnesstype{ x.A}$}} };
      \node at (0,1.5) {\dots};&&

      \node {$A \bowtie \bar{A}$}[grow=up] child {node {$\bowtie$}
        child {node {$\bar A$}} child {node {$A$}}
      };\\
    };

  \end{tikzpicture}

  \noindent\hrulefill
  \caption{Typing trees for $\alel$ terms}
  \label{fig:extrees}
\end{figure}
\begin{figure}
\noindent\hrulefill
  \centering
  \[
  \begin{prooftree}
    \vDash \bigvee_{j=1}^n P_j
    \justifies
    \finset{i} : P_1, \dots,  \finset{i} : P_n
    \using i
  \end{prooftree} 
  \]\vspace{0.7em}
  \[
  \begin{prooftree}
    \turnstile F, \ t : A[x:=a] 
    \justifies F,  \ \alpha [a].t : \forall x.A
    \using \forallR
  \end{prooftree} \qquad \quad
  \begin{prooftree}
     F, \ t : A[x:= M] \justifies 
     F, \ (\varepsilon [M].t): \exists x.A
    \using \existsR
  \end{prooftree} \]\vspace{0.7em}
  \[
  \begin{prooftree}
     F, \ t: \exists x.A, \ s:\exists x.A
    \justifies 
     F, \ t + s : \exists x.A
    \using  \Cexists
  \end{prooftree}\qquad\qquad
\begin{prooftree}
     F, \ S: P, \ T: P
    \justifies 
     F, \ S \cup T : P
    \using \CP
  \end{prooftree}
  \]\vspace{0.7em}
  \[
  \begin{prooftree}
     F, \ t: A \quad  G, \ s: \bar{A} 
    \justifies 
     F, \ G, \ t \bowtie s : A \bowtie \bar{A}
    \using \Cut
  \end{prooftree}
  \]
\noindent\hrulefill
  \caption{$\LK_H$: An annotated sequent calculus for prenex 
classical logic}
  \label{fig:LKH}
\end{figure}
\subsection{Decorating sequent derivations with terms}
\label{sec:decor-sequ-deriv}
We now use $\alpha\varepsilon$ terms to decorate the formulae
appearing in sequent proofs of classical logic, just as one may assign
lambda terms to proofs of intuitionistic logic.  This provides an
elegant assignment of typed forests to proofs.  The proofs we annotate
will be of a particular form; we restrict the system $\LK_{mid}$ in
Figure~\ref{fig:midseq} to a subsystem in which weakening does not
appear and contraction is restricted to existential and
quantifier-free formulae, and add term annotations. This system with
term annotations, $\LK_H$ is given in Figure~\ref{fig:LKH}.  In the
next section we will show that this system is complete for prenex
formulae, and in so doing give a function assigning a typed forest to
any proof in $\LK_{mid}$.

The rules of $\LK_H$ operate not on sequents, but on the typed
forests introduced in the previous section. The rule
\[\begin{prooftree}
  \vDash \bigvee_{j=1}^n P_j  \justifies  (i) : P_1, \dots, (i):P_n
  \using i
\end{prooftree}
\]
\noindent is the \emph{tautology rule}; it allows us to use any
propositional tautology as an axiom, where the formulae $P_j$ are the
QFFs of $\Sigma$.  This is the tautology rule in the absence of 
axioms.  Given a finite set $\univax$ of universal axioms,
we can replace the tautology rule with 
\[\begin{prooftree}
  \univax \vDash \bigvee_{j=1}^n P_j \justifies  (i) : P_1, \dots, (i):P_n
  \using i
\end{prooftree}
\]
to give a calculus for proofs in the theory $\univax$.  Once we do
this, it is easy to see that the \emph{forgetful projection} of a rule
in $\LK_H$ (where we simply remove the annotating term from every
formula, and remove all cut-terms) is a rule of $\LK_{mid}$.

For the annotation of formulae to be well behaved (in a sense that
will be explained below), we must treat eigenvariables
\emph{strictly}: each instance of the universal quantifier should have
a unique associated eigenvariable, and that eigenvariable should only
appear free in the subproof above the rule introducing that
quantifier.  We will also insist that each instance of the tautology
rule has a unique index.

\begin{definition}
  \label{def:LKH-proof}
  A proof in $\LK_H$ is a tree built from rule instances from Figure
  \ref{fig:LKH}, with instances of the tautology rule at the leaves.
  A proof $\Phi$ is \emph{strict} if
  \begin{enumerate}
  \item[(i)] each tautology rule in $\Phi$ is labelled with a distinct
    index $i$,
  \item[(ii)] each $\alpha$ in $\Phi$ binds a distinct eigenvariable
    $a$,
  \item[(iii)] An eigenvariable $a$ may not appear free in the type of
    any sequent outside the subproof above the rule introducing
    $\alpha[a]$.
  \end{enumerate}
  We will write $\LK_H \turnstile F$ if there is a \emph{strict} proof
  in $\LK_H$ of $F$.
\end{definition}
Note that case (ii) in the above definition ensures that
eigenvariables are used strictly in the usual sense, and additionally
enforces the usual variable restriction on the rule for the universal
quantifier.

\begin{example}
  \label{ex:drinker-derivation}
  Let $\Sigma$ contain the unary predicate $A$ and a constant symbol
  $\mathsf{c}$.  Recall the drinker's term $D$
  (Example~\ref{ex:drinkers-term}):
  \begin{equation}
    D= (\varepsilon[\mathsf{c}].\alpha[a].\finset{1} + \ \ \varepsilon[a].\alpha[b].\finset{1}) :
    \exists x. \forall y (\bar A(x)\lor A(y))
    \label{eq:drinker}
  \end{equation} 
  $D$ is the conclusion of the derivation below:
  \[
  \begin{prooftree}
    \[
    \[
    \[
    \[
    \[
    \justifies  \finset{1}:\bar A(\mathsf{c})\lor A(a), \ \
    \finset{1}: \bar A(a)\lor A(b) \using 1
    \]
    \justifies  \finset{1}:\bar A(\mathsf{c})\lor A(a), \ \
    \alpha[b].\finset{1} : \forall y \bar A(a)\lor A(y) \using
    \forallR
    \]
    \justifies  \finset{1}:\bar A(\mathsf{c})\lor A(a), \ \
    (\varepsilon[a].\alpha[b].\finset{1}) : \exists x. \forall y (\bar
    A(x)\lor A(y)) \using \existsR
    \]
    \justifies  \alpha[a].\finset{1}: \forall y (\bar
    A(\mathsf{c})\lor A(y)), \ \ (\varepsilon[a].\alpha[b].\finset{1})
    : \exists x. \forall y (\bar A(x)\lor A(y)) \using \forallR
    \]
    \justifies 
    (\varepsilon[\mathsf{c}]\alpha[a].\finset{1}):\exists x. \forall y
    (\bar A(x)\lor A(y)), \ \ \ (\varepsilon[a].\alpha[b].\finset{1})
    : \exists x. \forall y (\bar A(x)\lor A(y)) \using \existsR
    \]
    \justifies 
    (\varepsilon[\mathsf{c}].\alpha[a].\finset{1} \ + \
    \varepsilon[a].\alpha[b].\finset{1}) : \exists x. \forall y (\bar
    A(x)\lor A(y)) \using \Cexists
  \end{prooftree}
  \]
\end{example}

The following example illustrates the cut rule and the contraction
rule on QFFs:
\begin{example}
  \label{ex:second-example}
  \[
  \begin{prooftree}
    \[
    \[
    \[
    \justifies \finset{1}:P, \ \finset{1}: \bar P, \ \finset{1}:P
    \using 1
    \]
  \qquad 
    \[ 
    \justifies 
    \ \finset{2}:P, \ \finset{2}: \bar P, \ \finset{2}:\bar{P}
    \using 2 
    \]
    \justifies
  \finset{1}:P, \ \finset{2}:P, \ \finset{1}:\bar{P}, \
  \finset{2}:\bar{P},\ \ \finset{1}\bowtie\finset{2}: P \bowtie \bar{P}
  \using \Cut
  \]
  \justifies
  \finset{1}:P, \ \finset{2}:P, \ \finset{1,2}:\bar{P}, \ \ \finset{1}\bowtie\finset{2}: P \bowtie \bar{P}
  \using
  \CP
  \]
\justifies
 \finset{1,2}: P, \ \finset{1,2}:\bar{P}, \ \ \finset{1}\bowtie\finset{2}: P \bowtie \bar{P}
\using \CP
\end{prooftree}
\]
\end{example}

\subsection{Annotated sequents }

We can view the conclusion of a strict proof as a normal sequent
\emph{annotated} with some information about the rules used in the
proof.  Clearly, not every typed forest can arise from annotating a
sequent proof.  For example, there cannot be two occurences of the
prefix $\alpha[a]$ in such a conclusion.  We characterize in the
current subsection the typed forests having ``the right shape'' to
arise from a proof: our equivalent of a proof-stucture/preproofnet.
We will call these typed forests \emph{annotated sequents}.  Then, in
Section 7, we will give a correctness criterion singling out among
these annotated sequents the \emph{Herbrand nets}: those annotated
sequents which do indeed arise from a sequent proof.

\newcommand{\subt}{\mathrm{Subt}} 


We will refer to the subtrees of a forest as its ``nodes''
using standard proof-net terminology.  We will refer
to nodes by the outermost term-constructor used to form
them: so a subterm of the form: 
\begin{itemize}
\item $S$ is a propositional node;
\item $\alpha[a].t$ is an $\alpha$-node;
\item $\varepsilon[M].t$ is an $\varepsilon$-node;
\item $(w_1+\dots w_n)$ is an \emph{expansion} node.
\end{itemize}

We begin by giving a notion of \emph{type} to typed forests; for this,
the standard set-with-multiplicities definition of multiset suffices:

\begin{definition}
The \emph{type} of a typed forest is the multiset consisting of the
types of its non-cut roots.
\end{definition}

\newcommand{\free}{\mathsf{free}_\alpha}
\newcommand{\bound}{\mathsf{bound}_\alpha}

A typed forest will be ``of the right shape'' if it can be determined
that, if it did come from an $\LK_H$ proof, that proof was strict;
that is, we need a notion of strictness for typed forests.  The
essence of strictness in forests is that $\alpha$ behaves as a
non-local \emph{binder}.  The first condition is obvious: each
$\alpha$ should have a unique eigenvariable, just as it does in a
strict proof.



It is a little harder to capture the idea other aspects of strictness:
we begin by defining the $\alpha$-bound and $\alpha$-free variables of
a typed term/forest, a concept distinct from the bound/free variables
of its conclusion:

\begin{definition}
\label{def:alpha-free-bound}
  Let $t:A$ be a typed term.  We define two sets of variables $\bound{(t:A)}$ (the
  variables \emph{$\alpha$-bound} in $t:A$) and {$\free{(t:A)}$} (the
  $\alpha$-free variables of $t:A$) as follows:
  \begin{enumerate}
  \item The variable $a$ is a member of $\bound(t:A)$ if and only if
$\alpha[a]$ is a subterm of $t$.
\item  The set $\free(t:A)$ is defined as 
follows:

\begin{itemize}
\item $\free(S : P) = \fv(P)$
\item $\free (\alpha[a].t:\forall x.A) = \free(t:A[x:=a]) \setminus \finset{a}$ 
\item $\free (\varepsilon[M].t:\witnesstype{x.B}) = \free(t:B[x:=M])\cup \var(M)$ 
\item $\free ( (t_1 + \dots + t_n):\exists x.B) = \free(t_1:\witnesstype{x.B})\cup \dots \cup \free(t_n :\witnesstype{x.B})$\vspace{3pt}
\item $\free (t \bowtie s: A \bowtie \bar{A}) =  \free(t:A)\cup \free(s:\bar{A})$\vspace{3pt}
\end{itemize}
  \end{enumerate}

\end{definition}

\begin{example}
\label{ex:1}
  For the typed expansion tree
  \[t:A = (\varepsilon[b].\alpha[a].(\varepsilon[a].\finset{1})):\exists x.\forall
  y. \exists z.P(x,y,z,w)\]
  \noindent $\free(t:A) = \finset{b, w}$ and $\bound(t:A)=\finset{a}$.
\end{example}

Finally, consider the annotated sequent 
\[(\varepsilon[a].\finset{1}):\exists x.\bar P, \quad
\alpha[a].\finset{1}:\forall x.P \] 
Suppose that this were the conclusion of a strict sequent proof: then
it ended with the application of a $\forallR$ rule with eigenvariable $a$,
so the variable $a$ is \emph{bound} in this typed forest (despite being free
in the term $(\varepsilon[a].\finset{1}):\exists x.\bar P$).  With
these intuitions in place, we define the $\alpha$-free and 
$\alpha$-bound variables of a multiset of typed terms:

\begin{definition}
\label{alpha-free-bound-in-multiset}
  Let $F$ be a typed forest.  A variable $a$ is
  \emph{$\alpha$-bound in $F$ ($a \in \bound(F)$)} if it is in
  $\bound(t:A)$, for some term $(t:A)$ in $F$.  The variable $a$ is
  \emph{$\alpha$-free in $F$} ($a \in \free(F)$) if it is in $\free
  (t:A)$, for some term $(t:A)$ in $F$, and not $\alpha$-bound in $F$.
\end{definition} 

\begin{example}
\label{ex:struct-semi-struct}
In the typed forest
\[ F = \alpha[a].\finset{1}:\forall x.P, \ \
(\varepsilon[a].(\varepsilon[b].\finset{1})): \exists y. \exists z. Q\]
$a \in \bound(F)$ and $b \in \free(F)$.  $a \notin
\free(F)$: although $a$ is $\alpha$-free in the second typed expansion
tree, it is $\alpha$-bound in the first.
\end{example}

We are now in a place to define strictness for typed forests:

\begin{definition}[Strictness]
\label{def:semi-structure}
A typed forest is \emph{strict} if
\begin{enumerate}
\item each $\alpha$ has an eigenvariable, and 
\item for each non-cut root
$t:A$ of $F$, $\bound(F) \cap \fv(A) = \emptyset$.
\end{enumerate}
\end{definition}

So far, we have allowed typed forests containing naked witnesses: of
course, the witness types are not part of first-order logic, and so we
are particularly interested in forests without naked witnesses, as
they will have types consisting of multisets of prenex formulae:

\
\begin{definition}
  \label{def:1}
  An \emph{annotated sequent} is a strict typed forest with no naked witnesses.
\end{definition}

\begin{example}
\label{ex:3}
  The typed forest
  \[ \alpha[a].(\finset{1}):\forall x.P, \quad
  \varepsilon[a].\finset{1}:[\exists x.P] \] is strict, but
  not an annotated sequent, as the typed term $\varepsilon[a].(1):[\exists
  x.P]$ is neither an expansion tree nor a cut.  Neither of the following is
  strict:
 \[\alpha[a].\finset{1}:\forall x.P, \quad \alpha[a].\finset{1}:\forall x.Q \quad \text{  (fails condition (a))}\]
 \[\alpha[a].\finset{1}: \forall x.P,\quad \finset{1}: Q[x:=a] \quad \text{ (fails condition (b))}\]
On the other hand, 
\[\alpha[a].\finset{1}: \forall x.P, \quad \finset{1} \bowtie
\finset{2}:Q[x:=a] \bowtie \bar{Q}[x:=a] \] 
\noindent \emph{is} strict (and so an annotated sequent), since an $\alpha$-bound variable may
appear free in the type of a cut without violating (b).
\end{example}

\begin{proposition}
\label{prop:LKH-derives-structures}
The conclusion of a strict $\LK_H$ derivation is an annotated sequent.
\end{proposition}
\begin{proof}
  By induction: the conclusion of a tautology rule is an annotated
  sequent, and each rule of $\LK_H$ takes annotated sequents to
  annotated sequents.
\end{proof}

\label{ex:cut-elim-example}
We can represent the forest structure of annotated sequents
graphically using the graphical representation of terms: the following
annotated sequent will be our principal example for demonstrating the
cut-elimination theory (see Section~\ref{sec:minim-reduct-not}):
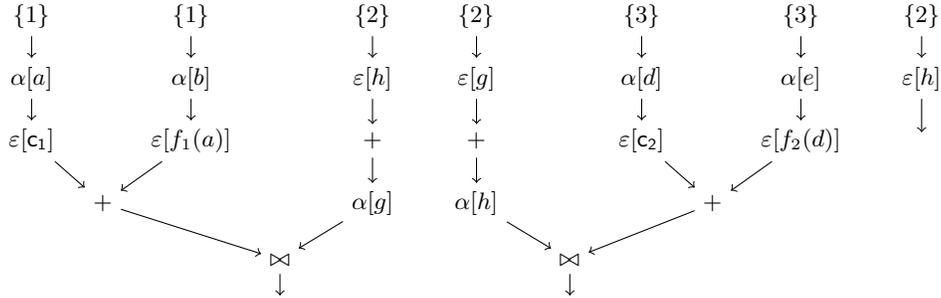
\begin{example}
 \[
\begin{tikzpicture}
\matrix[row sep=3mm,column sep=3mm]{
 \node (left prop) {$\finset{1}$}; &&
 \node (right prop) {$\finset{1}$}; &&& \node (left prop 2) {$\finset{2}$}; &&
 \node (right prop 2) {$\finset{2}$};&&& \node (left prop 3) {$\finset{3}$}; &&
 \node (right prop 3) {$\finset{3}$}; && \node(prop 4){$\finset{2}$};\\
\node (left alpha) {$\alpha[a]$}; &&
 \node (right alpha) {$\alpha[b]$}; &&& 
\node (left epsilon 2) {$\varepsilon[h]$}; &&
 \node (right epsilon 2) {$\varepsilon[g]$};&&&
 \node (left alpha 3) {$\alpha[d]$}; &&
 \node (right alpha 3) {$\alpha[e]$}; && \node(conclusion){$\varepsilon[h]$};\\
\node (left epsilon) {$\varepsilon[\mathsf{c_1}]$}; &&
 \node (right epsilon) {$\varepsilon[f_1(a)]$}; &&&
\node (left plus 2) {+}; &&
 \node (right plus 2) {+}; &&&
\node (left epsilon 3) {$\varepsilon[\mathsf{c_2}]$}; &&
 \node (right epsilon 3) {$\varepsilon[f_2(d)]$}; && \node (out) {};\\
 & \node(plus){$+$}; &&&&
\node (left alpha 2) {$\alpha[g]$}; &&
 \node (right alpha 2) {$\alpha[h]$}; &&&&
\node(plus 3){$+$};\\
&&& \node (left cut) {$\bowtie$}; &&&&&& \node (right cut) {$\bowtie$};\\
&&& \node (left out) {}; &&&&&& \node (right out) {};\\
};



\draw[->](prop 4) to (conclusion);
\draw[->](conclusion) to (out);
\draw[->](plus) to (left cut);
\draw[->](left alpha 2) to (left cut);
\draw[->](right alpha 2) to (right cut);
\draw[->](plus 3) to (right cut);
\draw[->](right cut) to (right out);
\draw[->](left cut) to (left out);
\draw[->](left prop) to (left alpha);
\draw[->](right prop) to (right alpha);
\draw[->](left alpha) to (left epsilon);
\draw[->](right alpha) to (right epsilon);
\draw[->](left epsilon) to (plus);
\draw[->](right epsilon) to (plus);

\draw[->](left prop 2) to (left epsilon 2);
\draw[->](right prop 2) to (right epsilon 2);
\draw[->](left epsilon 2) to (left plus 2);
\draw[->](right epsilon 2) to (right plus 2);
\draw[->](left plus 2) to (left alpha 2);
\draw[->](right plus 2) to (right alpha 2);

\draw[->](left prop 3) to (left alpha 3);
\draw[->](right prop 3) to (right alpha 3);
\draw[->](left alpha 3) to (left epsilon 3);
\draw[->](right alpha 3) to (right epsilon 3);
\draw[->](left epsilon 3) to (plus 3);
\draw[->](right epsilon 3) to (plus 3);

\end{tikzpicture}
\]
\end{example}



\subsection{Alpha renaming} 
We have mentioned that the $\alpha$ should be though of as a non-local
binder; so, in fact, should the indices used to annotate the
conclusions of the tautology rule.  As such, we will need to be able
to rename eigenvariables and indices.  We will use the notation
$[a\leftarrow b]$ to denote the renaming of an $\alpha$-bound
variable, and $[i \leftarrow j]$ for the renaming of an index $i$.

\begin{definition}
\label{def:renaming}
\begin{enumerate}
\item Let $i$ and $j$ be indices.  The operation ${[i
    \leftarrow j]}$ (tautology renaming) is defined as follows:
  \begin{align*}
    S[i \leftarrow j] & = \begin{cases} S & i \notin  S \\
    \ ((S \setminus \finset{i} )\cup \finset{j}) & i \in S \end{cases}\\
    (\alpha[d].t)[i\leftarrow j] & = 
    \alpha[d].(t[i\leftarrow j]) \\
    (\varepsilon[M].t)[i\leftarrow j] & =  \varepsilon[M].(t[i\leftarrow j])\\
    (t_1 +\dots+ t_n)[i\leftarrow j]&  =  (t_1[i\leftarrow j], \dots t_n[i\leftarrow j]) \\
    (t \bowtie s)[i\leftarrow j] & =  (t[i\leftarrow j] \bowtie s[i\leftarrow j])\\
  \end{align*}
 \item Let $a$ and $b$ be members of $\variables$.  The
  operation $[a \leftarrow b]$ (variable renaming) is defined as
  follows:
  \begin{align*}
    S[a \leftarrow b] & = S\\
    (\alpha[a].t)[a\leftarrow b] & =
    \alpha[b].(t[a\leftarrow b]) \\
    (\alpha[d].t)[a\leftarrow b] & =
    \alpha[d].(t[a\leftarrow b]) \qquad d \neq a\\
    (\varepsilon[M].t)[a\leftarrow b] & = 
    \varepsilon[M[a:=b]].(t[a\leftarrow b])\\
    (t_1 + \dots + t_n)[a\leftarrow b] & =  \finset{t_1[a\leftarrow b], \dots t_n[a\leftarrow b]} \\
    (t \bowtie s)[a\leftarrow b] & =  (t[a\leftarrow b] \bowtie s[a\leftarrow b])\\
  \end{align*}
\end{enumerate}
\end{definition}

Renaming respects typing in the following sense:
\begin{proposition}
\label{prop:renaming-preserve-typing}
If $t:T$ is well-typed then $t[a\leftarrow b]:T[a:=b]$ and
$t[i\leftarrow j]:T$. 
\end{proposition}
\begin{proof}
By induction on the typing derivation for $t:A$. 
\end{proof}

We use the shorthand $(t:T)[a\leftarrow b]$ for $t[a\leftarrow b]:T[a:=b]$.
We define the renaming of a variable in an annotated sequent pointwise on its roots: 
\begin{definition}
\label{def:renaming structure}
  Let $F = t_1:T_1, \dots, t_n:T_n$ be a typed forest.  Define
\[ F[a \leftarrow b] := (t_1:T_1)[a:=b], \dots, (t_n : T_n):[a:=b] \]
\noindent and 
\[ F[i \leftarrow j] := (t_1:T_1)[i:=j], \dots, (t_n : T_n):[i:=j]. \]
\end{definition}





In the process of eliminating admissible structural rules, we must
rename many eigenvariables and indices present in a subproof.  We give
now some notation for such a compound renaming:

\begin{definition}
\label{def:renamefunvars}
Let $V=v_1, \dots, v_n$ and  $x_1, \dots, x_n$, be 
two sequences of variable symbols.  Then define 
\[ \tau^{v_1, \dots v_n}_{x_1, \dots x_n}t := t[v_1 \leftarrow x_1]\dots[v_n \leftarrow x_n] \]
\end{definition}

\begin{definition}
\label{def:renamefuninds}
Let $V=i_1, \dots, i_n$ and  $j_1, \dots, j_n$, be 
two sequences of indices.  Then define 
\[ \tau^{i_1, \dots i_n}_{j_1, \dots j_n}t := t[i_1 \leftarrow j_1]\dots[i_n \leftarrow j_n] \]
\end{definition}

\subsection{Substitution}
\label{sec:substitution}
Suppose that $F$ is a typed forest containing a cut $\alpha[a].t
\bowtie (\varepsilon[M].s)$.  The intuititive explanation of this term
is a pending communication: at some point $\varepsilon[M]$ should
communicate its witnessing term, $M$, to $\alpha[a]$.  This is 
what happens, on the level of annotations, during a single step
of cut-reduction.  To carry out this operation, we must substitute
a first-order term $M$ for an $\alpha$-free variable $a$ in an
annotated sequent.  We define that operation now.

\begin{definition}
\label{def:first-order-substitution}
  We define an operation $[a := M]$ (substitute $M$ for $a$) on
  typed forests $F$ such that $a \notin \bound F$.  

  On witnessesing terms, of the form $\varepsilon[M].t$, the substitution
  applies inside the instantiating first-order term $M$ and in the remaining 
subterm $t$:
\[ (\varepsilon[N].t)[a:=M] = \varepsilon[N[a:=M]](t[a:=M])\]
Substitution is pushed past all the other term constructors, as
follows:
  \begin{align*}
    S[a:=M] &= S\\
    (\alpha[d].t)[a:=M] &= 
    \alpha[d].(t[a:=M]) \\
    (t_1 + \dots + t_n)[a:=M] &=  (t_1[a:=M] +  \dots + t_n[a:=M]) \\
    (t \bowtie s)[a:=M] &=  t[a:=M] \bowtie s[a:=M]\\
  \end{align*}
Finally, $F[a:=M]$ is defined as the pointwise substitution of $M$ for
$a$ in each term of $F$.
\end{definition}

\noindent By induction on the structure of typing derivations, we obtain:

\begin{proposition}
\label{prop:substitution-with-types}
  If $t$ can be assigned type $A$, then $t[a:=M]$ can be
  assigned type $A[a:=M]$.
\end{proposition}

Having defined substitution, we can formally define the operation of
communicating a witness across a cut, which will be one of our cut-reduction 
operations on Herbrand nets:
\begin{definition}
\label{def:communication-reduction}
  Let 
  \[ G = F, \alpha[a].t\bowtie\finset{\varepsilon[M].s}:
  \forall x.A \bowtie \exists x. \bar{A} \] be an annotated sequent. 
The
  \textsc{Comm} reduct of $G$ is 
  \[F[a:=M], (t \bowtie s)[a:=M] \ : A[x:=M] \bowtie \bar{A}[x:=M]\]
\end{definition}

\section{Cut-free completeness of $\LK_H$}
Consider the ``forgetful projection'' of $\LK_H$, where we simply
delete the term annotations.  This yields a standard sequent system
which is a subsystem of the calculus given in Figure~\ref{fig:midseq}.
Thus, to prove cut-free completeness of $\LK_H$, we need only show the
rules weakening and (general) contraction admissible. 
\label{sec:cut-free-compl}
By cut-free completeness of $\LK_H$, we mean the following:

\begin{theorem}
\label{thm:cut-free-completeness}
  Fix a signature $\Sigma$, containing at least one constant symbol.
  For every closed prenex formula $A$ over that signature, valid in classical
  predicate logic, there is an expansion tree $t$, such that $\LK_H \turnstile t:A$.
\end{theorem}

\begin{remark}
  The requirement that our signature contains a constant is related to
  the usual assumption in classical predicate logic that domains are
  non-empty: without it, weakening is not admissible below the midsequent.  
\end{remark}

The following demonstration of invertibility
will be essential: it is precisely the invertibility of the universal
rule which allows admissibility of contraction:

\begin{lemma}
\label{lem:invert-alpha}
The rule
\[
\begin{prooftree}
 F,\ t: A[x:=a]
\justifies
 F,\ \alpha [a].t : \forall x.A
\end{prooftree}
\]
is invertible -- that is, $ \LK_H \turnstile F, \ \alpha [a].t : \forall
x.A$ if and only if ${\LK_H \turnstile F, \ t: A[x:=a]}$.
\end{lemma}
\begin{proof}
  By induction on proof height.  Since its type contains instances of
  quantifiers, $F , \alpha [a].t : \forall x.A$ cannot be a conclusion of
  the tautology rule.  Suppose that the inversion
  holds for all proofs of height~$<n$, and let $\Phi$ be a proof of
  height $n$ of $F, \alpha [a].t : \forall x.A$.  We proceed by a case
  analysis on the last rule $\rho$ of $\Phi$.

  If $\Phi$ has the form
\[
\begin{prooftree}
  \using \Phi'
\proofdotseparation=1.2ex 
\leadsto
\[
  F', \ \alpha [a].t : \forall x.A
\justifies
  F, \ \alpha [a].t : \forall x.A
\using \rho
\]
\end{prooftree}
\]
then we may apply the induction hypothesis to $\Phi'$, which has
height $<n$, to obtain a proof of ${F', \ t : A[x:=a]}$, to which we may
then apply $\rho$.  

Otherwise, $\alpha [a].t : \forall x.A$ is the principal
formula of $\rho$, and $\Phi$ has the form
\[
\begin{prooftree}
  \using \Phi'
\proofdotseparation=1.2ex 
\leadsto
\[
  F, \  t : A[x:=a]
\justifies
  F, \ \alpha [a].t : \forall x.A
\using \forall
\]
\end{prooftree}
\]
and then $\Phi'$ is the desired proof.  
\end{proof}

\newcommand{\meld}{\mathsf{M}}

We show now that weakening and contraction are admissible in $\LK_H$.
We explain briefly what this means in the presence of annotations: let
$\rho$ be an instance of an ordinary sequent rule with premise
$\Gamma$ and conclusion $\Gamma'$.  Then the rule-instance is
admissible in $\LK_H$ if, given a proof whose conclusion $F$ has type
$\Gamma$, there exists a proof of an annotated sequent $G$ with type
$\Delta$.  In fact, for both weakening and contraction admissibility
we prove stronger results, in the sense that there is a close relation
between the annotations in premise and conclusion of the admissible
rule.

 To prove weakening
admissible, we must ensure that the formula introduced by weakening
does not contain any free occurences of eigenvariables; otherwise we
will violate strictness.  This does not, of course, impact
completeness, since we may always rename bound variables before
weakening.

\begin{lemma}
\label{lem:weakening-admissible}
If $\LK_H \turnstile F$, , and $\fv(A) \cap \bound(F)=\emptyset$, then
there is an expansion tree $t_A$ such that 
\[\LK_H \turnstile F, \ t_A:A.\]
\end{lemma}
\begin{proof}
  By induction on the structure of $A$, and on the length of a
  proof $\Phi$ of $F$. First suppose that $A$ is a quantifier-free
  formula $P$.  If $\Phi$ is an instance of the tautology rule
  labelled with $i$, then $F, (i):P$ is also an instance of the
  conclusion of the tautology rule.  By induction on the length of a
  proof of $F$ we may now show that if $\LK_H \turnstile F$, then
  $\LK_H \turnstile F, (i):P$ where $i$ is the index of a tautology in
  the derivation of $F$.
  
  We now show weakening admissible for general $A$, by induction on
  the rank of $A$.  For an induction hypothesis, suppose that all
  formulae of rank $<rk(A)$ admit weakening.  Now suppose that $A =
  \forall x.B$.  By the induction hypothesis, whenever we have a
  proof of $F$, we have a proof of $F, t_B:B[x:=z]$, for
  $z$ a fresh variable, i.e. not appearing in $\free(F)$ or
  $\bound(F)$.  Apply $\forallR$ to obtain a proof of $F,
  \alpha[z].t_B:A$.  

  Finally, suppose $A = \exists x.B$.  Let $*$ be a constant in
  $\Sigma$.  Then if we have a proof of $F$, we have a proof of $F,
  t_B:B[x:=*]$, from which we derive a proof of $F,
  (\varepsilon[*].t_B):A$ by an application of $\existsR$. 
\end{proof}

To define the contraction of two expansion treess $t:A$ and $s:A$, we
must find a ``merge'' $\meld(t,s)$ of the two expansion trees, and a
sequence of variable renamings to be made in the context.

\begin{proposition}[Admissible contraction]
\label{prop:admissible-contraction}
Given any annotated sequent $F, \ t: A,\  s:A$, there is 
a expansion tree $\meld(t,s): A$ and two sequences $(x_i), (y_i)$ 
of variables such that  
\[
\begin{prooftree}
 F, \ t: A,\  s:A
\justifies
 \tau^{(x_i)}_{(y_i)}(F), \ \meld(t,s): A
\end{prooftree}
\]
is admissible in $\LK_H$; that is, $\LK_H \turnstile F, \ t: A,\  s:A$
implies  $\LK_H \turnstile \tau^{(x_i)}_{(y_i)}(F), \ \meld(t,s): A$
\end{proposition}

\begin{proof}
  We proceed by induction on the rank of the formula $A$.  Suppose
  first that the rank of $A$ is zero.  Then $A$ is a QFF $P$, and the
  result trivially holds by an application of $\CP$, with empty sequences of 
variables.

Suppose now that the lemma holds for all $B$ of rank $n-1$, 
and let $A$ have rank $n$.  If $A = \exists x.B$,
then we have contraction on $A$ by the rule $\Cexists$: $\meld(t,s) = t+s$
and again the two sequences of variables are empty.  The
interesting case is where $A = \forall x.B$.  Suppose we have a  proof of
$\turnstile F, \ \alpha[a].t: \forall x.A,
\alpha[b].s:\forall y.A$.  Apply the invertibility of $\forall$ twice
to obtain a proof of

\[ F, \ t: A[x:=a], s:A[x:=b] .\]

If $e$ is a fresh free variable, let ${t' = t[a\leftarrow
  d].[b\leftarrow e]}$ and ${s' = s[a\leftarrow e].[b\leftarrow e]}$.
Making renaming substitution $[a\leftarrow e],[b\leftarrow e]$ inside
the proof we obtain a proof of

\[ F[a\leftarrow e][b\leftarrow e], \ t': A[x\leftarrow e], s':A[x \leftarrow e] .\]

Now apply the induction hypothesis to obtain a term $\meld(t',s'):
A[x:=e]$ and sequences $(x_i), (y_i)$ such that

\[\LK_H \turnstile \tau(F[a\leftarrow e].[b\leftarrow e]), \ \ \meld(t',s'): A[x:=e] .\]

An application of $\forall$ yields a  proof of
\[ \tau(F[a\leftarrow e].[b\leftarrow e]), \ \ \alpha[e].\meld(t',s'): \forall x.A .\]
\noindent as required, with $\meld(t,s) =\alpha[e].\meld(t',s')$ and
sequences $(a, b, x_i), (e, e, y_i)$ .
\end{proof}

This completes the proof of completeness for $\LK_H$.
We can view the cut-free completeness of this calculus as an
alternative strong statement of Herbrand's theorem, since each cut-free
$\LK_H$-proof gives rise to an Herbrand proof:

\begin{proposition}
\label{prop:nets-into-herb-proofs}
Let $A$ be a formula in prenex normal form over a signature $\Sigma$
containing at least on constant symbol.  There is a expansion tree $t$ such that
$\LK_H \turnstile t:A$, if and only if $A$ has an Herbrand proof.
\end{proposition}
\begin{proof}
  One direction is just cut-free completeness: if there is an Herbrand
  proof of $A$, then $A$ is provable, and so there is a expansion tree $t$ such
  that $\LK_H \turnstile t:A$.
  For the other direction, we must extract, from a derivation of $t:A$
  in $\LK_H$, an Herbrand proof of $A$. 

  Let $\Phi$ be a strict proof of $t:A$ in $\LK_H$.  Associate to each
  instance of $\existsR$ in $\Phi$ a distinct variable not occuring in
  $t:A$, and decorate the corresponding witness with that variable -- that is,
we replace each occurence of the $\existsR$ rule with the rule:
    \[ \begin{prooftree}
     F, \ t : A[x:= M_u] \justifies 
     F, \ \finset{\varepsilon_u [M_u].t}: \exists x.A
    \using \existsR_u
  \end{prooftree} \] 
with a different $u$ for each occurence.

 \newcommand{\Deep}{\mathrm{Deep}}
  Given a typed expansion tree $s:B$, labelled as above, we 
  extract an expansion $\Deep(s:B)$ of $B$ as follows:
  \begin{align*}
    \Deep(S:P) &= P \\
    \Deep(\alpha[a].t:\forall x.A) &= \forall a.\Deep(t:A[x:=a])\\
    \Deep((w_1 + \dots + w_n):\exists x.A) &= \Deep(w_1:[\exists x.A])\lor \dots \lor \Deep(w_n:[\exists x.A])\\
    \Deep(\varepsilon_u[M_u].t:[\exists x.A]) &= \exists u.(\Deep(t:A[x:=M_u])
    \end{align*}
    \noindent Let $V$ be the set of bound variables in $\Deep(t:A)$.
    Each member of $V$ is either a label of an instance of $\existsR$
    in $\Phi$ or the eigenvariable of an instance of $\forallR$ in
    $\Phi$.  Since $\Phi$ does not branch (it has no cuts), it imposes
    linear order $x_1, x_2, \dots x_n$ on $V$.  Let $A^*$ be the
    matrix of $\Deep(t:A)$: then $Q_1x_1\dots Q_nx_nA^*$ is a
    prenexification of $\Deep(t:A)$. Let $\sigma_\Phi$ be the sequence
    of first-order terms (eigenvariables and witnessing terms) induced
    by this ordering; then $\sigma_\Phi$ is a witnessing substitution
    for $Q_1x_1\dots Q_nx_nA^*$.  Thus $A$ has an Herbrand proof.
\end{proof}






\section{Herbrand nets}
\label{sec:herbrand-nets}
The definition of annotated sequent is such that ever conclusion of a
strict $\LK_H$ proof is an annotated sequent; the converse is not
true.  For example, an annotated sequent can contain cuts of the form
$\alpha[a].t\bowtie \varepsilon[M(a)].s$. To formulate a criterion
excluding such terms, we treat annotated sequents as \emph{proof
  structures} (known elsewhere as \emph{pre-proofnets}), in the sense
of Girard~\cite{Gir96ProNetPar}.  The sequent calculus cut is, in
form, a linear logic tensor; the standard techniques of multiplicative
proof-net correctness provide precisely the tools to decide whether
such a tensor can be the last rule of a sequent derivation.  We will
define a correctness criterion singling out, among the annotated
sequents, those arising as the conclusion of an $\LK_H$ derivation,
and we will call these annotated sequents \emph{Herbrand nets}.
Finally, we develop \emph{sequentialization}: any annotated sequent
$F$ which is the endsequent of a derivation contains itself enough
information to reconstruct a derivation $\LK_H \turnstile F$
(although, of course, we might reconstruct a different derivation to
the one we originally used to derive $F$).

\subsection{Prelude: expansion-tree proofs}
We begin by giving necessary and sufficient conditions for a cut-free
annotated sequent to be the conclusion of a proof.  First, let $t:A$
be an annotated sequent: that is to say, $t$ is an expansion tree and
$\bound(t) \cap \fv(A)=\emptyset$. Given a proof $\pi$ of $t:A$, we
have linear ordering of the $\alpha$ and $\varepsilon$ nodes in $t$;
this ordering is compatible with the a notion of \emph{dependency} on
those nodes, where the a subterm depends on its predecessor, and a
subterm $\varepsilon[M].t$ depends on the subterm $\alpha[a].s$ if the
variable $a$ appears free in $M$.  This natural notion of dependency
provides a simple condition for deciding whether a given $t:A$ is
provable in $LK_H$; this condition is a special case of the condition
given by Miller for his expansion-tree-proofs:
\begin{proposition}
  Let $t:A$ be an annotated sequent.  Then $t:A$ is provable in $LK_H$
  if and only if
  \begin{enumerate}
  \item the transitive closure $\vartriangleleft$ of the above-defined
    dependency relation is irreflexive, and
  \item $\mathrm{Deep}(t:A)$, as defined in the proof of
    Proposition~\ref{prop:nets-into-herb-proofs}, is a tautology.
  \end{enumerate}

\end{proposition}
The proof of this proposition is easy; a sequent proof of $t:A$ is
(essentially) a linear order on the non-propositional nodes of $t$
extending $\vartriangleleft$, which exists if and only if
$\vartriangleleft$ is irreflexive.  

It is not hard generalize this condition to one necessary and
sufficient for \emph{any} cut-free annotated sequent, by generalizing
$\mathrm{Deep}$ to work on sequents instead of just typed expansion
trees.  The irreflexivity of dependency on is nodes is not enough,
however, when we introduce cuts (we cannot, by this route, exclude
cuts of the form $\alpha[a].t\bowtie \varepsilon[M(a)].s$).  For this
reason, we turn to the the techniques of proof-net correctness.

\subsection{Typed forests as proof structures}
We have mentioned, already, that we consider proof nets to be forests
with a linking structure.  The forest structure of an annotated
sequent has already been discussed at length, and we move now onto the
linking structure.  This will consist of the dependency discussed
above, plus a linking which generalizes the usual axiom links of proof
nets.  The usual sequent calculus axioms are replaced in $\LK_H$ by
tautology rules.  Similarly, the usual axiom links of proof nets,
linking two dual formulae, are replaced in Herbrand nets by
something more general: the information contained at the leaves of a
typed forest plays the role of generalized axiom links.  This
generalization is two-fold: each ``tautology link'' (each index
appearing in a set at some leaf) may have an arbitrary (finite) number
of conclusions, and (because of contraction) each leaf may be
connected to several such links.  For the purposes of correctness,
these links play a dual role.  Later we will see that they play a part
in the switching criterion (which ensures that the structure of
substitution and cuts can be sequentialized).  In addition, the links
are necessary to check whether the \emph{propositional} information in
a structure is correct: if the disjunction of the formulae arising
from a tautology index is really a tautology:

\begin{definition}
Let $F$ be a typed forest, and let $i$ be a tautology index appearing in 
$F$.  The formula $F_i$ is defined as follows:
\[ F_i = \bigvee \finset{ A \ | \ (S): A \ \text{is a propositional
      node in F}, \ i \in S} \]
\end{definition}

\begin{definition}[Herbrand Structure]
\label{def:herbrand-structure}
Let $F$ be an annotated sequent over a theory $(\Sigma, \univax)$.
$F$ is an \emph{Herbrand structure} if, for each tautology $i$ in $F$,
we have $\univax \vDash F_i$.
\end{definition}

The tree structure of a typed forest $F$ defines a natural directed
graph structure, with vertices given by the nodes of $F$ and edges
directed from child to parent.  The linking structure on a typed
forest is given using \emph{jumps} \cite{Gir96ProNetPar} -- extra directed edges in the
graph. For each tautology index in $F$, we add a vertex, and
we add jumps (directed edges) from a tautology index to the leaves
where it appears.  If the variable $x$ appears free a first-order term
$M$, we will make a jump from each $\varepsilon[M]$ to the alpha node
binding $x$.  This jump indicates that, in a sequent proof of $F$, the
existential rule introducing the $\varepsilon[M]$ must occur before
the universal rule introducing the $\alpha[a]$ Less obviously, we also
need jumps from cuts: if the variable $a$ is free in the type of a
cut, then that cut must occur above the rule binding $a$. We will call this graph with jumps the \emph{dependency graph}
of the forest. 

\newcommand{\dependson}{{\color{red}\curvearrowright}}

\newcommand{\Dep}{\mathrm{Dep}} 

\begin{definition}
\label{def:dependency-graph}
  Let $F$ be a semistructure.  The \emph{dependency
    graph} $\Dep(F)$ of $F$ is a labelled directed graph whose vertices are:

\begin{enumerate}
\item The instances of subterms of $F$, plus
\item one node for each tautology index $i$ in $F$, labelled with the index.
\end{enumerate}

The edges of  $\Dep(F)$ are the edges of $F$ considered as a directed graph, plus the \emph{jumps}:
\begin{itemize}
\item An edge from $\varepsilon[M]:B$ to $\alpha[a]$ whenever $a \in
  \fv(M)$;
\item An edge from $t \bowtie s: A \bowtie \bar{A}$ to $\alpha[a]$ whenever $a \in \fv(A)$
\item An edge from the vertex $i$ to each leaf $S$ of $F$ with $i \in S$.
\end{itemize} 
\end{definition}

We use red curved arrows to represent jumps in the dependency
graph, and red labels for the tautology vertices; the black, straight
arrows and black vertices continue to represent the underlying forest.

\begin{example}
\label{ex:dependency-graph-drinker}
The dependency graph of the annotated drinker's term D is
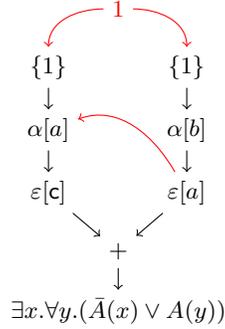
\[
\begin{tikzpicture}
\matrix[row sep=3mm,column sep=-10mm]{
& \node [red](tautology) {$1$}; \\
 \node (left prop) {$\finset{1}$}; &&
 \node (right prop) {$\finset{1}$}; \\
\node (left alpha) {$\alpha[a]$}; &&
 \node (right alpha) {$\alpha[b]$}; \\
\node (left epsilon) {$\varepsilon[\mathsf{c}]$}; &&
 \node (right epsilon) {$\varepsilon[a]$}; \\
 & \node(plus){$+$};\\
 & \node(conclusion){$\exists x. \forall y. (\bar{A}(x) \lor A(y))$} ;\\
};
\draw[->](left prop) to (left alpha);
\draw[->](right prop) to (right alpha);
\draw[->](left alpha) to (left epsilon);
\draw[->](right alpha) to (right epsilon);
\draw[->](left epsilon) to (plus);
\draw[->](right epsilon) to (plus);
\draw[->](plus) to (conclusion);

\draw[red,->](tautology) to [out=180, in = 90] (left prop);
\draw[red,->](tautology) to [out=0, in = 90] (right prop);
\draw[red,->](right epsilon) to [out=120, in = 20] (left alpha);
\end{tikzpicture}
\]
\noindent the dependence of $\varepsilon[b]$ on $\alpha[b]$
is indicated by the upwards-pointing grey arrow.   
\end{example}

\begin{example}
\label{ex:2}
The dependency graph of the annotated sequent derived in example is
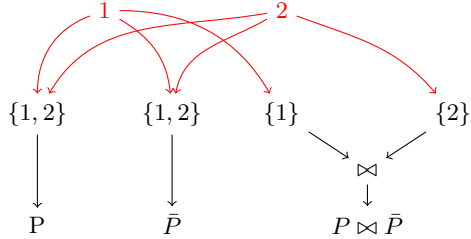
\[
\begin{tikzpicture}
 \matrix[row sep=3mm,column sep=2mm]{ 
 & 
 \node [red](First tautology) {$1$}; & 
                                & & &
 \node [red](Second tautology) {$2$}; & \\[6mm]
 \node (First proposition) {$\finset{1,2}$};  &&
 \node (Second proposition) {$\finset{1,2}$}; &&&
 \node (Left cut) {$\finset{1}$};             &
                                             &
 \node (Right cut){$\finset{2}$};             & \\[0mm]
 &&&& &&\node (Cut) {$\bowtie$};  \\
 \node (First conclusion) {P};                & &
 \node (Second conclusion){$\bar{P}$};        & & &&
 \node (Cut conclusion){$P \bowtie \bar{P}$};                                \\ 
 };
 \draw [red,->] (First tautology)  to [out=200, in=90] (First proposition);
 \draw [red,->] (First tautology) to [out=-30, in=95] (Second proposition);
 \draw [red,->] (First tautology) to [out=0, in=120] (Left cut);
 \draw [red,->] (Second tautology)  to [out=190, in=60] (First proposition);
 \draw [red,->] (Second tautology) to [out=210, in=80] (Second proposition);
 \draw [red,->] (Second tautology) to [out=-20, in=130] (Right cut);
 \draw [<-] (Cut) to (Left cut);
 \draw [<-] (Cut) to (Right cut);
 \draw [->] (Cut) to (Cut conclusion);
 \draw [->] (First proposition) to (First conclusion);
 \draw [->] (Second proposition) to (Second conclusion);
\end{tikzpicture}
\]
\end{example}

\begin{definition}
\label{def:dependency relation}
Let $F$ be a typed forest.  The \emph{dependency relation}
on $F$ is defined to be the relation $\vartriangleright = (\rightarrow
\cup \dependson)^*$, restricted to the nodes of $F$.
\end{definition}


\subsection{Correctness}
\label{sec:correctness}
We will use a variation on the well-established Danos-Regnier ACC
(acyclic-connected) correctness criterion~\cite{DanReg89StrMul}, since
it is well-known and easily stated.  The criterion as given is
exponential (we can decide in exponential time if a given Herbrand
structure is a net), but it is known that correctness for this kind of
proof-net is actually NL-complete \cite{DBLP:conf/csl/NauroisM07}.  Of
course, checking that a given annotated sequent $F$ is an Herbrand
structure is coNP, since we must check that that each $F_i$ is a
tautology.


The crucial notions in ACC correctness are the switching and the
switching graph, which in our setting are defined (for semistructures)
as follows:

\begin{definition}
\label{def:switching}
Let $F$ be a typed forest.
\begin{enumerate}
  \item The \emph{switched nodes} of $F$ are the subterms of the form $\alpha[a].t'$,
    $(t_1 + \cdots + t_n)$, or $S$. All other nodes of $F$ are unswitched.
  \item A \emph{switching} $\sigma$ of $F$
  a choice of, for each switched link $t$ of $F$, exactly one incoming edge
  for $t$ in $\Dep(F)$.  
  \item The \emph{switching graph} $F_\sigma$ of a
  switching $\sigma$ is the undirected graph derived from $\Dep(F)$ by
  deleting, for each switched node $t$, all edges coming into $t$ except
  that chosen by the switching, and then forgetting directedness of edges.
\end{enumerate}
\end{definition}
 
\begin{definition}
\label{def:nets}
A typed forest $F$ is \emph{ACC-correct} (or just ACC), if for each
switching $\sigma$, $F_\sigma$ is connected and acyclic.  A \emph{net}
is an annotated sequent (i.e. with no naked witnesses) which is ACC.
\end{definition}

An Herbrand net is a structure satisfying both kinds of correctness: the
correctness of the tautology links and the switching correctness:

\begin{definition}
\label{def:Herbrand-net}
  An annotated sequent $F$ is an \emph{Herbrand net} if is is an Herbrand
structure and a net.
\end{definition}

\begin{proposition}
  \label{prop:removing-nodes}
\begin{enumerate}
\item If $F = F', \alpha[a].t: \forall x.A$ is ACC then  $G= F', t: A[x:=a]$
  is ACC.
\item $F = F', (w_1 + \cdots w_n):\exists x.A$ is ACC iff
$G=F', w_1:[\exists x.A], \dots 
  w_n:[\exists x.A]$ is ACC.
\item $F = F', S:P$ is ACC iff $F'$ is ACC.
\end{enumerate}
\end{proposition}
\begin{proof}
An easy application of the definition of correctness; passing from $F$ to $G$ removes,
in each case, a switched node which is a root in $F$. This cannot affect either
connectedness or cyclicity of the switching graph.
\end{proof}

\begin{proposition}
\label{LKH-proofs-give-Herbrand-nets}
The conclusion of any $\LK_H$ proof is an Herbrand net.
\end{proposition}
\begin{proof}
By induction on the tree-structure of an $\LK_H$ proof.
\end{proof}

We now justify our claim that Herbrand nets generalize the first-order
expansion-tree proofs of \cite{Mil87ComRep} by showing that, in the
absense of cuts, we may replace our switching notion of correctness by
irreflexivity of $\vartriangleright$.  First, we see that correctness
guarantees that $\vartriangleright$ is a strict partial order among the
alpha and epsilon nodes:

\begin{proposition}
\label{prop:irreflexive-dependency}
  Let $F$ be an an Herbrand net. Then its dependency relation
  $\vartriangleright$ is irreflexive (i.e.  the dependency graph is
  acyclic.)
\end{proposition}
\begin{proof}
  If the dependency graph of $F$ has a cycle, then it has an
  \emph{elementary} cycle $C = t_0, t_1, \dots, t_{n}, t_{n+1}=t_0$, in which
  each vertex is only visited once.  Choose a switching $\sigma$ of
  $F$ in which we choose for each switched node $t_{i+1}$ on that elementary
  cycle the vertex $t_i$ as its switching.  Then the cycle also appears in
  $F_\sigma$, and so $F$ is not ACC. 
\end{proof}

In the cut-free case we derive the converse:
\begin{proposition}
\label{prop:ex-tree-correctness}
  Let $F$ be a cut-free Herbrand structure.  $F$ is an Herbrand net if
  and only if it contains a single tautology index and its dependency
  relation $\vartriangleleft$ is irreflexive.
\end{proposition}

\begin{proof} $F$ is an Herbrand net if and only if each switching
  graph is acyclic and connected.  Assume that $F$ is an Herbrand net:
  by the previous lemma we know that it has an irreflexive dependency
  relation.  Suppose $F$ contains two or more tautology nodes and no
  cuts; then each switching graph of $F$ is disconnected, and $F$ is
  not ACC.

  We now show the converse.  Suppose first that $F$ has a switching
  $\sigma$ for which $F_\sigma$ has a cycle.  Because $F$ contains no
  cuts, each node in $F_\sigma$ has at most one path to a
  propositional node in $F_\sigma$, and thus this cycle cannot pass
  through any propositional or tautology node.  It follows that the
  cycle in $F_\sigma$, restricted to logical nodes, gives a cycle in
  the dependency graph of $F$, and thus a reflexive node for
  $\vartriangleleft.$

  Finally, suppose that each switching graph of $F$ is acyclic, but
  $F$ is not an Herbrand net.  Then, for each switching, each
  non-propositional, nontautology node has a unique path to a
  propositional node in $F$.  This switching graph will be connected
  unless there are at least two tautology nodes in $F$. 
\end{proof}



\subsection{Subnets of Herbrand Nets}
\label{sec:subn-herbr-nets}
To show that the Herbrand nets are precisely the proof structures
arising from $\LK_H$ proofs, we will adapt the notion of \emph{subnets}
of a proof net\cite{212898}; this will also be important in defining
cut-reduction on Herbrand nets. 

\begin{definition}[Substructure]
\label{def:substructure}
  Let $F$ be an ACC forest.  A \emph{substructure} of $F$ is a subset $G$ of
  the nodes of $F$ closed under dependency: that is, if $s \in G$ and
  $t \vartriangleright s$ then $t \in G$.
\end{definition}

Since a substructure $G$ of an ACC forest $F$ is always closed under
subterms, it defines a typed forest.  The following is therefore
well-defined.

\begin{definition}[Subnet]
\label{def:subnet}
Let $G$ be a substructure of the ACC forest $F$.  $G$ is a 
\emph{subnet} of $F$ if it is ACC.
\end{definition}

\begin{figure}
\noindent\hrulefill
\begin{center}
\begin{tikzpicture}
\matrix[row sep=3mm,column sep=-10mm]{
 \node (left prop) {$\finset{1}$}; &&
 \node (right prop) {$\finset{1}$}; \\
\node (left alpha) {$\alpha[a]$}; &&
 \node (right alpha) {$\alpha[b]$}; \\
\node (left epsilon) {$\varepsilon[\mathsf{c}]$}; &&
 \node (right epsilon) {$\varepsilon[a]$}; \\
 & \node(plus){$+$};\\
 & \node(conclusion){$\exists x. \forall y. (\bar{A}(x) \lor A(y))$} ;\\
};
\draw[->](left prop) to (left alpha);
\draw[->](right prop) to (right alpha);
\draw[->](left alpha) to (left epsilon);
\draw[->](right alpha) to (right epsilon);
\draw[->](left epsilon) to (plus);
\draw[->](right epsilon) to (plus);
\draw[->](plus) to (conclusion);

\begin{pgfonlayer}{background}
\node [fill=black!20,fit=(left epsilon) (right epsilon) (left prop)] {};
\end{pgfonlayer}
\end{tikzpicture} \quad
\begin{tikzpicture}
\matrix[row sep=3mm,column sep=-10mm]{
 \node (left prop) {$\finset{1}$}; &&
 \node (right prop) {$\finset{1}$}; \\
\node (left alpha) {$\alpha[a]$}; &&
 \node (right alpha) {$\alpha[b]$}; \\
\node (left epsilon) {$\varepsilon[\mathsf{c}]$}; &&
 \node (right epsilon) {$\varepsilon[a]$}; \\
 & \node(plus){$+$};\\
 & \node(conclusion){$\exists x. \forall y. (\bar{A}(x) \lor A(y))$} ;\\
};
\draw[->](left prop) to (left alpha);
\draw[->](right prop) to (right alpha);
\draw[->](left alpha) to (left epsilon);
\draw[->](right alpha) to (right epsilon);
\draw[->](left epsilon) to (plus);
\draw[->](right epsilon) to (plus);
\draw[->](plus) to (conclusion);

\begin{pgfonlayer}{background}
\node [fill=black!20,fit=(left alpha)  (right prop)] {};
\node [fill=black!20,fit=(right epsilon) (right prop)] {};
\end{pgfonlayer}
\end{tikzpicture} \quad
\begin{tikzpicture}
\matrix[row sep=3mm,column sep=-10mm]{
 \node (left prop) {$\finset{1}$}; &&
 \node (right prop) {$\finset{1}$}; \\
\node (left alpha) {$\alpha[a]$}; &&
 \node (right alpha) {$\alpha[b]$}; \\
\node (left epsilon) {$\varepsilon[\mathsf{c}]$}; &&
 \node (right epsilon) {$\varepsilon[a]$}; \\
 & \node(plus){$+$};\\
 & \node(conclusion){$\exists x. \forall y. (\bar{A}(x) \lor A(y))$} ;\\
};
\draw[->](left prop) to (left alpha);
\draw[->](right prop) to (right alpha);
\draw[->](left alpha) to (left epsilon);
\draw[->](right alpha) to (right epsilon);
\draw[->](left epsilon) to (plus);
\draw[->](right epsilon) to (plus);
\draw[->](plus) to (conclusion);

\begin{pgfonlayer}{background}
\node [fill=black!20,fit=(right epsilon)  (right prop)] {};
\end{pgfonlayer}
\end{tikzpicture}

\noindent\hrulefill
\end{center}
\caption{Three subnets of the drinker's term}
\label{fig:subnetsdrinker}
\end{figure}
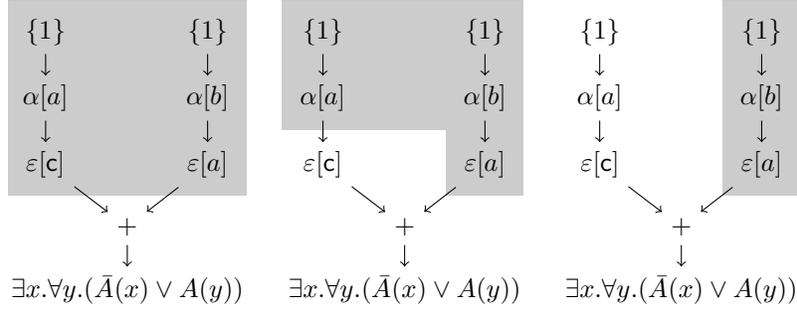

Notice that we do not require that a subnet of an Herbrand net is an
Herbrand net: we do not even require that it is a net.
Figure~\ref{fig:subnetsdrinker} shows three subnets of the drinker's
term, none of which are nets.  As another example, consider the
following immediate consequence of the definition of subnet

\begin{proposition}
  Let $F$ be an $ACC$ forest, and $S$ a leaf of $F$.  Then the subset
$\finset{S}$ of the nodes of $F$ is a subnet of $F$.
\end{proposition}

There is a strong connection between subnets of an Herbrand net and
subproofs of its sequentializations, but developing this is beyond the
scope of the current paper: we will treat the notion of subnet as
a tool to define sequentialization and cut-elimination.

Apart from the oddity that a subnet of an Herbrand net is not
necessarily an Herbrand net, the standard properties of subnets follow
much as for $\MLLminus$ nets.  We state these properties here: the more
technical proofs are summarized in the appendices.

\begin{proposition}
\label{prop:intersection-union-subnets}
Let $G_1$ and  $G_2$ be  subnets of an ACC forest.  
\begin{enumerate}
\item $G_1 \cap G_2$ is a subnet of $F$ if and only if it is nonempty. 
\item If $G_1 \cap G_2$ is nonempty, then $G_1 \cup G_2$ is a subnet.
\end{enumerate}
\end{proposition}
\begin{proof}
\begin{enumerate}
\item Suppose $G = G_1 \cap G_2$ to be nonempty but not a subnet of $F$.  It
is clearly a substructure, so to fail to be a subnet there
must be a switching $\sigma$ for which $G_\sigma$  is disconnected.  But then
either ${G_1}_\sigma$ or ${G_2}_\sigma$ must be disconnected.
\item Now suppose that $G_1 \cap G_2$ is nonempty, but that $G = G_1 \cup
  G_2$ is not a subnet of $F$.  Again, there must be a switching $\sigma$
for which $G_\sigma$ is disconnected.  But since $G = G_1 \cap G_2$ is nonempty,
there is a node $t$ in $G_\sigma$ present in both ${G_1}_\sigma$ and ${G_2}_\sigma$,
and thus connected to each node of $G_\sigma$.
\end{enumerate} 
\end{proof}

\begin{definition}\label{def:empire-kingdom}
  Let $F$ be an ACC forest, and let $t$ be a node in $F$.  The
  \emph{empire} of $t$ in $F$ is the largest subnet of $F$ having $t$
  as a root.  The kingdom of $t$ in $F$ is the smallest subnet having $t$ as
  a root.  
\end{definition}

The kingdom of a node, if it exists, has a particular structure:
\begin{proposition}
  \label{prop:shape-of-the-kingdom}
  Let $t$ be a node of an ACC forest $F$, and let $G, t:A$ be its kingdom.  Then the roots of
$G$ are either witnesses and cuts.
\end{proposition}
\begin{proof}
  By Proposition~\ref{prop:removing-nodes}, if $s:B$ a root of $G$ has
  any other form, we can find a smaller subnet with $t$ as a
  root.  
\end{proof}

By Proposition~\ref{prop:intersection-union-subnets}, if the set of
subnets having a node $t$ as a root is nonempty, $t$ has an empire and
a kingdom.

\begin{proposition}
\label{prop:existence-of-empires}
Let $F$ be an ACC forest.
For each node $t$ in $F$, there is a subnet having $t$ as a root.
\end{proposition}
\begin{proof}
  A variation on the proof in \cite{212898} of the same result for
  \MLL, which we relegate to the appendix ( see Section
  \ref{sec:prop-herbr-nets}).  
\end{proof}

\begin{corollary}
  Every node in $F$ has a kingdom and an empire.
\end{corollary}

The following relation will be the key to our sequentialization
and cut-elimination results.

\begin{definition}
  \label{def:kingdom-ordering}
  Let $F$ be an ACC forest.  We define a relation $\ll$ on the nodes
  of $F$ as follows: $t \ll s$ if $t \in k(s)$.
\end{definition}

If $t$ is a node of an Herbrand net $F$, we can think of the nodes $s$
such that $s \ll t$ as the inference steps that must happen before we may
``do'' the step at the root of $t$.  As one might hope, this relation is an
order extending $\vartriangleleft$:

\begin{proposition}
  \label{prop:kingdom-ordering}
  The relation $\ll$ is a partial order on the terms of a
  ACC forest.
\end{proposition}
\begin{proof}
See Section~\ref{sec:prop-herbr-nets}.
\end{proof}


\subsection{Sequentialization}
\label{sec:sequentialization}
We now seek to establish that every Herbrand net arises as the
conclusion of an $\LK_H$ derivation.  The proof that this is the case
will be an induction using the following measures:

\begin{definition}
\label{def:measure}
Let $F$ be an Herbrand net.
\begin{enumerate}
\item The \emph{size} $s(F)$ of $F$ is the number of $\alpha$, $\varepsilon$
and $\bowtie$ nodes in $F$.
\item The \emph{width} $w(t)$ of an expansion node $t = (w_1 + \dots + w_n)$ in $F$ is
  $n$.  The \emph{width} $w(s)$ of a propositional node $s=S$ in $F$ is the
  cardinality of $S$.  
\end{enumerate}
The $w-rank$ $w(F)$ of an Herbrand net $F$ is $\sum_t(w(t)-1)$, where $t$ ranges over all
expansion nodes and propositional nodes of $F$. 
\end{definition}

We show that all nets may be sequentialized by induction on
$s(F)+w(F)$.  Our base case is where $s(F)=0$ (in which case $w(F)$ is
also $0$):

\begin{proposition}
\label{prop:3}
If $F$ is an Herbrand net of size $0$ (i.e. it contains no $\alpha$,
$\varepsilon$ or $\bowtie$ nodes) it is the conclusion of the
tautology rule of $\LK_H$.
\end{proposition}
\begin{proof}
  Since $F$ contains no $\bowtie$ nodes, and is a net,
   it can contain only one tautology index $i$. So $F$ has the form
  $\finset{1}:P_1, \dots, \finset{1}:P_n$, with $\bigvee P_i$ a
  tautology (since $F$ is an Herbrand structure). 
\end{proof}

In case of non-zero measure, we look for a rule of $\LK_H$ whose
conclusion is $F$ and whose premisses are also Herbrand nets.  This
can be seen as a guided form of proof-search.  In some cases
(corresponding to invertible sequent rules) this is easy:
\begin{proposition}
\label{prop:semantic-inversions}
Let $F$ be an Herbrand net.
\begin{enumerate}
\item If $F = F', \alpha[a].t: \forall x.A$, then $G= F', t: A[x:=a]$
  is also an Herbrand net.
\item If $F = F', s_1+s_2:\exists x.A$, then $G=F', s_1:\exists x.A,
  s_2:\exists x.A$ is also an Herbrand net.
\item If $F = F', S_1\cup S_2:P$ then $G=F', S_1:P, S_2:P$ is also an
  Herbrand net.
\end{enumerate}
\end{proposition}
\begin{proof}
  In each case, it is easy to verify that $F$ satisfies the conditions
  for being an Herbrand net if and only if $G$ does. 
\end{proof}

\noindent These three ``inversions'' produce a net $G$ of lower
measure than $F$.  It is less obvious when to apply the non-invertible
rules of $\LK_H$: the existential rule and the cut-rule.  For example,
in the annotated sequent
\[ (\varepsilon[\mathsf{c}].\alpha[a].\finset{1}): \exists x.\forall
y.A(y), (\varepsilon[a].\finset{1}): \exists z.A(z)\] we cannot
instantiate the rightmost existential until both other quantifier
rules have been instantiated.  For cuts the situation is even more
complicated.  In an annotated sequent of the form
\[ F, \alpha[a].t_1 \bowtie t_2 :A\bowtie \bar{A}], \varepsilon[a].t_3
\bowtie \alpha[b].\varepsilon[a].t_4:B\bowtie  \bar{B}\] the right-hand
cut cannot be decomposed before the left-hand cut, as the
eigenvariable $a$ is used in both sides of that cut.  Even where we do
know which cut to decompose, we must find the correct splitting of the
context.

What this means is that some roots of an Herbrand net can be
decomposed, yielding one or two smaller Herbrand nets, others can not.
We will call the roots of a stucture which admit immediate
decomposition its \emph{gates}:

\begin{definition}
\label{def:gates}
  Let $F$ be an Herbrand net, and let $t:A$ be a root of $F$.  The $t$
  is a \emph{gate} of $F$ if and only if there is a rule instance of
  $\LK_H$, with $F$ as conclusion, with $t:A$ as the active succedent,
  and with premisses that are also Herbrand nets.
\end{definition}

The main work of the rest of this section will be to show that each
Herbrand net has a gate.  We will use the notions of kingdom, empire,
and the relation $\ll$, defined in the previous section.

\begin{proposition}
\label{prop:gate}
Let $F$ be an Herbrand net
\begin{enumerate}
\item \label{itemone} Every root of $F$ of the form $\alpha[a].t$, $\finset{w_1, \dots w_n}$ or a non-singleton set $S$
  is a gate.
\item a root $t$ of the form $s_1 \bowtie s_2$ is a gate if and
  only if it is $\ll$-maximal.
\item a root $t:T$ of the form $(\varepsilon[M].s): \exists x.A$
  is a gate if and only if  $\varepsilon[M]:\witnesstype{ x.A}$ is
  $\ll$-maximal in $F, \varepsilon[M]:\witnesstype{x.A}$.
\end{enumerate}
\end{proposition}
We can immediately see that (a) holds, by Proposition
\ref{prop:semantic-inversions}. Before proving parts (b) and (c), let
us observe that this characterization of gates is enough to show that
every net of nonzero size has a gate:

\begin{proposition}
\label{prop:nets-have-gates}
Let $F$ be an Herbrand net.  Either $F$ is the conclusion of the
tautology rule, or it has a gate.
\end{proposition}

\begin{proof}
  If $F$ has size zero and width zero, $F$ is a conclusion of the
  tautology rule.  Now assume that $F$ has nontrivial
  size/width; by Lemma~\ref{lem:ll-is-an-order}, $\ll$ is a partial
  order on the nodes of $F$, so $F$ has at least one $\ll$-maximal
  node $t$: this node is also, by definition, $\leftarrow$-minimal,
  and so $t$ (with appropriate type) is a root of $F$.  If $t$ is a
  gate, we are done.  Suppose that $t$ is not a gate: then it is of
  the form $\finset{i}$ or $(\varepsilon[M].t)$.  Suppose the former:
  since $F = G, \finset{i}:P$ has nonzero size, so does $G$, and $G$
  is a net: thus $G$ has a gate $t:A$.  This is also a gate of $F$,
  since $t \notin k(\finset{i})$.
  
  Finally, suppose that all $\ll$-maximal nodes of $F$ are of the form
  $(\varepsilon[M_i].s_i)$, for $1\leq i \leq n$; so
\[F = G, \ \ 
  (\varepsilon[M_1].s_1):\exists x_1.A_1,
  \dots,(\varepsilon[M_n].s_n):\exists x_n.A_n \]
 The typed ACC forest 
\[F' = G, \ \ 
  \varepsilon[M_1].s_1:\witnesstype{x_1.A_1},
  \dots,\varepsilon[M_n].s_n\witnesstype{x_n.A_n} \] has an
  $\ll$-maximal node, and it must be $
  \varepsilon[M_j].s_j:\witnesstype{x_j.A_j}$, for some $j$.  This
  node is also $\ll$-maximal in 
  \[G, \ \ (\varepsilon[M_1].s_1):\exists x_1.A_1,
  \dots,\varepsilon[M_j].s_j:\witnesstype{x_j.A_j}, \dots,
  (\varepsilon[M_n].s_n):\exists x_n.A_n, \] (where we have placed a
  $+$ below all the naked witnesses but $\varepsilon[M_j].s_j$) and so
  ($\varepsilon[M_j].s_j):\exists x_j.A_j$ is a gate of $F$. 
\end{proof}

\begin{theorem}[Sequentialization]
\label{thm:sequentialization}
  An annotated sequent $F$ is an Herbrand net if and
  only if it is the endsequent of an $\LK_H$ proof $\pi$. We call
  $\pi$ a sequentialization of $F$.
\end{theorem}

\begin{proof}
  By induction on the $s(F)+w(F)$.  If this measure is zero, $F$ 
is the conclusion of the tautology rule.  Otherwise, $F$ has a gate,
and there is a sequent rule which decomposes $F$ into smaller Herbrand
nets, each of which can be sequentialized by the induction hypothesis.
\end{proof}

The following cases of Proposition~\ref{prop:gate} remain to be proved:

\begin{lemma}[Splitting $\bowtie$]
  \label{lem:splitting-bowtie}
  Let $F = F', t \bowtie s:A\bowtie \bar{A}$ be an ACC forest; then  $t
  \bowtie s$ $\ll$-maximal in $F$ iff there is a  partition $F' =
  F_1, F_2$ such that $F_1, t:A$ and $F_2, s:\bar{A}$ are ACC.
  If, further, $F$ is an Herbrand net, then $F_1, t:A$ and $F_2,
  s:\bar{A}$ are Herbrand nets.
\end{lemma}
\begin{proof}
This is a variation on the standard ``splitting tensor'' theorem for $\MLL$
proof nets: see Section \ref{sec:prop-herbr-nets} for the proof.
\end{proof}

\begin{lemma}
\label{lem:removing-epsilons}
Let $F = G, (\varepsilon[M].t):\exists x.A$ be ACC (resp. an Herbrand
net).  Then $ F' = G, t:A[x:=M]$ is also ACC (resp. an Herbrand net)
if and only if $\varepsilon[M].t:\witnesstype{ x.A}$ is $\ll$-maximal
in $F'' = G, \varepsilon[M].t:\witnesstype{ x.A}$.
\end{lemma}
\begin{proof}

  Suppose $F$ is ACC, and $F'$ is also ACC, and suppose for a
  contradiction that $(\varepsilon[M].t)$ is a member of $k(X)$ for
  some other node $X$ of $F$.  But then the kingdom of $X$ in $F'$ is
  also a subnet of $F$, smaller than $k(X)$, contradicting minimality
  of the kingdom.

  Suppose now that $F'' = G, \varepsilon[M].t:[\exists x.A]$ is ACC
  with $\ll$-maximal node
  $\varepsilon[M].t:\varepsilon[M].t:\witnesstype{ x.A}$. We show that
  $F'$ is ACC.  Since $F'$ is a subgraph of $F''$, all its switching
  graphs are acyclic: we must show that they are also connected.
  Observe that $\var(M) \subseteq \free(F)$.  For otherwise, there is
  a variable $a$ with $a \in \var(M)$, $a \notin \free(F)$; then there
  is a node of $F$ of the form $\alpha[a].s$, and $(\varepsilon[M].t)
  \in k(\alpha[a].s)$, contradicting the fact that
  $(\varepsilon[M].t)$ is a gate.  Thus the node $\varepsilon[M].t$ is
  connected to each switching graph only by its unique successor in
  the forest structure of $F''$, and so removing it cannot disconnect
  any switching graph.

  Finally, notice that $F$ is an Herbrand structure if and only if
  $F'$ is an Herbrand structure, since $F$ and $F'$ have the same
  leaves.  
\end{proof}



\section{Cut-elimination}
\label{sec:cut-elimination}
The sequentialization theorem immediately gives us a way to access a
very weak cut-elimination theorem for Herbrand nets.  Given an
Herbrand net $F$ with type $\G$, it is the endsequent of some proof in
$\LK_H$.  Since we have shown that cut-free $\LK_H$ is complete, there
is a cut free proof of an annotation of $\G$; that proof gives rise to
an Herbrand net $F'$ which is $\bowtie$-free and has the same type as
$F$.  

In this section we will show a system of reductions (``Minimal reduction'')
such that any Herbrand net may be transformed into a cut-free Herbrand net
using these reductions.  The reductions are inspired by reductions in the sequent
calculus, so we begin with a discussion of how one might prove cut-elimination
in $\LK_H$.

\subsection{The basic cut-reduction steps}
\label{subsec:cutredlkh} 

Cut-reduction in classical logic, done stepwise, typically consists of
three kinds of operation, applied to subproofs of the proof being
normalized:

\begin{enumerate}
\item Rank-reducing steps, perfomed on \emph{logical cuts},
where both cut-formulae are the principal formula in the proof-tree above
the cut, and where the rules introducing the cut-formulae are logical; for example:
\begin{equation}
\begin{prooftree}
  \[
  \[
  \leadsto
   \turnstile \Gamma_1, A[x:=a] \quad 
  \using \pi_1
  \]
  \justifies
   \turnstile \Gamma_1, \forall x.A 
  \using \forallR\] 
  \qquad 
  \[
  \[
  \leadsto
  \turnstile  \Gamma_2, \bar{A}[x:=M] \quad
  \using \pi_2
  \]
  \justifies
  \turnstile \Gamma_2, \exists x.\bar{A} 
  \using \existsR\] 
  \using \Cut 
  \justifies
  \turnstile \Gamma_1, \Gamma_2
\end{prooftree}
\qquad \to \qquad
\begin{prooftree}
  \[
  \leadsto
  \turnstile \Gamma_1, A[x:=M] \quad 
  \using \pi_1[x:=M]
  \] 
  \quad 
  \[
  \leadsto
  \turnstile \Gamma_2, \bar{A}[x:=M] \quad
  \using \pi_2
  \] 
  \using \Cut 
  \justifies
  \turnstile \Gamma_1, \Gamma_2
\end{prooftree}
\end{equation}
\item Structural steps, where one of the cut formulae is the principal
  formula in the proof-tree above the cut, and the rule introducing it
  is a structural rule; for example, contraction:
\begin{equation}
\begin{prooftree}
\[\leadsto
  \turnstile \Gamma_1, A
\using \pi_1
\] 
  \qquad 
  \[
\[
\leadsto
   \turnstile\Gamma_2, \bar{A}, \bar{A} \quad
\using \pi_2
\]
  \justifies
   \turnstile\Gamma_2, \bar{A} 
  \using \C \] 
  \using \Cut 
  \justifies
   \turnstile\Gamma_1, \Gamma_2
\end{prooftree}
\qquad \to \qquad
\begin{prooftree}
\[
  \[ 
  \leadsto
  \turnstile\Gamma_1, A 
  \using \pi_1
  \]\qquad  
  \[
  \[
\leadsto
   \turnstile\Gamma_1, A
  \using \pi_1
 \] 
  \qquad 
\[
\leadsto
  \turnstile\Gamma_2, \bar{A}
  \using \pi_2
\]
  \justifies
   \turnstile\Gamma_1,\Gamma_2, \bar{A}
  \using \Cut \] 
  \using \Cut 
  \justifies
  \turnstile\Gamma_1,\Gamma_1,\Gamma_2
\]
\justifies
 \turnstile\Gamma_1,\Gamma_2
\using \C^*
\end{prooftree}
\end{equation}
\item Operations to manipulate the proof tree (rule permutations), in
  order to place a cut in one of the two forms above.
\end{enumerate}

To discover the basic cut-reduction operations of Herbrand nets, we
can examine the action these cut-reduction steps have on the
\emph{annotated} conclusions of $\LK_H$ proofs.

For (c), the rule permutations, this question is easy to answer:
one great advantage of working with (box-free) proof nets is that we
no longer need the third item on this list, and indeed it is easy to
show that two $\LK_H$ proofs differing by a permutation of rules have
the same annotated sequent as a conclusion.

%


Let us consider now the logical cut in $\LK_H$: the annotated
redex is:

\[
\begin{prooftree}
  \[
  F_1, t: A[x:=a]
  \justifies
  F_1, \alpha[a].t: \forall x.A 
  \using \forallR\] 
  \qquad 
  \[
  F_2, s: \bar{A}[x:=M]
  \justifies
  F_2, \finset{\varepsilon[M].s}: \exists x.\bar{A} 
  \using \existsR\] 
  \using \Cut 
  \justifies
  F_1, F_2, \ \alpha[a].t \bowtie
  (\varepsilon[M].s) : \forall x.A \bowtie \exists x. \bar{A}
\end{prooftree}
\]
and the annotated reduct, where the term $M$ has been ``communicated''
across the cut:  
\[
\begin{prooftree}
  F_1[a:=M], \ t[a:=M]: A[x:=M]
  \qquad 
  F_2, \ s: \bar{A}[x:=M]
  \using \Cut 
  \justifies
  F_1[a:=M], \ F_2, \ t[a:=M] \bowtie s : A[x:=M]\bowtie \bar{A}[x:=M]
\end{prooftree}
\]

Since permutations of a sequent proof in $\LK_H$ leave the
annotated conclusion unchanged, we obtain the following result:
\begin{lemma}
  If 
  \[\LK_H \turnstile F, \quad \alpha[a].t \bowtie
  (\varepsilon[M].s): \forall x. A \bowtie \exists x.\bar{A}, \]  then
 \[ \LK_H \turnstile F[a:=M], \quad t[a:=M] \bowtie
  s: (A \bowtie \bar{A})[x:=M] \]
\end{lemma}
That is, the reduct of a logical cut does not
depend on the sequential structure of the proof containing it.

The redex of a cut against contraction has the form

\begin{equation}
\begin{prooftree}
\[\leadsto
  F, \alpha[a].t: \forall x.A
\using \Phi
\] 
  \qquad 
  \[
\[
\leadsto
  F', s_1: \exists x.\bar{A}, s_2: \exists x.\bar{A}
\using \Psi
\]
  \justifies
  F, s_1+s_2: \exists x.\bar{A} 
  \using C \] 
  \using \Cut 
  \justifies
  F, F', \ \alpha[a].t \bowtie
  s_1 + s_2 : \forall x.A \bowtie \exists x. \bar{A}
\end{prooftree}
\label{eq:structredex}
\end{equation}

In order to write down the reduct, maintaining strictness, we must
rename all $\alpha$-bound variables occurring in the duplicated copies
of $\Phi$.  If $\free (F, \alpha[a].t: \forall x.A)= x, y, \dots, z$,
then let $V_0= x_0, y_0, \dots, z_0$, and $V_1 = x_1, y_1, \dots z_1$.
If $I$ is the set of tautology indices in $\Phi$, define $I_0$ and
$I_1$ similarly.  Recalling the renaming functions of
Definitions~\ref{def:renamefunvars}~and~\ref{def:renamefuninds}, let Let $\tau_i(t) =
\tau^V_{V_i}(\tau^I_{I_i}(t))$.  The reduct of~\eqref{eq:structredex}
is
\[
\begin{prooftree}
\[
  \[ 
  \leadsto
 \tau_1(F), \alpha[a_1].\tau_1(t): \forall x.A 
  \using \Phi_1
  \]\qquad  
  \[
  \[
\leadsto
  \tau_0(F), \alpha[a_0].\tau_0(t): \forall x.A
  \using \Phi_0
 \] 
  \qquad 
\[
\leadsto
 F', s_0: \exists x.\bar{A}, s_1: \exists x.\bar{A}
  \using \Psi
\]
  \justifies
  \tau_0(F), F', \ \alpha[a_0]\tau_0(t).\bowtie
  s_0: forall x.A \bowtie \exists x. \bar{A}, s_1: \exists x.\bar{A}
  \using \Cut \] 
  \using \Cut 
  \justifies
    \tau_0(F),  \tau_1(F), F', \  \alpha[a_0]. \tau_0(t) \bowtie
  s_0: \forall x.A \bowtie \exists x. \bar{A}, \alpha[a_1].\tau_1(t) \bowtie
  s_1: \forall x.A \bowtie \exists x. \bar{A}
\]
\justifies
   \meld(\tau_0(F),  \tau_1(F)), F', \  \alpha[a_0]. \tau_0(t) \bowtie
  s_0: \forall x.A \bowtie \exists x. \bar{A}, \alpha[a_1].\tau_1(t) \bowtie
  s_1: \forall x.A \bowtie \exists x. \bar{A}
\using C^*
\end{prooftree}
\]
\noindent where here $\C^*$ denotes many applications of \emph{admissible}
contraction (and $\meld$ is extended to annotated sequents in an appropriate
fashion).

This ``reduction'' is highly problematic in Gentzen's sequent
calculus, as it is not possible to construct a measure which decreases
on its application -- not possible, since reduction using this rule
diverges.  The problem is fixed by Gentzen by moving to a more general
calculus with multicut, for which a rather complicated measure
calculated by following cut-formulae up the sequent-tree can be
calculated.  Using multicuts amounts to imposing a strategy on the
application of the above reduction, for which a measure is given by
looking upwards into the sequent proof to see how many times a cut
formula is duplicated by a contraction.  (In $\LK_H$, this is
precisely the information contained in the width of a term $s$
annotating a contracted formula $\exists x.A$, so giving a direct
cut-elimination result for $\LK_H$ is easier).

Unlike the communication reduction above, it is not immediately clear
how to apply this ``duplication'' reduction to an Herbrand net without
an accompanying derivation.  The results of the previous section give us
a way to proceed if the cut to be reduced is a gate: in that case,  
we can write the net as $F_1, F_2, \ \alpha[a].t \bowtie
  s_1 + s_2 : \forall x.A \bowtie \exists x. \bar{A}$
  \noindent where $F_1, \alpha[a].t:\forall x.A$ and $F_2, s_1+s_2:
  \exists x.\bar A$ are also Herbrand nets.  It is then easy to see that
\begin{equation} G' = \tau_0(F_1), \tau_1(F_1), F_2, \ \alpha[a_0].\tau_0(t) \bowtie s_1:
  \forall x.A \bowtie \exists x. \bar{A}, \ \alpha[a_1].\tau_1(t) \bowtie s_2:
  \forall x.A \bowtie \exists x. \bar{A} \label{eq:2}
\end{equation}
is also an Herbrand net. For each term in $F_1$, we have created two
copies, with the same type.  For each member of $u$ of $F_1$ which is
a expansion tree (i.e. not a cut), we can apply admissible contraction
to $\tau_0(u), \tau_1(u)$, obtaining a single term of the same type as
$u$.  Thus we may recover a net of the same type as $G$ by
applications of admissible contraction, as in the sequent calculus.

The annotated conclusion of the derivation
resulting from applying the $\LK_H$ reduction is dependent on the
derivation, and not just its annotated conclusion; unlike the
communication reduction, this reduction does not commute with rule
permutations.  For any Herbrand net
\[G = F, \ \alpha[a].t \bowtie_X s_1 + s_2 : \forall x.A \bowtie \exists
x. \bar{A}\]
\noindent there are many different sequentializations, and so many
different possible results of reducing the cut labelled $X$.  Consider,
for example, the following proof, and the cut labelled $X$ within it.  
\begin{equation}
\begin{prooftree}
\[
\[\leadsto
  F, \alpha[a].t: \forall x.A, \alpha[b].t': \forall x.B
\using \pi_1
\] 
  \qquad 
  \[
\[
\leadsto
  F', s_1: \exists x.\bar{A}, s_2: \exists x.\bar{A}
\using \pi_2
\]
  \justifies
  F, s_1+s_2: \exists x.\bar{A} 
  \using C \] 
  \using \Cut_X
  \justifies
  F, F', \ \alpha[b].t': \forall x.B \ \alpha[a].t \bowtie
  s_1 + s_2 : \forall x.A \bowtie_X \exists x. \bar{A}
\] \qquad 
\[
\leadsto 
\using \pi_3
F'', s':\exists x.B
\]
\justifies 
F, F', F'', \ \alpha[a].t \bowtie
  s_1 + s_2 : \forall x.A \bowtie_X \exists x. \bar{A}, \
  \alpha[b].t' \bowtie s' : \forall x.B \bowtie \exists x. \bar{B}
\using \Cut
\end{prooftree}
\end{equation}
If we apply one reduction step to the cut labelled $X$, the subproof
$\pi_1$ will be duplicated, but not the subproof $\pi_3$. If, however,
we perform a permutation of the instances of $\Cut$, we obtain a 
derivation with the same annotated conclusion, but different cut-reduction
behaviour:
\begin{equation}
\begin{prooftree}
 \[
    \[
    \leadsto
    F, \alpha[a].t: \forall x.A, \ \alpha[b].t': \forall x.B
    \using \pi_1
    \] 
 \qquad 
    \[
    \leadsto 
    F'', s':\exists x.B
    \using \pi_3
    \]
    \justifies
  F', F'', \ 
  \alpha[b].t' \bowtie s' : \forall x.B \bowtie \exists x. \bar{B}
  \using \Cut
  \]
\qquad 
  \[
       \[
       \leadsto
       F', s_1: \exists x.\bar{A}, s_2: \exists x.\bar{A}
       \using \pi_2
       \]
  \justifies
  F, s_1+s_2: \exists x.\bar{A} 
  \using C 
  \] 
  \using \Cut_X
\justifies 
F, F', F'', \ \alpha[a].t \bowtie
  s_1 + s_2 : \forall x.A \bowtie_X \exists x. \bar{A},
  \ \alpha[b].t' \bowtie s' : \forall x.B \bowtie \exists x. \bar{B}
\end{prooftree}
\end{equation}
Here the reduction of cut $X$ will cause the duplication of both $\pi_1$
and $\pi_3$.

In the following section we will see how to find a
canonical such reduction, derived from the \emph{kingdom} of the cut.

\subsection{Substitution and reduction triples}
\label{sec:reduction-triples}
\newcommand{\fterm}{f_{\mathrm{root}}}
\newcommand{\ftaut}{f_{\mathrm{taut}}}
\newcommand{\gterm}{g_\mathrm{root}}
\newcommand{\gtaut}{g_\mathrm{taut}}
\newcommand{\hterm}{h_\mathrm{root}}
\newcommand{\htaut}{h_\mathrm{taut}}
\newcommand{\eterm}{h_\mathrm{root}}
\newcommand{\etaut}{h_\mathrm{taut}} 

We saw in the previous section how to reduce a cut in an Herbrand net,
in the special case where the cut is a gate of the net.  This
corresponds to the case, in the sequent calculus, of a cut which is
the last inference in a derivation.  In the sequent calculus, this is
easily extended to the general case by reducing the cut within a
sub-proof.  In Herbrand nets this is a harder notion to define: we do
so in this section by means of an auxilliary notion: a
\emph{substitution triple}.

We begin by attempting to define the duplication reduction on an
arbitrary ACC forest; that is, a forest
\[ G = F,\ \alpha[a].t \bowtie s_1 + s_2 : \forall x.A \bowtie
\exists x. \bar{A}\]
\noindent where $F$ possibly contains naked witnesses.  Suppose that
the cut to be reduces is $\ll-maximal$ in $G$; then $G= F_1, F_2 ,\
\alpha[a].t \bowtie s_1 + s_2 : \forall x.A \bowtie \exists
x. \bar{A}$, with $F_1, \alpha[a].t:\forall x.A$ and $F_2, s_1 +
s_2:\exists x.\bar{A}$ also ACC.  As before, we set
\begin{equation} G' = \tau_0(F_1), \tau_1(F_1), F_2, \
  \alpha[a_0].\tau_0(t) \bowtie s_1: \forall x.A \bowtie \exists
  x. \bar{A}, \ \alpha[a_1].\tau_1(t) \bowtie s_2: \forall x.A \bowtie
  \exists x. \bar{A} 
\end{equation}
Where possible, we now try to make the type of $G'$ the same as the
type of $G$: form a new ACC forest $G''$ by applying admissible
contraction to $\tau_0(u), \tau_1(u)$ in $G'$ wherever $u$ is a
expansion tree. Unlike before, $G'$ also contains naked witnesses to
which admissible contraction cannot be applied.  


Let $\fterm$ be the function from the roots of $G$ to the roots of
$G''$, which is the identity on $F_2$, maps all roots of expansion
trees resulting from admissible contraction to their corresponding
roots in $F_1$, and maps $\tau_i(t:T)$ to $t:T$ where $t$ is a cut or
naked witness in $F_1$.

Now suppose that $G$ was a subnet of some Herbrand net $H$.  How can
we define a reduct for $H$, given $G''$?  We seek a net $H'$ with the
same conclusion as $H$, in which the subnet $G$ has been ``replaced''
by $G''$; that $G$ is deleted, and then $G''$ is ``wired in'' in its
place.

$\fterm$ satisfies the following properties:
\begin{itemize}
\item $\fterm$ is surjective;  
\item $\fterm$ is injective on expansion trees
\item $\fterm$ is type-preserving.
\end{itemize}

This function allows us to see which roots of $G''$ should be the
wired where in $H$; that is, which node of $H\setminus G$ should
be the predecessor of a given root of $G''$ in $H'$.

Having replaced $G$ by $G''$ in $H$, there will be a mismatch between
the tautology indices appearing in $H'$; there will be indices in $G$
which were duplicated, but which also occur in $H \setminus G$.
Similarly to $\fterm$, we can define a function $\ftaut$ mapping
tautology indices appearing in $G''$ to tautology indices appearing in
$G$; it maps all tautology indices appearing in $F_2$ to themselves,
and the two copies of a copied index back to the original.  It is
clear that
\[
G''_i \leftrightarrow G_{\ftaut(i)}
\]
where, we recall, $F_i = \bigvee\finset{ A \ | \ S: A \ \text{is a
    leaf in $F$}, \ i \in S}$.  This then allows us to conclude that,
when replacing $G$ by $G''$ in $H$, we do not violate the
Herbrand-structure condition if we replace every leaf $S$ of $H
\setminus G$ by $\bigcup_{i \in S} \ftaut(i)$.  This concludes our
informal description of how to replace $G$ by $G''$ in $H$.

The forest $G''$, plus the functions $\fterm$ and $\ftaut$ form the
first example of what we will call a \emph{substitution triple}:

\begin{definition}
\label{def:reduction-triples}
  Let $F$ be an ACC forest.  A \emph{substitution triple} for $F$ is a
  triple $(F', \fterm, \ftaut)$ such that
\begin{itemize}
\item $F'$ is an ACC forest.
\item $\fterm$ is a function from the roots of $F'$ to
  the roots of $F$ which
  \begin{enumerate}
  \item preserves types on non-cut roots
  \item is surjective
  \item is injective when restricted to the roots of expansion trees
    of $F'$
  \end{enumerate}
\item $\ftaut$ is a function from the tautology indices of $F'$ to the 
  tautology indices of $F$ such that \[F'_i \leftrightarrow F_{\ftaut(i)}.\]
\end{itemize}
\end{definition}

The intended meaning of a substitution triple $(G', \fterm, \ftaut)$
for $G$ is that whenever $G$ is a subnet of an ACC forest $F$, we may
``substitute'' $G'$ for $G$ in $F$.  This has been explained above in
terms of ``wiring''; we now make this notion formal, for the
particularly easy case where every non-cut root of $G$ is a naked
witness: in this case the substitution triple gives a surjective
function from the roots of $G'$ to the roots of $G$, which specifies
to which expansion-node of $F\setminus G$ the non-cut roots of $G'$
should be connected.  This is the content of the following definition
and subsequent lemma: if we have a substitution triple for a subnet of
an $ACC$ forest $F$, it can be extended to a substitution triple for
$F$:

\begin{definition}
\label{def:using-reduction-triples}
Given a subnet $G$ of an ACC forest $F$, all of whose non-cut roots are
witnesses, and a substitution triple $(G', \gterm, \gtaut)$ for $G$, we
construct a triple $(F[G'/G], \fterm, \ftaut)$ (which we prove in the
following lemma, is a substitution triple for $F$).  

The function $\ftaut$ is defined as follows: $\ftaut(i) = \gtaut(i)$
if $i$ is in the domain of $\gtaut$, and $\ftaut(i)=i$ otherwise.

We now define $F[G'/G]$, and the function $\fterm$.  Any roots of $F$
that are also roots of $G$ are either naked witnesses or cuts: let $t$
be such a root: then the members
of $\gterm^{-1}(w)$ are roots of $F[G'/G]$; if $z \in \gterm^{-1}(w)$
then $\fterm(z)=w$. 

If $t$ is a root of $F$ but not a root of $G$, it may still contain
witnesses which are roots of $G$, which we replace by their pre-images
under $\fterm$.  In addition, it will contain tautaology indices which
must be replaced by their images under $\ftaut$.  More formally, if
$t$ is a root in $F$ but not a root of $G$, construct a new term
$R(t)$ as follows:
\begin{align*}
  R(\finset{S}) &= \bigcup_{i \in S} \ftaut^{-1}(i)\\
R(s \bowtie t) &= R(s) \bowtie R(t) \\
R(\alpha [a]. t) & = \alpha [a]. R(t)\\
R(w_1+ \cdots + w_n) & = R(w_1)+ \cdots +  R(w_n)\\
R(\varepsilon[M].t) & =
\begin{cases}
  \varepsilon[M].R(t) & \varepsilon[M].t \notin G \\
  \Sigma(\fterm^{-1}(\varepsilon[M].t)) & \varepsilon[M].t \in G\\
\end{cases}
\end{align*}
$R(t)$ is a root of $F[G'/G]$ and  $\gterm(R(t)) = t$. 
\end{definition}

It remains to show that this defines a substitution triple:

\begin{proposition}
 \label{prop:6}
The triple $(F[G'/G], \fterm, \ftaut)$ is a substitution triple for $F$.
\end{proposition}
\begin{proof}
  By contruction, the functions $\fterm$ and $\ftaut$ have the
  required properties.  The only nontrivial observation is that
  $F[G'/G]$ is an ACC forest.  To see this, consider informally how we
  constructed $F[G'/G]$. Essentially, we do three things:
  \begin{enumerate}
    \item  delete all the nodes of $G$ from $F$, 
    \item wire every witness $w$ in $G'$ to the $+$ which was the parent of
      $\gterm(w)$, and
    \item replace every leaf of the resulting structure with the union
      of the inverse images of its members under $\ftaut$.
\end{enumerate}
Using this intuition, we prove the claim.  Suppose that there was a
switching $\sigma$ of $F[G'/G]$ such that the switching graph of
$F[G'/G]$ contained a cycle $p$.  We will ``project'' that cycle onto
$F$.

The cycle $p$ must intersect the nodes of $G'$, else it would also be
a cycle in $F$.  Let $p'$ be a maximal subpath of $p$ not intersecting
with the nodes of $G'$, and let $X$ and $Y$ be its endpoints.  let
$X'$ and $Y'$ be the nodes in $G'$ to which the subpath $p'$ may be
extended in $p$.  $X'$ is either a witness, and the root of $G'$, or a
tautology index in $G'$.  The same holds for $Y'$.  Via either
$\ftaut$ or $\fterm$, there are corresponding nodes $X''$ and $Y''$ in
$G$, and since $G$ is a subnet, for any switching $\sigma'$ there is a
path from $X''$ to $Y''$ in the switching graph of $F$.  We already
know that there is a path $p'$ from $X$ to $Y$ in $F$.  By choosing a
switching such that $X$ jumps to $X''$ and $Y$ jumps to $Y''$, we find
a cycle in the switching graph of $F$.  
\end{proof}

The substitution triples we are interested in arise from the basic
cut-reduction operations of communication and duplication, closed
under subnets and composition: we will call these triples
\emph{reduction-triples}
\begin{definition}
  \label{def:2}
  A reduction triple for an ACC forest F has one of the following forms: 
  \begin{enumerate}
  \item  If $F = F_1, F_2, \alpha[a].t \bowtie s:\forall x.A\bowtie
    \exists x.\bar{A} $,
and  $F_1, \alpha[a].t:\forall x.A$, $F_2, s:\exists x.\bar{A}$ are also ACC forests, then
\begin{enumerate}
\item (Communication) if $s = \varepsilon[M].s'$, then
  \[(F_1[a:=M], F_2, \ t[a:=M] \bowtie s': A[x:=M] \bowtie
  \bar{A}[x:=M], \quad \fterm, \ftaut)\] is a reduction triple, where
  $\fterm$ and $\ftaut$ are the evident bijections between the
  roots/indices.
\item (Duplication) if $s = s_0 + s_1$, and all roots of $F_1$, $F_2$
  are either cuts or naked witnesses, then $(F', \fterm, \ftaut)$ is
a reduction triple, where
  \[F' = \tau_0(F_1), \tau_1(F_1), F_2,\ \tau_0(t) \bowtie s_0:\forall
  x.A\bowtie \exists x.\bar{A},\ \tau_1(t) \bowtie s_1:\forall
  x.A\bowtie \exists x.\bar{A}\] where $\fterm$ is the identity on roots/indices coming from
  $F_2$, maps $\tau_i(t)$ to $t$, and maps the two new cuts to the cut
  reduced: similarly $\ftaut$.
\end{enumerate}
\item (composition) If $(F', \fterm,\ftaut)$ is a reduction triple for
$F$, and $(F'', \fterm', \ftaut')$ is a reduction triple for $F'$, then
$(F'', \fterm\circ\fterm', \ftaut\circ\ftaut')$ is a reduction triple for
$F$.
\item (reduction in a subnet) If $K$ is a subnet of $F$, with all
  roots of $K$ being cuts or naked witnesses, and $(K, \gterm,
  \gtaut)$ is a reduction triple for $K$, then $(F[K'/K], \fterm,
  \ftaut)$ is a reduction triple for $F$.
  \end{enumerate}
  
\end{definition}

As an example of the above, we will now look at the reduction of a
structural cut in an Herbrand net which is not $\ll$-maximal.  Recall
that, in the sequent calculus this reduction is sensitive to
permutations of inference steps; we would like to define a reduction
directly on the Herbrand net.  To do this, we can look at a subnet in
which it \emph{is} $\ll$-maximal.  Such a net always exists: we can
take the \emph{kingdom} of the cut.  The following is an immediate
consequence of Proposition~\ref{prop:kingdom-ordering}:

\begin{proposition}
\label{prop:gate-of-own-kingdom}
A node $t$ in an ACC forest $F$ is  $\ll$-maximal in $k(t)$.
\end{proposition}

\noindent Since in $k(X)$ the cut $X$ is $\ll$-maximal, and since all
the roots of $k(X)$ are either naked witnesses or cuts, we can apply
the duplication reduction inspired by the sequent calculus, to obtain,
not just a reduct $K'$ of $k(X)$, but also the functions $\fterm$ and
$\ftaut$: a reduction triple for $k(X)$.  In addition, the
noncut-roots of $k(X)$ are all witnesses
(Proposition~\ref{prop:shape-of-the-kingdom}) and so we may apply the
construction in Definition \ref{def:using-reduction-triples}, giving a
reduction triple for the ACC forest $F$.  Since this is a rather
important operation on Herbrand nets, we will take the trouble to
define it directly:

\begin{definition}[The duplication reduction \textsc{Dup}]
\label{def:duplication-reduction}
Let $G = F, \alpha[x].t\bowtie_Y (s_1 + s_2): A\bowtie \bar{A}$ be an
Herbrand net.  Let $X$ be the unique node labeled with $\alpha[x]$ in
$G$, and let $k(X)$ be its kingdom in $G$.  Let $V$ be the variables
bound in $\alpha$ binders in $k(X)$, and $I$ be the tautology nodes in
$k(X)$.  Let the functions $\tau_0$ and $\tau_1$ be renaming functions
as before for the sets $V$ and $I$.  We
define the minimal duplication reduct of the
\begin{align*}
  D_x(i) &= \finset{\tau_0(i), \tau_1(i)};\\
  D_x(S) &= [\bigcup_{i \in S} D_x(i)]\\
  D_x(t \bowtie s) &= \begin{cases}
    D_x(t) \bowtie D_x(s) & t\bowtie s \notin k(x)\\
    \tau_0(t \bowtie s), \ \ \tau_1(t \bowtie s) & t\bowtie s \in k(x)\\
\end{cases}\\
D_x(\alpha [a]. t) & = \alpha [a]. D_x(t)\\
D_x(t_1+ \cdots + t_n) & = (D_x(t_1))+ \cdots +  (D_x(t_n))\\
D_x(\varepsilon[M].t) & =
\begin{cases}
  \varepsilon[M].D_x(t) & \varepsilon[M] \notin k(x)\\
  \tau_0(\varepsilon[M].t) +  \tau_1(\varepsilon[M].t)& \varepsilon[M] \in k(x)\\
\end{cases}
\end{align*}
\noindent (where by ``$t$ is in $k(x)$'' we mean ``the node labelled
with $t$ is in k(x)'').

\noindent Define $D_x(F)$ pointwise on the roots
of $F$.  Then $F$ \textsc{Dup}-reduces to
\[ D_x(F), \quad \alpha[x_0].\tau_0(t) \bowtie_S s_0, \quad   \alpha[x_1].\tau_1(t) \bowtie_S s_1 \]
\end{definition}

\subsection{The principal lemma for partial cut-elimination}
\label{sec:princ-lemma}
In this section we give state and prove the following reduction lemma:
\begin{lemma}
\label{lem:principal-lemma}
  Let $F = G, \ t\bowtie s :A \bowtie  \bar{A}$ be an ACC forest, where all
  cuts appearing in $G$ are of rank 0.  Then $F$ has a reduction triple $(F',
  \fterm, \ftaut)$ such that $F'$ contains no cuts of nonzero rank.
\end{lemma}
This is a generalization of the following, which says that we
can remove a single cut of non-zero rank from a net:
\begin{corollary}
Let $F = G, \ t\bowtie s :A \bowtie \bar{A}$ be an Herbrand net,
and let $G$ contain only cuts of rank 0.  There is an Herbrand net
$F'$, with the same type as $F$, containing only cuts of rank 0.
\end{corollary}

\begin{proof}
  Consider the reduction triple $(F', \fterm, \ftaut)$ provided for
  $F$ by the principal lemma.  The function $\fterm$ is surjective,
  and injective on expansion trees.  Since all the non-cut terms of $F$ are
  expansion trees, there is a 1-1 correpsondence between non-cut terms of $F$
  and non-cut terms of $F'$. Since $\fterm$ preserves types,
  $F'$ is a net and $F$ and $F'$ have the same type.  Finally, if $i$
  is a quantifier index in $F'$, $A^{F'}_i \leftrightarrow
  F_{\ftaut(i)}$, and $F_{\ftaut(i)}$ is a tautology (since $F$ is
  an Herbrand net).  So $F'$ is also an Herband net. 
\end{proof}

The proof of the reduction lemma is
strikingly close to Gentzen's original demonstration of
cut-elimination for the classical sequent calculus, with two
adjustments.  These adjustments both arise from the lack of tree
structure in a proof.  First, we can no longer speak of the
``topmost'' cut in a proof; instead, we eliminate cuts which are
potentially topmost; maximal cuts with respect to the order $\ll$.
Second, we cannot use any notion of height as an induction measure:
instead we use a more natural measure of the complexity of a cut: the
number of witnesses taking place in it (its ``width'').

\begin{proof}(Of the reduction lemma)
  Our proof proceeds by an induction over three measures, ordered
  lexicographically: the first is the size of the ACC forest, meaning the
  number of nodes it has.  The second is the rank of the cut
  apppearing in the ACC forest.  The final measure is the ``width'' of
  the cut: if the cut-term decorating the cut is $\alpha[a]\bowtie s$,
  then the width of the cut is the width of $s$ -- otherwise the 
  width of the cut is 0.

  Our base case is where all cuts are of rank 0; there is no work
  to be done, and we can set the $F=F'$ and both functions $\fterm$
  and $\ftaut$ to be the identity.

  Suppose now that there is a nonzero cut of rank $n$, but that $F$ is
  not equal to the kingdom of that cut.  Then we can find a smaller
  ACC forest $K$ containing the cut: the kingdom of the cut.  By the
  induction hypothesis, we obtain a nonzero-cut-free reduction triple
  $(K', \fterm, ftaut)$ for $K$, and hence by Proposition~\ref{prop:6} a
  nonzero-cut-free reduction triple for $F$.
   
  Now suppose that the nonzero cut is a gate of $F$, and that $F$ is
  the kingdom of that cut.  Then we may write $F$ as
  \[ F_1, \alpha[a].t \bowtie s:(\forall x.A \ | \ \exists x.\bar{A}),
  F_2\] where $F_1, \alpha[a].t:A$ and $F_2, s: \bar{A}$ are also
  semi-nets, with gates $\alpha[a].t$ and $s$ respectively.  We
  proceed by case analysis on the structure of $s$.

  If $s=(\varepsilon[M].s')$, there is a reduction triple between $F$ and 
  \[E = F_1[a:=M], t[a:=M] \bowtie s':(A[x:=M] \ | \ \bar{A}[x:=M]),
  F_2\] which has measure less than that of $F$.  By the induction
  hypothesis, there is a reduction triple $(E', \gterm, \gtaut)$ for
  $E$. By composition, there is a reduction triple between $F$ and $E'$.

  Finally, suppose that $s$ has the form $\varepsilon[M_1].s_1+ \cdots
  +\varepsilon[M_n].s_n$.  Since the relation $\ll$ is a partial order
  on the nodes of $F$, there must be a $\ll$-minimal $\varepsilon[M_i].s_i$
  among the components of $s$, such that $s =  \varepsilon[M_i].s_i+s'$.
  There is a reduction triple between $F$ and 
  $E = E', \alpha[a_0].t_0 \bowtie_Y \varepsilon[M_i].s_i,\alpha[a_1].t_1 \bowtie_X s'.$ 

  Consider now the kingdom $K=k(X)$ of the cut $X$ in $E$.  Since we
  picked $\varepsilon[M_i].s_i$ to be $\ll$-minimal among the
  components of $S$, it does not appear in $k(s')$, and thus does not
  appear in $K$. Since $\varepsilon[M_i].s_i$ is not a member of $K$,
  neither is the cut $Y$.  $K$ is, therefore, an ACC forest of lower
  measure than $F$ (it contains a single cut of nonzero rank, with the
  same rank but lower width than that appearing in $F$) and thus there
  is a reduction triple $(K', \gterm, \gtaut)$ for $K$.  We may substitute $K'$
  for $K$ in $E$, yielding a ACC forest $E[K'/K]$ and functions
  $\hterm$ and $\htaut$ forming a reduction-triple for $E$.  The ACC
  forest $E[K'/K]$ now contains a single nonzero cut of width $1$:
  since $\varepsilon[M_i].s_i$ was not in $K$, the width of this cut
  in $E[K'/K]$ is the same as that in $E$.  $E[K'/K]$ is thus subject
  to the induction hypothesis, which yields a triple $(F', \eterm,
  \etaut)$ for $E[K'/K]$.  We may now compose these 
three reduction triples to obtain
  the required reduction triple for $F$. 
\end{proof} 

As a corollary to the principal lemma, we obtain partial cut-elimination.

\begin{theorem}[Partial cut-elimination]
\label{thm:partial-cut-elimination}
Let $F$ be an Herbrand net.  There is a net $F'$, containing only
cuts of rank zero, with the same conclusion as $F$.
\end{theorem}
\begin{proof}
  By induction on the number of nonzero cuts in an ACC forest $F$.  If there are
  none, we are done.  Now suppose we may remove the nonzero cuts from
  an ACC forest containing $n-1$ nonzero cuts, and let $F$ contain $n$
  nonzero cuts.  Let $X$ be a $\ll$-maximal nonzero cut in $F$, and
  consider $K=k(X)$, it's kingdom.  By the previous lemma, there is a
  nonzero-cut-free reduction triple $(K', \fterm, \ftaut)$ for $K$.
  The ACC forest $F[K'/K]$ has the same type as $F$, but has $n-1$
  nonzero cuts; we may apply the induction hypothesis to obtain our
  nonzero-cut-free net. 
\end{proof}

\subsection{From Partial to Full cut-elimination}
\label{sec:partial-to-full}
Usually, when one performs partial cut-elimination, it is because the
remaining cuts cannot be eliminated.  Here this is not the case: the
cuts of rank zero may very easily be eliminated, but in a way that
interferes with the notion of reduction triple.  The reader
might suspect that here we find a source of nondeterminism in the
reductions: a term $S:P$ where $S$ has cardinality $n$ greater than one,
represents an $n-1$-fold contraction, and so, since we may form cuts
$S\bowtie T$, one might expect to have to make duplications, and to
have to choose a direction in which the cut should be reduced.

In fact, we can avoid such issues, owing to the following lemma:
 
\begin{lemma}
\label{lem:trivial-propositional-cuts}
Let $F= G, (S)\bowtie(T): P \bowtie \bar{P}$ be an Herbrand net, with $G$
cut-free: then $S$ and $T$ are disjoint singleton sets.
\end{lemma}
\begin{proof}
  By the definition of correctness: alternatively, observe that as $F$
  is an Herbrand net it must be the conclusion of an $\LK_H$ derivation
  containing one cut, and thus two branches. 
\end{proof}

\noindent Such cuts are easy to eliminate
\begin{lemma}
\label{lem:eliminating-trivial-propositional-cuts}
  Let $F, (i)\bowtie(j)$ be an Herbrand net.  Then $F[i\leftarrow j]$
  is an ACC forest.
\end{lemma}
\begin{proof}
  By induction on the height of a proof of $F, (i)\bowtie(j)$ in
  $\LK_H$.  Since the proof contains a cut, it cannot have height $1$ -
  the minimal height is $2$, with the proof having the form
\[
\begin{prooftree}
\[
\justifies (i):P_1, \dots (i):P_n, (i):P
\using
i
\]
\qquad
\[
\justifies (j):Q_1, \dots (j):Q_m, (j):\bar{P}
\using
j
\]
\justifies
 (i):P_1, \dots (i):P_n, (j):Q_1, \dots (j):Q_m, (i)\bowtie(j):(P \ |\ \bar{P})
\using
\Cut 
\end{prooftree}
\]
It follows that $\bigvee P_i \vee \bigvee Q_i$ is a tautology, and so 
\[ (i):P_1, \dots, (i):P_n, (i):Q_1, \dots, (i):Q_m \] is the
conclusion of a tautology rule.  The remainder of the proof is a
simple induction on the height of a proof, relying on the fact that
any other rule in $\LK_H$ can be pushed below a cut of the form $(i)
\bowtie (j)$.
\end{proof}

\begin{corollary}
Let $F$ be an Herbrand net containing only cuts of rank $0$.
Then there is an Herbrand net $F'$ of the same type which is cut-free,
which can be obtained by applying the transformation
\[ \text{\textsc{Prop}}:F, (i)\bowtie(j)  \leadsto F[i:= j] \]
\end{corollary}
\begin{proof}
  By induction on the number of cuts in $F$.  Suppose that we may
  remove $n-1$ cuts of zero rank from a net.  Then if $F$ contains $n$
  cuts, it in particular contains one cut of the form
  $\finset{i}\bowtie\finset{j}$, which may be removed by the above
  lemma.  The remaining proof contains $n-1$ cuts and so falls under
  the induction hypothesis.
\end{proof}

This is enough to give us cut-elimination, via the transformations in 
Figure \ref{fig:red}.
\begin{theorem}[Weak Normalization for Herbrand nets]
\label{thm:weak-normalization} 
Let $F$ be an Herbrand net with type $\G$.  By applying rules from
Figure~\ref{fig:red} we may produce an Herbrand net $F'$ which is
cut-free.
\end{theorem}

\begin{figure}
\noindent\hrulefill
\[ \text{\textsc{Prop}}:F, (i)\bowtie(j)  \leadsto F[i:= j] \]
\[ \text{\textsc{Comm}}:F, \ \  
\alpha[a].t \bowtie \finset{\varepsilon[M].s} \leadsto F[a:=M], \ \
t[a:=M] \bowtie s\] 
\[ \text{\textsc{Dup}}: F, \alpha[x].t\bowtie (s_0 + s_1) \leadsto D_x(F), \quad \alpha[x_0].\tau_0(t) \bowtie s_0, \quad   \alpha[x_1].\tau_1(t) \bowtie s_1 \]
\noindent\hrulefill
\caption{Minimal reduction on Herbrand nets}
\label{fig:red}
\end{figure}

\section{Minimal reduction is not confluent}
\label{sec:minim-reduct-not}
Conventional wisdom is that the non-confluence of the sequent-calculus
is a result of structural rules meeting in a cut.  In particular, the
Weakening--Weakening-, or Lafont-example \cite{GLT89} constructs,
given arbitrary proofs $\Phi$ and $\Psi$ of a sequent $\G$, a third
proof $\Phi*\Psi$ of $\G$ which reduces to both $\Phi$ and $\Psi$.  A
Contraction--Contraction cut similarly produces a situation in which
we may reduce a cut in two different directions -- here of course it
is difficult, without a good notion of equality on proofs, to say
whether reduction is nonconfluent.

At first sight, it might appear that we avoid non-confluence in $\LK_H$
and Herbrand nets: each cut reduction rule in Figure~\ref{fig:red} has
a unique orientation.  The only choice we are offered is how to split
the existential side of a \textsc{Dup}-reduct, but we can remove this
difficulty by insisting that we completely decompose such cuts in one
go.  This was, indeed, the original motivation for considering such a
system.

Nevertheless, it emerges that the minimal reduction system on Herbrand
nets is nonconfluent: the nonconfluence arises between, not within,
cuts: that is, the choice we are asked to make is not how to reduce
one particular cut, but instead which cut we should reduce.  

In the remainder of this section we work over a signature and theory
axiomatizing a successor function: $\Sigma = (\mathcal{X},
\finset{\mathsf{0}, \mathsf{s}}, \finset{ \mathsf{iszero}})$ with
$\mathsf{0}$ a constant, $\mathsf{s}$ a unary function symbol, and
$\mathsf{iszero}$ a unary relation symbol.  The universal axiom set
$\univax$ for this theory consists of the single open formula $\lnot
\mathsf{iszero}(\mathsf{s}(x))$ -- a successor is never zero.
The signature is necessary, not to obtain a nonconfluent example, but
to exhibit the nature of the nonconfluence. We give a proof
of a $\Sigma_1$ formula such that, when we eliminate the cuts in
two different ways, we get two different sets of witnesses to that formula:
the sets differ in size, and with the help of our signature, also in 
what terms they contain.

Our example net is the net from
Example~\ref{ex:cut-elim-example}, whose dependency graph is the
following:

 
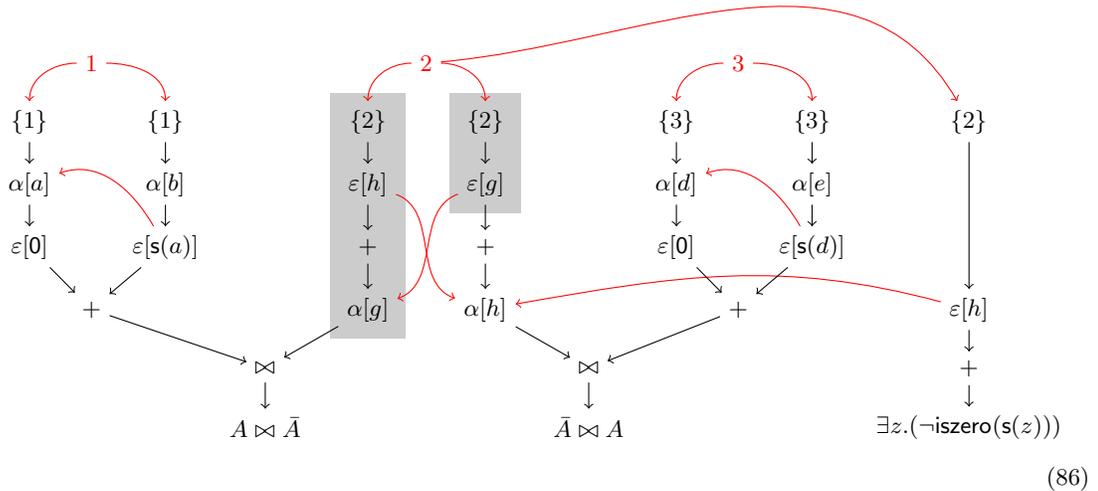
\begin{equation}
\label{eq:cutelimmexample}
\begin{tikzpicture}
\matrix[row sep=3mm,column sep=2mm]{
& \node [red](tautology) {$1$}; &&&&& 
\node [red](tautology 2) {$2$}; &&&&& 
\node [red](tautology 3) {$3$}; \\
 \node (left prop) {$\finset{1}$}; &&
 \node (right prop) {$\finset{1}$}; &&& \node (left prop 2) {$\finset{2}$}; &&
 \node (right prop 2) {$\finset{2}$};&&& \node (left prop 3) {$\finset{3}$}; &&
 \node (right prop 3) {$\finset{3}$}; & \node(prop out){$\finset{2}$};\\
\node (left alpha) {$\alpha[a]$}; &&
 \node (right alpha) {$\alpha[b]$}; &&& 
\node (left epsilon 2) {$\varepsilon[h]$}; &&
 \node (right epsilon 2) {$\varepsilon[g]$};&&&
 \node (left alpha 3) {$\alpha[d]$}; &&
 \node (right alpha 3) {$\alpha[e]$}; & \\
\node (left epsilon) {$\varepsilon[\mathsf{0}]$}; &&
 \node (right epsilon) {$\varepsilon[\mathsf{s}(a)]$}; &&&
\node (left plus 2) {+}; &&
 \node (right plus 2) {+}; &&&
\node (left epsilon 3) {$\varepsilon[\mathsf{0}]$}; &&
 \node (right epsilon 3) {$\varepsilon[\mathsf{s}(d)]$}; &\\
 & \node(plus){$+$}; &&&&
\node (left alpha 2) {$\alpha[g]$}; &&
 \node (right alpha 2) {$\alpha[h]$}; &&&&
\node(plus 3){$+$}; && \node (epsilon out) {$\varepsilon[h]$};\\
&&& \node (left cut) {$\bowtie$}; &&&&&& \node (right cut) {$\bowtie$};
&&&&   \node(plus out){$+$};\\
&&& \node (left out) {$A \bowtie \bar{A}$}; &&&&&& \node (right out) {$\bar{A} \bowtie A$}; &&&& 
\node(out out){$\exists z. (\lnot \mathsf{iszero}(\mathsf{s}(z)))$};\\
};



\draw[->] (prop out) to (epsilon out);
\draw[->] (epsilon out) to (plus out);
\draw[->] (plus out) to (out out);
\draw[->](plus) to (left cut);
\draw[->](left alpha 2) to (left cut);
\draw[->](right alpha 2) to (right cut);
\draw[->](plus 3) to (right cut);
\draw[->](right cut) to (right out);
\draw[->](left cut) to (left out);
\draw[->](left prop) to (left alpha);
\draw[->](right prop) to (right alpha);
\draw[->](left alpha) to (left epsilon);
\draw[->](right alpha) to (right epsilon);
\draw[->](left epsilon) to (plus);
\draw[->](right epsilon) to (plus);

\draw[->](left prop 2) to (left epsilon 2);
\draw[->](right prop 2) to (right epsilon 2);
\draw[->](left epsilon 2) to (left plus 2);
\draw[->](right epsilon 2) to (right plus 2);
\draw[->](left plus 2) to (left alpha 2);
\draw[->](right plus 2) to (right alpha 2);

\draw[red,->](tautology 2) to [out=180, in = 90] (left prop 2);
\draw[red,->](tautology 2) to [out=0, in = 90] (right prop 2);
\draw[red,->](tautology 2) to [out=5, in = 120] (prop out);
\draw[red,->](right epsilon 2) to [out=200, in = 20] (left alpha 2);
\draw[red,->](left epsilon 2) to [out=-20, in = 160] (right alpha 2);
\draw[red,->](epsilon out) to [out=165, in = 10] (right alpha 2);

\draw[red,->](tautology) to [out=180, in = 90] (left prop);
\draw[red,->](tautology) to [out=0, in = 90] (right prop);
\draw[red,->](right epsilon) to [out=120, in = 20] (left alpha);

\draw[->](left prop 3) to (left alpha 3);
\draw[->](right prop 3) to (right alpha 3);
\draw[->](left alpha 3) to (left epsilon 3);
\draw[->](right alpha 3) to (right epsilon 3);
\draw[->](left epsilon 3) to (plus 3);
\draw[->](right epsilon 3) to (plus 3);

\draw[red,->](tautology 3) to [out=180, in = 90] (left prop 3);
\draw[red,->](tautology 3) to [out=0, in = 90] (right prop 3);
\draw[red,->](right epsilon 3) to [out=120, in = 20] (left alpha 3);
\begin{pgfonlayer}{background}
\node [fill=black!20,fit=(left alpha 2) (left prop 2)] {};
\node [fill=black!20,fit=(right epsilon 2) (right prop 2)] {};
\end{pgfonlayer}
\end{tikzpicture}
\end{equation}
\noindent (The greyed
nodes indicate the kingdom $K$ of the node $\alpha[g]$: we will later
use this subnet to begin the elimination of cuts from this net).
This net is a proof, in the system mentioned above, that there is a
nonzero element.  Our interest will be in the result of eliminating
the cuts from this proof: which witnesses to the theorem do we obtain?

To see that \eqref{eq:cutelimmexample} is an Herbrand net, observe
that it is the result of cutting together three cut-free Herbrand nets
(one could, of course, check the correctness criterion directly).  Two
are of the form of the drinker's formula, but with the addition of a
function symbol:
\[
(\varepsilon[\mathsf{0}].\alpha[a].\finset{1}+\varepsilon[\mathsf{s}(a)].\alpha[b].\finset{1}):A\]

\noindent and 
\[
(\varepsilon[\mathsf{0}].\alpha[d].\finset{3}+\varepsilon[\mathsf{s}(d)].\alpha[e].\finset{3}): A \]
\noindent while the last is the following
\[ \alpha[g].(\varepsilon[h].\finset{2}):\bar{A}, \quad  \alpha[h].(\varepsilon[g].\finset{2}):\forall x. \exists y. \bar{A}, \quad(\varepsilon[h].\finset{3}):\exists z.(\lnot \mathsf{iszero}(z)).\]

We leave it as a simple exercise to check that these are Herbrand nets. 
We will reduce this net in two ways, obtaining distinct cut-free nets.

To begin, we reduce \eqref{eq:cutelimmexample}  by a \textsc{Dup}-reduction applied
to the left-hand cut, which duplicates the shaded subnet $K$, the kingdom of the
node $\alpha[g]$.  The following net
is the result, with the two copies of the shaded subnet displayed:

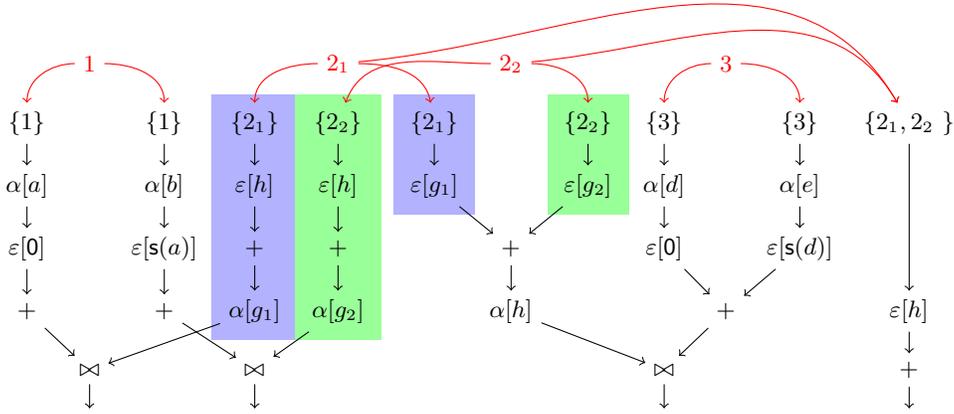
\[
\begin{tikzpicture}
\matrix[row sep=3mm,column sep=2mm]{
& \node [red](tautology) {$1$}; &&&
\node [red](tautology 21) {$2_1$}; &&&\node [red](tautology 22) {$2_2$}; &&&
\node [red](tautology 3) {$3$}; \\
 \node (left prop) {$\finset{1}$}; &&
 \node (right prop) {$\finset{1}$}; & \node (left prop 21) {$\finset{2_1}$};
& \node (left prop 22) {$\finset{2_2}$}; &&
 \node (right prop 21) {$\finset{2_1}$};&&\node (right prop 22) {$\finset{2_2}$};& \node (left prop 3) {$\finset{3}$}; &&
 \node (right prop 3) {$\finset{3}$}; & \node (prop out){$\finset{2_1,2_2$}};\\
\node (left alpha) {$\alpha[a]$}; &&
 \node (right alpha) {$\alpha[b]$}; &\node (left epsilon 21) {$\varepsilon[h]$};&
\node (left epsilon 22) {$\varepsilon[h]$}; &&
 \node (right epsilon 21) {$\varepsilon[g_1]$};&&\node (right epsilon 22) {$\varepsilon[g_2]$};&
 \node (left alpha 3) {$\alpha[d]$}; &&
 \node (right alpha 3) {$\alpha[e]$}; \\
\node (left epsilon) {$\varepsilon[\mathsf{0}]$}; &&
 \node (right epsilon) {$\varepsilon[\mathsf{s}(a)]$}; &\node (left plus 21) {+};&
\node (left plus 22) {+}; &&&
 \node (right plus 2) {+}; &&
\node (left epsilon 3) {$\varepsilon[\mathsf{0}]$}; &&
 \node (right epsilon 3) {$\varepsilon[\mathsf{s}(d)]$};\\
 \node(left plus 1){$+$}; &&\node(left plus 2){$+$};&\node (left alpha 21) {$\alpha[g_1]$}; &
\node (left alpha 22) {$\alpha[g_2]$}; &&&
 \node (right alpha 2) {$\alpha[h]$}; &&&
\node(plus 3){$+$}; && \node(epsilon out){$\varepsilon[h]$}; \\
& \node(left cut 1) {$\bowtie$}; && \node (left cut 2) {$\bowtie$}; &&&&&& \node (right cut) {$\bowtie$}; &&& \node(plus out){$+$};\\
&\node (left out 1) {};&& \node (left out 2) {}; &&&&&& \node (right out) {}; &&& \node(out out) {};\\
};


\draw[->] (prop out) to (epsilon out);
\draw[->] (epsilon out) to (plus out);
\draw[->] (plus out) to (out out);
\draw[->](left plus 1) to (left cut 1);
\draw[->](left plus 2) to (left cut 2);
\draw[->](left alpha 21) to (left cut 1);
\draw[->](left alpha 22) to (left cut 2);
\draw[->](right alpha 2) to (right cut);
\draw[->](plus 3) to (right cut);
\draw[->](right cut) to (right out);
\draw[->](left cut 1) to (left out 1);
\draw[->](left cut 2) to (left out 2);
\draw[->](left prop) to (left alpha);
\draw[->](right prop) to (right alpha);
\draw[->](left alpha) to (left epsilon);
\draw[->](right alpha) to (right epsilon);
\draw[->](left epsilon) to (left plus 1);
\draw[->](right epsilon) to (left plus 2);

\draw[->](left prop 21) to (left epsilon 21);
\draw[->](left prop 22) to (left epsilon 22);
\draw[->](right prop 21) to (right epsilon 21);
\draw[->](right prop 22) to (right epsilon 22);
\draw[->](left epsilon 21) to (left plus 21);
\draw[->](left epsilon 22) to (left plus 22);

\draw[->](right epsilon 21) to (right plus 2);
\draw[->](right epsilon 22) to (right plus 2);
\draw[->](left plus 21) to (left alpha 21);
\draw[->](left plus 22) to (left alpha 22);
\draw[->](right plus 2) to (right alpha 2);

\draw[red,->](tautology 21) to [out=180, in = 90] (left prop 21);
\draw[red,->](tautology 22) to [out=180, in = 70] (left prop 22);
\draw[red,->](tautology 21) to [out=0, in = 100] (right prop 21);
\draw[red,->](tautology 22) to [out=0, in = 90] (right prop 22);
\draw[red,->](tautology 21) to [out=5, in = 120] (prop out);
\draw[red,->](tautology 22) to [out=5, in = 120] (prop out);
\draw[red,->](tautology) to [out=180, in = 90] (left prop);
\draw[red,->](tautology) to [out=0, in = 90] (right prop);

\draw[->](left prop 3) to (left alpha 3);
\draw[->](right prop 3) to (right alpha 3);
\draw[->](left alpha 3) to (left epsilon 3);
\draw[->](right alpha 3) to (right epsilon 3);
\draw[->](left epsilon 3) to (plus 3);
\draw[->](right epsilon 3) to (plus 3);

\draw[red,->](tautology 3) to [out=180, in = 90] (left prop 3);
\draw[red,->](tautology 3) to [out=0, in = 90] (right prop 3);
\begin{pgfonlayer}{background}
\node [fill=blue!30, fit=(left alpha 21) (left prop 21)] {};
\node [fill=blue!30, fit=(right epsilon 21) (right prop 21)] {};
\node [fill=green!40,fit=(left alpha 22) (left prop 22)] {};
\node [fill=green!40,fit=(right epsilon 22) (right prop 22)] {};
\end{pgfonlayer}
\end{tikzpicture}
\]

(From this point on, we will only show jumps where they serve to
clarify the situation).  The net $K$ had, in addition to $\alpha[g]$,
one root $\varepsilon[g]$, the child of an expansion node outside of
$K$. Notice that, in the reduct, both copies of $\varepsilon[g]$ are
children of that expansion node.  Notice also that the node labelled
$\finset{2}$ outside of $K$ becomes $\finset{2_1, 2_2}$ after
reduction.

To continue the reduction of this net, we perform four \textsc{Comm}
reductions, in which the $\varepsilon$ nodes transmit their
first-order terms to the corresponding $\alpha$ nodes.  It does not
matter in which order these communications take place.  Here is the
result after two communications:
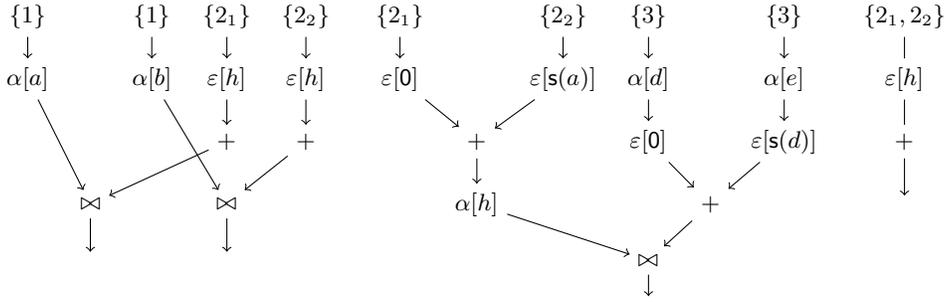
\[
\begin{tikzpicture}
\matrix[row sep=3mm,column sep=2mm]{
 \node (left prop) {$\finset{1}$}; &&
 \node (right prop) {$\finset{1}$}; & \node (left prop 21) {$\finset{2_1}$};
& \node (left prop 22) {$\finset{2_2}$}; &&
 \node (right prop 21) {$\finset{2_1}$};&&\node (right prop 22) {$\finset{2_2}$};& \node (left prop 3) {$\finset{3}$}; &&
 \node (right prop 3) {$\finset{3}$}; &&  \node (prop out) {$\finset{2_1,2_2}$}; \\
\node (left alpha) {$\alpha[a]$}; &&
 \node (right alpha) {$\alpha[b]$}; &\node (left epsilon 21) {$\varepsilon[h]$};&
\node (left epsilon 22) {$\varepsilon[h]$}; &&
 \node (right epsilon 21) {$\varepsilon[\mathsf{0}]$};&&\node (right epsilon 22) {$\varepsilon[\mathsf{s}(a)]$};&
 \node (left alpha 3) {$\alpha[d]$}; &&
 \node (right alpha 3) {$\alpha[e]$}; && \node(epsilon out){$\varepsilon[h]$};\\
&&&\node (left plus 21) {+};&
\node (left plus 22) {+}; &&&
 \node (right plus 2) {+}; &&
\node (left epsilon 3) {$\varepsilon[\mathsf{0}]$}; &&
 \node (right epsilon 3) {$\varepsilon[\mathsf{s}(d)]$};&& \node(plus out){$+$};\\
 & \node(left cut 1) {$\bowtie$}; && \node (left cut 2) {$\bowtie$}; &
 &&&
 \node (right alpha 2) {$\alpha[h]$}; &&&
\node(plus 3){$+$}; &&& \node(out out){};\\
&\node (left out 1) {};&& \node (left out 2) {}; &&&&&& \node (right cut) {$\bowtie$};\\
&&&&&&&&& \node (right out) {};\\
};



\draw[->](prop out)--(epsilon out)--(plus out)--(out out);

\draw[->](left alpha) to (left cut 1);
\draw[->](right alpha) to (left cut 2);
\draw[->](left plus 21) to (left cut 1);
\draw[->](left plus 22) to (left cut 2);
\draw[->](right alpha 2) to (right cut);
\draw[->](plus 3) to (right cut);
\draw[->](right cut) to (right out);
\draw[->](left cut 1) to (left out 1);
\draw[->](left cut 2) to (left out 2);
\draw[->](left prop) to (left alpha);
\draw[->](right prop) to (right alpha);

\draw[->](left prop 21) to (left epsilon 21);
\draw[->](left prop 22) to (left epsilon 22);
\draw[->](right prop 21) to (right epsilon 21);
\draw[->](right prop 22) to (right epsilon 22);
\draw[->](left epsilon 21) to (left plus 21);
\draw[->](left epsilon 22) to (left plus 22);

\draw[->](right epsilon 21) to (right plus 2);
\draw[->](right epsilon 22) to (right plus 2);
\draw[->](right plus 2) to (right alpha 2);



\draw[->](left prop 3) to (left alpha 3);
\draw[->](right prop 3) to (right alpha 3);
\draw[->](left alpha 3) to (left epsilon 3);
\draw[->](right alpha 3) to (right epsilon 3);
\draw[->](left epsilon 3) to (plus 3);
\draw[->](right epsilon 3) to (plus 3);


\end{tikzpicture}
\]

\noindent And after two more \textsc{Comm}s:

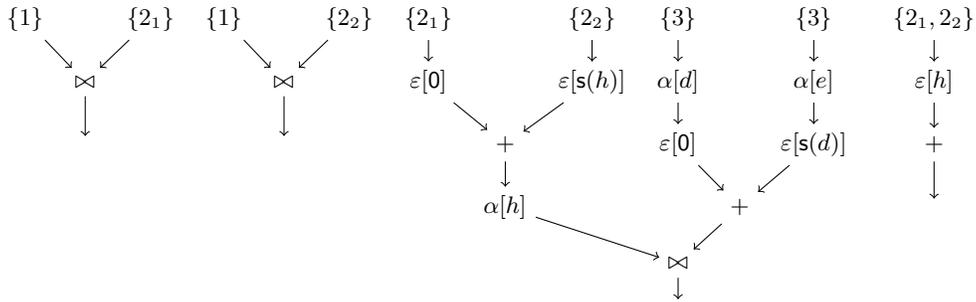
\[ 
\begin{tikzpicture}
\matrix[row sep=3mm,column sep=2mm]{
 \node (left prop) {$\finset{1}$}; &&
 \node (left prop 21) {$\finset{2_1}$}; & \node (right prop) {$\finset{1}$};
&& \node (left prop 22) {$\finset{2_2}$}; &
 \node (right prop 21) {$\finset{2_1}$};&&\node (right prop 22) {$\finset{2_2}$};& \node (left prop 3) {$\finset{3}$}; &&
 \node (right prop 3) {$\finset{3}$}; 
&& \node (out prop){$ \finset{2_1, 2_2}$};\\
& \node(left cut 1) {$\bowtie$}; &&& \node (left cut 2) {$\bowtie$};&
&
 \node (right epsilon 21) {$\varepsilon[\mathsf{0}]$};&&\node (right epsilon 22) {$\varepsilon[\mathsf{s}(h)]$};&
 \node (left alpha 3) {$\alpha[d]$}; &&
 \node (right alpha 3) {$\alpha[e]$}; && \node(out epsilon){$\varepsilon[h]$};\\
 &\node (left out 1) {};&&& \node (left out 2) {}; &
&&
 \node (right plus 2) {+}; &&
\node (left epsilon 3) {$\varepsilon[\mathsf{0}]$}; &&
 \node (right epsilon 3) {$\varepsilon[\mathsf{s}(d)]$}; && \node (out plus){$+$};\\
&&&&
 &&&
 \node (right alpha 2) {$\alpha[h]$}; &&&
\node(plus 3){$+$};&&& \node(out out){};\\
&&& &&&&&& \node (right cut) {$\bowtie$};\\
&&&&&&&&& \node (right out) {};\\
};



\draw[->](out prop) to (out epsilon);
\draw[->](out epsilon) to (out plus);
\draw[->](out plus) to (out out);
\draw[->](right alpha 2) to (right cut);
\draw[->](plus 3) to (right cut);
\draw[->](right cut) to (right out);
\draw[->](left cut 1) to (left out 1);
\draw[->](left cut 2) to (left out 2);
\draw[->](left prop) to (left cut 1);
\draw[->](right prop) to (left cut 2);

\draw[->](left prop 21) to (left cut 1);
\draw[->](left prop 22) to (left cut 2);
\draw[->](right prop 21) to (right epsilon 21);
\draw[->](right prop 22) to (right epsilon 22);

\draw[->](right epsilon 21) to (right plus 2);
\draw[->](right epsilon 22) to (right plus 2);
\draw[->](right plus 2) to (right alpha 2);



\draw[->](left prop 3) to (left alpha 3);
\draw[->](right prop 3) to (right alpha 3);
\draw[->](left alpha 3) to (left epsilon 3);
\draw[->](right alpha 3) to (right epsilon 3);
\draw[->](left epsilon 3) to (plus 3);
\draw[->](right epsilon 3) to (plus 3);


\end{tikzpicture}
\]

\noindent Two applications of the tautology reduction leave a net with only one 
cut remaining, replacing the three tautologies $1$, $2_1$ and $2_2$
with a single tautology $4$.

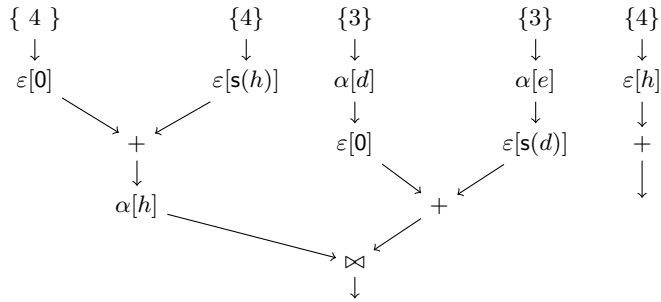
\begin{equation}
\begin{tikzpicture}
\matrix[row sep=3mm,column sep=5mm]{
\node (right prop 21) {\finset{4}}; &&\node (right prop 22) {$\finset{4}$};& \node (left prop 3) {$\finset{3}$}; &&
 \node (right prop 3) {$\finset{3}$}; & \node(out prop){$\finset{4}$};\\
 \node (right epsilon 21) {$\varepsilon[\mathsf{0}]$};&&\node (right epsilon 22) {$\varepsilon[\mathsf{s}(h)]$};&
 \node (left alpha 3) {$\alpha[d]$}; &&
 \node (right alpha 3) {$\alpha[e]$}; & \node(out epsilon){$\varepsilon[h]$};\\
&
 \node (right plus 2) {+}; &&
\node (left epsilon 3) {$\varepsilon[\mathsf{0}]$}; &&
 \node (right epsilon 3) {$\varepsilon[\mathsf{s}(d)]$};& \node(out plus){$+$};\\
 &
 \node (right alpha 2) {$\alpha[h]$}; &&&
\node(plus 3){$+$}; && \node(out out){};\\
&&& \node (right cut) {$\bowtie$};\\
&&& \node (right out) {};\\
};



\draw[->] (out prop) to (out epsilon);
\draw[->] (out epsilon) to (out plus);
\draw[->] (out plus) to (out out);
\draw[->](right alpha 2) to (right cut);
\draw[->](plus 3) to (right cut);
\draw[->](right cut) to (right out);
\draw[->](right prop 21) to (right epsilon 21);
\draw[->](right prop 22) to (right epsilon 22);

\draw[->](right epsilon 21) to (right plus 2);
\draw[->](right epsilon 22) to (right plus 2);
\draw[->](right plus 2) to (right alpha 2);



\draw[->](left prop 3) to (left alpha 3);
\draw[->](right prop 3) to (right alpha 3);
\draw[->](left alpha 3) to (left epsilon 3);
\draw[->](right alpha 3) to (right epsilon 3);
\draw[->](left epsilon 3) to (plus 3);
\draw[->](right epsilon 3) to (plus 3);


\end{tikzpicture}
\end{equation}

\noindent To reduce the remaining cut, we must first apply \textsc{Dup}, duplicating the
kingdom of $\alpha[h]$:
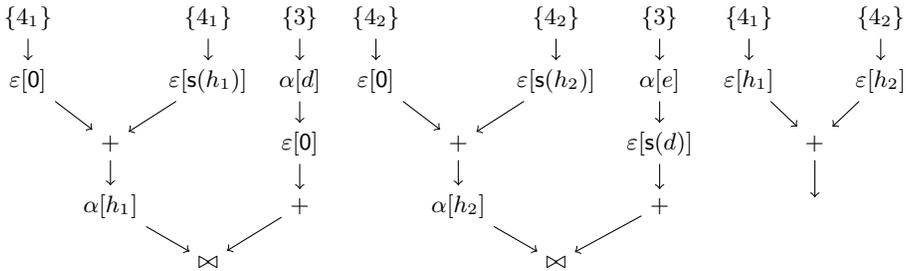
\[
\begin{tikzpicture}
\matrix[row sep=3mm,column sep=2mm]{
\node (left prop 1) {$\finset{4_1}$}; && \node (right prop 1) {$\finset{4_1}$};
&\node (left prop 3) {$\finset{3}$}; &
\node (left prop 2) {$\finset{4_2}$}; && \node (right prop 2) {$\finset{4_2}$};&  
 \node (right prop 3) {$\finset{3}$}; & \node (left prop out) {$\finset{4_1}$};
&&\node (right prop out) {$\finset{4_2}$};\\
 \node (left epsilon 1) {$\varepsilon[\mathsf{0}]$};
&&\node (right epsilon 1) {$\varepsilon[\mathsf{s}(h_1)]$}; & \node (left alpha 3) {$\alpha[d]$};&
 \node (left epsilon 2) {$\varepsilon[\mathsf{0}]$};&&\node (right epsilon 2) {$\varepsilon[\mathsf{s}(h_2)]$};&
 \node (right alpha 3) {$\alpha[e]$}; & \node(left epsilon out){$\varepsilon[h_1] $}; 
&& \node(right  epsilon out){$\varepsilon[h_2]$};\\
&
 \node (left plus 1) {+}; &&\node (left epsilon 3) {$\varepsilon[\mathsf{0}]$};
&& \node (left plus 2) {+}; &&\node (right epsilon 3) {$\varepsilon[\mathsf{s}(d)]$}; && \node(plus out){$+$};
 \\
 &
 \node (left alpha 1) {$\alpha[h_1]$}; &&\node(left plus 3){+};&& \node (left alpha 2) {$\alpha[h_2]$};&&
\node(right plus 3){$+$}; && \node(out out){};\\
&& \node (left cut) {$\bowtie$};&&&&\node (right cut) {$\bowtie$};\\
&& \node (left out) {}; &&&& && \node (right out) {};\\
};

\draw[->] (left prop out) to (left epsilon out);
\draw[->] (right prop out) to (right epsilon out);
\draw[->] (left epsilon out) to (plus out);
\draw[->] (right epsilon out) to (plus out);
\draw[->] (plus out) to (out out);
\draw[->] (left prop 1) to (left epsilon 1);
\draw[->] (left prop 2) to (left epsilon 2);
\draw[->] (right prop 1) to (right epsilon 1);
\draw[->] (right prop 2) to (right epsilon 2);
\draw[->] (left epsilon 1) to (left plus 1);
\draw[->] (left epsilon 2) to (left plus 2);
\draw[->] (right epsilon 1) to (left plus 1);
\draw[->] (right epsilon 2) to (left plus 2);
\draw[->] (left plus 1) to (left alpha 1);
\draw[->] (left plus 2) to (left alpha 2);
\draw[->] (left alpha 1) to (left cut);
\draw[->] (left alpha 2) to (right cut);

\draw[->] (left prop 3) to (left alpha 3);
\draw[->] (right prop 3) to (right alpha 3);
\draw[->] (left alpha 3) to (left epsilon 3);
\draw[->] (right alpha 3) to (right epsilon 3);
\draw[->] (left epsilon 3) to (left plus 3);
\draw[->] (right epsilon 3) to (right plus 3);
\draw[->] (left plus 3) to (left cut);
\draw[->] (right plus 3) to (right cut);




\end{tikzpicture}
\]

\noindent We are here presented with a choice of which cut to reduce:
we pick the rightmost, as it will involve fewer reduction steps.  A
Comm reduction leads us to
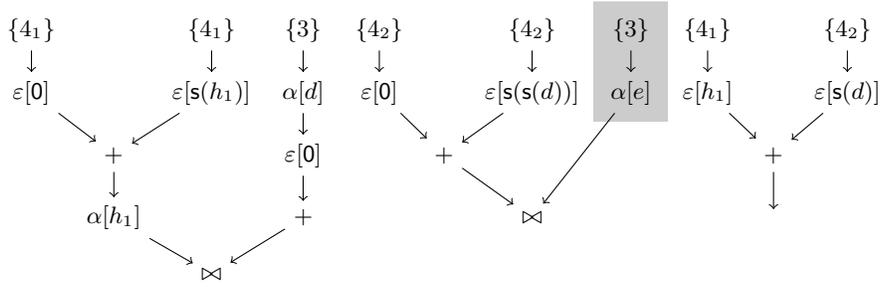
\[ 
\begin{tikzpicture}
\matrix[row sep=3mm,column sep=2mm]{
\node (left prop 1) {$\finset{4_1}$}; && \node (right prop 1) {$\finset{4_1}$};
&\node (left prop 3) {$\finset{3}$}; &
\node (left prop 2) {$\finset{4_2}$}; && \node (right prop 2) {$\finset{4_2}$};&  
 \node (right prop 3) {$\finset{3}$};  & \node (left prop out) {$\finset{4_1}$};
&&\node (right prop out) {$\finset{4_2}$};\\
 \node (left epsilon 1) {$\varepsilon[\mathsf{0}]$};
&&\node (right epsilon 1) {$\varepsilon[\mathsf{s}(h_1)]$}; & \node (left alpha 3) {$\alpha[d]$};&
 \node (left epsilon 2) {$\varepsilon[\mathsf{0}]$};&&\node (right epsilon 2) {$\varepsilon[\mathsf{s}(\mathsf{s}(d))]$};&
 \node (right alpha 3) {$\alpha[e]$}; & \node(left epsilon out){$\varepsilon[h_1] $}; 
&& \node(right  epsilon out){$\varepsilon[\mathsf{s}(d)]$};\\
&
 \node (left plus 1) {+}; &&\node (left epsilon 3) {$\varepsilon[\mathsf{0}]$};
&&\node (left plus 2) {+};&&&& \node(plus out){$+$};
 \\
 &
 \node (left alpha 1) {$\alpha[h_1]$}; &&\node(left plus 3){+};&&&\node (right cut) {$\bowtie$}; 
&&& \node(out out){};\\
&& \node (left cut) {$\bowtie$};&&&&\\
&& \node (left out) {}; &&&& && \node (left out) {};\\
};

\draw[->] (left prop out) to (left epsilon out);
\draw[->] (right prop out) to (right epsilon out);
\draw[->] (left epsilon out) to (plus out);
\draw[->] (right epsilon out) to (plus out);
\draw[->] (plus out) to (out out);

\draw[->] (left prop 1) to (left epsilon 1);
\draw[->] (left prop 2) to (left epsilon 2);
\draw[->] (right prop 1) to (right epsilon 1);
\draw[->] (right prop 2) to (right epsilon 2);
\draw[->] (left epsilon 1) to (left plus 1);
\draw[->] (left epsilon 2) to (left plus 2);
\draw[->] (right epsilon 1) to (left plus 1);
\draw[->] (right epsilon 2) to (left plus 2);
\draw[->] (left plus 1) to (left alpha 1);
\draw[->] (left alpha 1) to (left cut);
\draw[->] (left plus 2) to (right cut);

\draw[->] (left prop 3) to (left alpha 3);
\draw[->] (right prop 3) to (right alpha 3);
\draw[->] (left alpha 3) to (left epsilon 3);
\draw[->] (left epsilon 3) to (left plus 3);
\draw[->] (left plus 3) to (left cut);
\draw[->] (right alpha 3) to (right cut);



\begin{pgfonlayer}{background}
\node [fill=black!20,fit=(right alpha 3) (right prop 3)] {};
\end{pgfonlayer}
\end{tikzpicture}
\]

\noindent where again we must make a duplication; since however, the
eigenvariable $e$ appears nowhere else in the proof, the elimination
of the middle cut has little effect on the remaining proof: the reader
may verify that, after one duplication, two communications and two
tautology reductions, we arrive at the following proof:

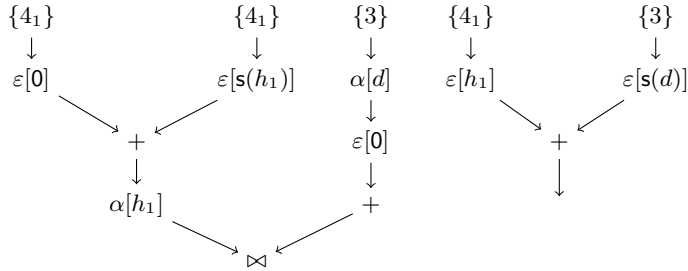
\[
\begin{tikzpicture}
\matrix[row sep=3mm,column sep=5mm]{
\node (left prop 1) {$\finset{4_1}$}; && \node (right prop 1) {$\finset{4_1}$};
&\node (left prop 3) {$\finset{3}$}; & \node (left prop out) {$\finset{4_1}$};
&&\node (right prop out) {$\finset{3}$}; \\
 \node (left epsilon 1) {$\varepsilon[\mathsf{0}]$};
&&\node (right epsilon 1) {$\varepsilon[\mathsf{s}(h_1)]$}; & \node (left alpha 3) {$\alpha[d]$};&
  \node(left epsilon out){$\varepsilon[h_1] $}; 
&& \node(right  epsilon out){$\varepsilon[\mathsf{s}(d)]$};\\
& \node (left plus 1) {$+$}; && \node (left epsilon 2) {$\varepsilon[\mathsf{0}]$};
&& \node(plus out){$+$};\\
 &
 \node (left alpha 1) {$\alpha[h_1]$}; &&\node(left plus 3){+};&& \node(out out){};\\
&& \node (left cut) {$\bowtie$};&&&&\\
&& \node (left out) {}; &&&& && \node (left out) {};\\
};

\draw[->] (left prop out) to (left epsilon out);
\draw[->] (right prop out) to (right epsilon out);
\draw[->] (left epsilon out) to (plus out);
\draw[->] (right epsilon out) to (plus out);
\draw[->] (plus out) to (out out);
\draw[->] (right prop 1) to (right epsilon 1);
\draw[->] (left prop 1) to (left epsilon 1);
\draw[->] (left epsilon 1) to (left plus 1);
\draw[->] (right epsilon 1) to (left plus 1);
\draw[->] (left plus 1) to (left alpha 1);
\draw[->] (left prop 3) to (left alpha 3);
\draw[->] (left alpha 3) to (left epsilon 2);
\draw[->] (left epsilon 2) to (left plus 3);
\draw[->] (left plus 3) to (left cut);
\draw[->] (left alpha 1) to (left cut);




\end{tikzpicture}
\]

\noindent We now communicate the term $\mathsf{0}$ into the eigenvariable $h_1$:
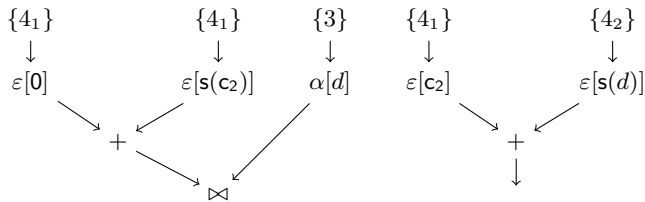
\[
\begin{tikzpicture}
\matrix[row sep=3mm,column sep=5mm]{
\node (left prop 1) {$\finset{4_1}$}; && \node (right prop 1) {$\finset{4_1}$};
&\node (left prop 3) {$\finset{3}$}; & \node (left prop out) {$\finset{4_1}$};
&&\node (right prop out) {$\finset{4_2}$}; \\
 \node (left epsilon 1) {$\varepsilon[\mathsf{0}]$};
&&\node (right epsilon 1) {$\varepsilon[\mathsf{s}(\mathsf {c_2})]$}; & \node (left alpha 3) {$\alpha[d]$};&
 \node(left epsilon out){$\varepsilon[\mathsf {c_2}] $}; 
&& \node(right  epsilon out){$\varepsilon[\mathsf{s}(d)]$}; \\
& \node (left plus 1) {$+$}; &&&& \node(plus out){$+$};\\
&& \node (left cut) {$\bowtie$};&&& \node(out out){};&\\
&& \node (left out) {}; &&&& && \node (left out) {};\\
};

\draw[->] (left prop out) to (left epsilon out);
\draw[->] (right prop out) to (right epsilon out);
\draw[->] (left epsilon out) to (plus out);
\draw[->] (right epsilon out) to (plus out);
\draw[->] (plus out) to (out out);
\draw[->] (right prop 1) to (right epsilon 1);
\draw[->] (left prop 1) to (left epsilon 1);
\draw[->] (left epsilon 1) to (left plus 1);
\draw[->] (right epsilon 1) to (left plus 1);
\draw[->] (left plus 1) to (left cut);
\draw[->] (left prop 3) to (left alpha 3);
\draw[->] (left alpha 3) to (left cut);




\end{tikzpicture}
\]

\noindent The resulting net is of a rather simple form: one
application of \textsc{Dup}, two applications of \textsc{Comm} and two
applications of \textsc{Prop} result in a cut-free net: inuitively, we
substitute both of the terms $\mathsf{0}$ and $\mathsf{s}\mathsf{0}$
for $d$:
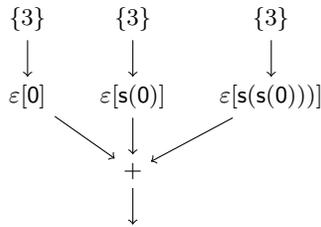
\[
\begin{tikzpicture}
\matrix[row sep=5mm,column sep=5mm]{
\node (output tautology 1){$\finset{3}$}; & \node (output tautology 2){$\finset{3}$};& \node (output tautology 3){$\finset{3}$}; \\
\node (epsilon 1) {$\varepsilon[\mathsf{0}]$}; &
\node (epsilon 2) {$\varepsilon[\mathsf{s}(\mathsf{0})]$}; & \node (epsilon 3) {$\varepsilon[\mathsf{s}(\mathsf{s}(\mathsf{0})))]$}; \\
&\node (plus) {+};\\
&\node (out) {};\\
};
\draw [->] (output tautology 1) to (epsilon 1);
\draw [->] (output tautology 2) to (epsilon 2);
\draw [->] (output tautology 3) to (epsilon 3);
\draw [->] (epsilon 1) to (plus);
\draw [->] (epsilon 2) to (plus);
\draw [->] (epsilon 3) to (plus);
\draw [->] (plus) to (out);
\end{tikzpicture}
\]

\noindent The result of eliminating the cut is a net comprising of three
distinct witnesses.

We now sketch the reduction beginning instead with the right-hand cut.
Rather than repeat the steps above we summarize the reduction as follows:
beginning with \eqref{eq:cutelimmexample}, we instead duplicate the 
kingdom of $\alpha[h]$.  After four \textsc{Comm}s and two aplications
of \textsc{Prop} the resulting net is

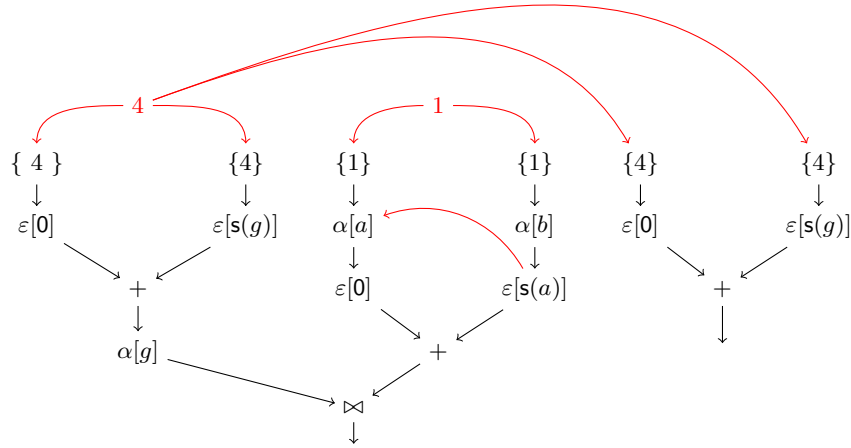
\[
\begin{tikzpicture}
\matrix[row sep=3mm,column sep=5mm]{
&\node [red](tautology 22) {$4$}; &&&
\node [red](tautology 3) {$1$}; \\
\node (right prop 21) {\finset{4}}; &&\node (right prop 22) {$\finset{4}$};& \node (left prop 3) {$\finset{1}$}; &&
 \node (right prop 3) {$\finset{1}$}; & \node(left out prop){$\finset{4}$};
&& \node(right out prop){$\finset{4}$};\\
 \node (right epsilon 21) {$\varepsilon[\mathsf{0}]$};&&\node (right epsilon 22) {$\varepsilon[\mathsf{s}(g)]$};&
 \node (left alpha 3) {$\alpha[a]$}; &&
 \node (right alpha 3) {$\alpha[b]$}; & \node(left out epsilon){$\varepsilon[\mathsf{0}]$};
&& \node(right out epsilon){$\varepsilon[\mathsf{s}(g)]$};\\
&
 \node (right plus 2) {+}; &&
\node (left epsilon 3) {$\varepsilon[\mathsf{0}]$}; &&
 \node (right epsilon 3) {$\varepsilon[\mathsf{s}(a)]$};&& \node(out plus){$+$};\\
 &
 \node (right alpha 2) {$\alpha[g]$}; &&&
\node(plus 3){$+$}; &&& \node(out out){};\\
&&& \node (right cut) {$\bowtie$};\\
&&& \node (right out) {};\\
};



\draw[->] (left out prop) to (left out epsilon);
\draw[->] (right out prop) to (right out epsilon);
\draw[->] (left out epsilon) to (out plus);
\draw[->] (right out epsilon) to (out plus);
\draw[->] (out plus) to (out out);
\draw[->](right alpha 2) to (right cut);
\draw[->](plus 3) to (right cut);
\draw[->](right cut) to (right out);
\draw[->](right prop 21) to (right epsilon 21);
\draw[->](right prop 22) to (right epsilon 22);

\draw[->](right epsilon 21) to (right plus 2);
\draw[->](right epsilon 22) to (right plus 2);
\draw[->](right plus 2) to (right alpha 2);

\draw[red,->](tautology 22) to [out=180, in = 90] (right prop 21);
\draw[red,->](tautology 22) to [out=20, in = 120] (right out prop);
\draw[red,->](tautology 22) to [out=20, in = 120] (left out prop);
\draw[red,->](tautology 22) to [out=0, in = 90] (right prop 22);

\draw[->](left prop 3) to (left alpha 3);
\draw[->](right prop 3) to (right alpha 3);
\draw[->](left alpha 3) to (left epsilon 3);
\draw[->](right alpha 3) to (right epsilon 3);
\draw[->](left epsilon 3) to (plus 3);
\draw[->](right epsilon 3) to (plus 3);

\draw[red,->](tautology 3) to [out=180, in = 90] (left prop 3);
\draw[red,->](tautology 3) to [out=0, in = 90] (right prop 3);
\draw[red,->](right epsilon 3) to [out=120, in = 20] (left alpha 3);

\end{tikzpicture}
\]

After an application of \textsc{Dup}, we arrive at the net

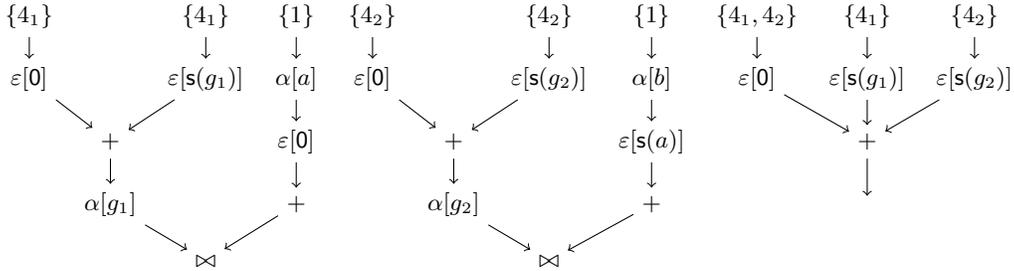
\[
\begin{tikzpicture}
\matrix[row sep=3mm,column sep=2mm]{
\node (left prop 1) {$\finset{4_1}$}; && \node (right prop 1) {$\finset{4_1}$};
&\node (left prop 3) {$\finset{1}$}; &
\node (left prop 2) {$\finset{4_2}$}; && \node (right prop 2) {$\finset{4_2}$};&  
 \node (right prop 3) {$\finset{1}$}; &  \node(left prop out){$\finset{4_1,4_2}$};
& \node(mid prop out){$\finset{4_1}$}; &\node(right prop out){$\finset{4_2}$};\\
 \node (left epsilon 1) {$\varepsilon[\mathsf{0}]$};
&&\node (right epsilon 1) {$\varepsilon[\mathsf{s}(g_1)]$}; & \node (left alpha 3) {$\alpha[a]$};&
 \node (left epsilon 2) {$\varepsilon[\mathsf{0}]$};&&\node (right epsilon 2) {$\varepsilon[\mathsf{s}(g_2)]$};&
 \node (right alpha 3) {$\alpha[b]$}; & \node(left epsilon out){$\varepsilon[\mathsf{0}]$};&
\node(mid epsilon out){$\varepsilon[\mathsf{s}(g_1)] $}; 
& \node(right  epsilon out){$\varepsilon[\mathsf{s}(g_2)]$};\\
&
 \node (left plus 1) {+}; &&\node (left epsilon 3) {$\varepsilon[\mathsf{0}]$};
&& \node (left plus 2) {+}; &&\node (right epsilon 3) {$\varepsilon[\mathsf{s}(a)]$}; && \node(plus out){$+$};
 \\
 &
 \node (left alpha 1) {$\alpha[g_1]$}; &&\node(left plus 3){+};&& \node (left alpha 2) {$\alpha[g_2]$};&&
\node(right plus 3){$+$}; && \node(out out){};\\
&& \node (left cut) {$\bowtie$};&&&&\node (right cut) {$\bowtie$};\\
&& \node (left out) {}; &&&& && \node (right out) {};\\
};

\draw[->] (left prop out) to (left epsilon out);
\draw[->] (mid prop out) to (mid epsilon out);
\draw[->] (right prop out) to (right epsilon out);
\draw[->] (left epsilon out) to (plus out);
\draw[->] (mid epsilon out) to (plus out);
\draw[->] (right epsilon out) to (plus out);
\draw[->] (plus out) to (out out);
\draw[->] (left prop 1) to (left epsilon 1);
\draw[->] (left prop 2) to (left epsilon 2);
\draw[->] (right prop 1) to (right epsilon 1);
\draw[->] (right prop 2) to (right epsilon 2);
\draw[->] (left epsilon 1) to (left plus 1);
\draw[->] (left epsilon 2) to (left plus 2);
\draw[->] (right epsilon 1) to (left plus 1);
\draw[->] (right epsilon 2) to (left plus 2);
\draw[->] (left plus 1) to (left alpha 1);
\draw[->] (left plus 2) to (left alpha 2);
\draw[->] (left alpha 1) to (left cut);
\draw[->] (left alpha 2) to (right cut);

\draw[->] (left prop 3) to (left alpha 3);
\draw[->] (right prop 3) to (right alpha 3);
\draw[->] (left alpha 3) to (left epsilon 3);
\draw[->] (right alpha 3) to (right epsilon 3);
\draw[->] (left epsilon 3) to (left plus 3);
\draw[->] (right epsilon 3) to (right plus 3);
\draw[->] (left plus 3) to (left cut);
\draw[->] (right plus 3) to (right cut);




\end{tikzpicture}
\]

\noindent and eliminating the right-hand cut, we obtain
 
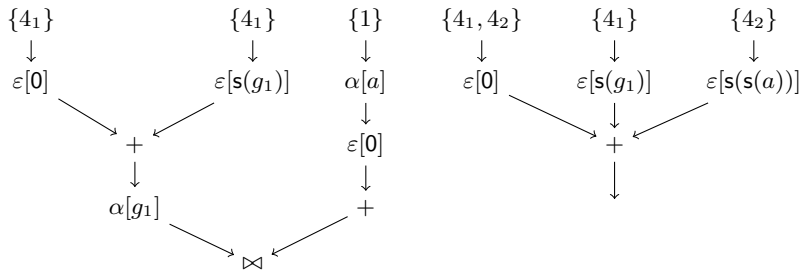
\[
\begin{tikzpicture}
\matrix[row sep=3mm,column sep=5mm]{
\node (left prop 1) {$\finset{4_1}$}; && \node (right prop 1) {$\finset{4_1}$};
&\node (left prop 3) {$\finset{1}$}; &  \node(left prop out){$\finset{4_1,4_2}$};
& \node(mid prop out){$\finset{4_1}$}; &\node(right prop out){$\finset{4_2}$}; \\
 \node (left epsilon 1) {$\varepsilon[\mathsf{0}]$};
&&\node (right epsilon 1) {$\varepsilon[\mathsf{s}(g_1)]$}; & \node (left alpha 3) {$\alpha[a]$};&
 \node(left epsilon out){$\varepsilon[\mathsf{0}]$};&
\node(mid epsilon out){$\varepsilon[\mathsf{s}(g_1)] $}; 
& \node(right  epsilon out){$\varepsilon[\mathsf{s}(\mathsf{s}(a))]$};\\
& \node (left plus 1) {$+$}; && \node (left epsilon 2) {$\varepsilon[\mathsf{0}]$};
&& \node(plus out){$+$};\\
 &
 \node (left alpha 1) {$\alpha[g_1]$}; &&\node(left plus 3){+};&& \node(out out){};\\
&& \node (left cut) {$\bowtie$};&&&&\\
&& \node (left out) {}; &&&& && \node (left out) {};\\
};

\draw[->] (left prop out) to (left epsilon out);
\draw[->] (mid prop out) to (mid epsilon out);
\draw[->] (right prop out) to (right epsilon out);
\draw[->] (left epsilon out) to (plus out);
\draw[->] (mid epsilon out) to (plus out);
\draw[->] (right epsilon out) to (plus out);
\draw[->] (plus out) to (out out);
\draw[->] (right prop 1) to (right epsilon 1);
\draw[->] (left prop 1) to (left epsilon 1);
\draw[->] (left epsilon 1) to (left plus 1);
\draw[->] (right epsilon 1) to (left plus 1);
\draw[->] (left plus 1) to (left alpha 1);
\draw[->] (left prop 3) to (left alpha 3);
\draw[->] (left alpha 3) to (left epsilon 2);
\draw[->] (left epsilon 2) to (left plus 3);
\draw[->] (left plus 3) to (left cut);
\draw[->] (left alpha 1) to (left cut);




\end{tikzpicture}
\]

\noindent The only expansion tree of our net already has three branches, and we have not finished 
cut-elimination.  In particular, it remains to evaluate the eigenvariable $a$,
and there are two witnesses with which we can evaluate it.  Thus the
final result of the cut-elimination is a net with \emph{four}
witnesses rather than three:
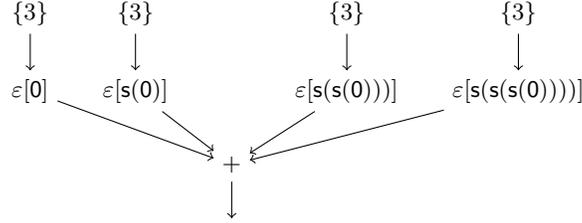
\[
\begin{tikzpicture}
\matrix[row sep=5mm,column sep=5mm]{
\node (output tautology 1){$\finset{3}$}; & \node (output tautology 2){$\finset{3}$};&& \node (output tautology 3){$\finset{3}$};& \node (output tautology 4){$\finset{3}$}; \\
\node (epsilon 1) {$\varepsilon[\mathsf{0}]$}; &
\node (epsilon 2) {$\varepsilon[\mathsf{s}(\mathsf{0})]$}; && \node (epsilon 3) {$\varepsilon[\mathsf{s}(\mathsf{s}(\mathsf{0})))]$};& \node (epsilon 4) {$\varepsilon[\mathsf{s}(\mathsf{s}(\mathsf{s}(\mathsf{0}))))]$}; \\
&&\node (plus) {+};\\
&&\node (out) {};\\
};
\draw [->] (output tautology 1) to (epsilon 1);
\draw [->] (output tautology 2) to (epsilon 2);
\draw [->] (output tautology 3) to (epsilon 3);
\draw [->] (output tautology 4) to (epsilon 4);
\draw [->] (epsilon 1) to (plus);
\draw [->] (epsilon 2) to (plus);
\draw [->] (epsilon 3) to (plus);
\draw [->] (epsilon 4) to (plus);
\draw [->] (plus) to (out);
\end{tikzpicture}
\]

Hence minimal reduction in Herbrand nets is not confluent.

\section{Other kinds of reduction}
The notion of a kingdom took a lot of effort to define, and is
(unfortunately) little known outside the community of specialists in linear
logic proof nets.  In this section we address (and reject) two possible
alternatives.
\subsection{Copying too little: dependent subforests}
Given an annotated sequent of the form
\[  F, \alpha[a].t \bowtie s_1+s_2 :A \bowtie \bar{A} \]
\noindent if we are to copy the subterm $\alpha[a].t$, to provide two
copies to cut against $s_1$ and $s_2$, what is the smallest subforest
(not necessarily a subnet) we must duplicate in order to still have an
annotated sequent?  A little thought suggests the \emph{dependent
  subforest}, consisting of all the subterms $t'$ such that
$\alpha[a].t \vartriangleleft t'$.  Since subnets are also closed
under dependency, we would never copy more than the kingdom, but in
general we copy much less. In addition, since the tautology jumps 
play no part in the dependency relation, we can
simply drop them, (being sure to replace the condition on being an
Herbrand structure with some other tautology checking condition).

Such a reduction was the subject of study by the author, and
independently by Willem Heijltjes (and others before us); it is
seductively simple and holds the promise of an elegant abstract
representation of classical proofs, but the system has a fatal flaw:
as observed by Heijltjes, by duplicating dependent subforests we may
reduce the example from the previous section to a forest containing a
cut of the following shape, where there is a jump ``across the cut'':

\begin{equation}
\label{eq:garbage}
\alpha[a] \bowtie \varepsilon[M(a)]
\end{equation}

Such a ``proof''    can, of course, never arise as the annotation of a sequent derivation.
This suggests, as is indeed the case, that the dependent-subforest duplicating reduction
does not preserve the property of being an Herbrand net.  

While we rejected this reduction in favour of Minimal reduction, which
does preserve the property, Heijltjes opts instead to treat such
redices as appear in \eqref{eq:garbage} as ``garbage'', and adds an
extra garbage collection reduction to remove them.  Since the
structure at tautology nodes is not needed for dependent subforest
duplication, Heijltjes's ``Proof Forests'' can derived from our
annotated sequents by forgetting the structure at the leaves.  His
correctness criterion is such that (the forgetful projection of) any
Herbrand structure is a correct Proof Forest.  Moreover, his weakly
normalizing reduction seems to yield the same results as ours, since
it always reduces an $\ll$-topmost cut (where the kindom and dependent
subforest coincide).  Nonetheless, there are correct Proof Forests
containing no ``garbage'' redices and yet corresponding to no
sequent-derivation. 

In the way they behave and are handled, Heijltjes's forests are
rather similar to Lamarche and Strassburger's proof nets 
for propositional classical logic~\cite{LamStra05NamProCla}.  We consider them an
interesting parallel strand of research to our own. 

\subsection{Copying too much: empires}
As mentioned above, the very natural concept of kingdom is little-mentioned
in the proof-net literature.  The concept of empire, by contrast, appears
in almost all introductions to the theory of proof nets for \MLLminus,
and played a central role in their development.  Moreover, the empire
of a node is easy to calculate; for \MLLminus nets, for example, it can
be calculated in time linear in the size of the net.

It is natural to ask, therefore, if this more familiar notion can be
the basis of a cut-elimination for Herbrand nets.  The following
counterexample shows this is not possible.  Let the underlying theory
be as for the strong normalization counterexample, and let $B =
\exists z. (\lnot \mathsf{iszero}(z)$.  In the following net, the
shaded subnet is the \emph{copyable part} of the empire of
$\alpha[g]$; the largest subnet of the empire of $\alpha[g]$ whose
roots, other than $\alpha[g]$, are all cuts or naked witnesses.

   
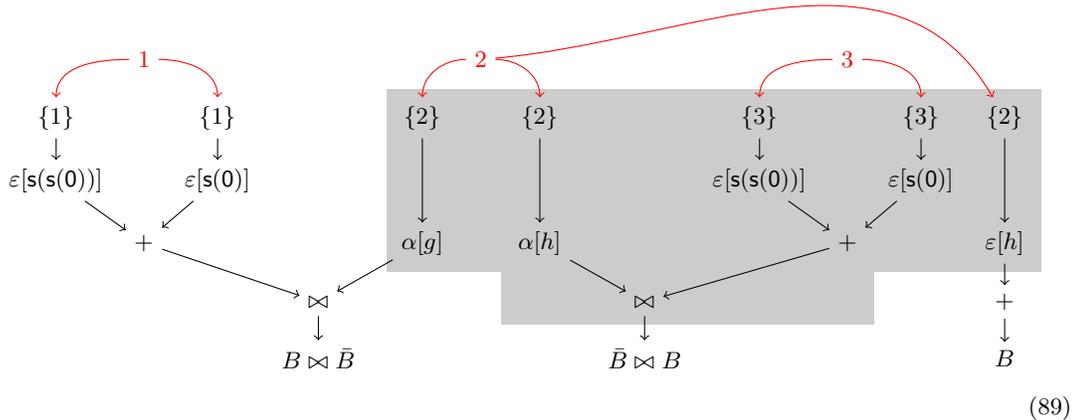
\begin{equation}
\label{eq:empireexample}
\begin{tikzpicture}
\matrix[row sep=3mm,column sep=2mm]{
& \node [red](tautology) {$1$}; &&&&& 
\node [red](tautology 2) {$2$}; &&&&& 
\node [red](tautology 3) {$3$}; \\
 \node (left prop) {$\finset{1}$}; &&
 \node (right prop) {$\finset{1}$}; &&& \node (left prop 2) {$\finset{2}$}; &&
 \node (right prop 2) {$\finset{2}$};&&& \node (left prop 3) {$\finset{3}$}; &&
 \node (right prop 3) {$\finset{3}$}; & \node(prop out){$\finset{2}$};\\
\node (left epsilon) {$\varepsilon[\mathsf{s}(\mathsf{s}(\mathsf{0}))]$}; &&
 \node (right epsilon) {$\varepsilon[\mathsf{s}(\mathsf{0})]$}; &&&
 &&
 &&&
\node (left epsilon 3) {$\varepsilon[\mathsf{s}(\mathsf{s}(\mathsf{0}))]$}; &&
 \node (right epsilon 3) {$\varepsilon[\mathsf{s}(\mathsf{0})]$}; &\\
 & \node(plus){$+$}; &&&&
\node (left alpha 2) {$\alpha[g]$}; &&
 \node (right alpha 2) {$\alpha[h]$}; &&\node (nix) {};&&
\node(plus 3){$+$}; && \node (epsilon out) {$\varepsilon[h]$};\\
&&& \node (left cut) {$\bowtie$}; &&&&&& \node (right cut) {$\bowtie$};
&&&&   \node(plus out){$+$};\\
&&& \node (left out) {$B \bowtie \bar{B}$}; &&&&&& \node (right out){$\bar{B} \bowtie B$}; &&&& 
\node(out out){$B$};\\
};



\draw[->] (prop out) to (epsilon out);
\draw[->] (epsilon out) to (plus out);
\draw[->] (plus out) to (out out);
\draw[->](plus) to (left cut);
\draw[->](left alpha 2) to (left cut);
\draw[->](right alpha 2) to (right cut);
\draw[->](plus 3) to (right cut);
\draw[->](right cut) to (right out);
\draw[->](left cut) to (left out);
\draw[->](left prop) to (left epsilon);
\draw[->](right prop) to (right epsilon);
\draw[->](left epsilon) to (plus);
\draw[->](right epsilon) to (plus);

\draw[->](left prop 2) to (left alpha 2);
\draw[->](right prop 2) to (right alpha 2);

\draw[red,->](tautology 2) to [out=180, in = 90] (left prop 2);
\draw[red,->](tautology 2) to [out=0, in = 90] (right prop 2);
\draw[red,->](tautology 2) to [out=5, in = 120] (prop out);

\draw[red,->](tautology) to [out=180, in = 90] (left prop);
\draw[red,->](tautology) to [out=0, in = 90] (right prop);

\draw[->](left prop 3) to (left epsilon 3);
\draw[->](right prop 3) to (right epsilon 3);
\draw[->](right epsilon 3) to (plus 3);
\draw[->](left epsilon 3) to (plus 3);

\draw[red,->](tautology 3) to [out=180, in = 90] (left prop 3);
\draw[red,->](tautology 3) to [out=0, in = 90] (right prop 3);

\begin{pgfonlayer}{background}
\node [fill=black!20, fit = (epsilon out) (left prop 2)] {};
\node [fill=black!20,fit=(right alpha 2)(plus 3) (right cut)] {};
\end{pgfonlayer}
\end{tikzpicture}
\end{equation}

The reader can verify that, if this subnet is copied in the obvious
way, and the resulting \textsc{Comm}/\textsc{Prop} redices reduced, the
resulting net contains \ref{eq:empireexample} as a subnet, and indeed,
it is not hard to prove that this net has no finite sequence of
reductions, ending in a cut-free net, if we replace the minimal
Duplication with the duplication of the (copyable part of) the empire.

\section{Conclusions and further work}
\label{sec:concl-furth-work}
We shown, in this paper, a system of proof nets for classical
first-order logic in prenex normal form, derived from Herbrand's
theorem.  The system has the minimal set of properties one might
expect of a proof system for classical logic --- like Gentzen's \LK it
has weakly normalizing cut-elimination.  We hope, of course, for more.
Surprisingly, given the polarization of connectives, (and thus the
avoidance of the contraction-contraction and weakening-weakening
problems detailed in \cite{Girard91NewCon}) cut-reduction in this
system is nonconfluent (a counterexample for Heijltjes' system, also
applicable to our system, was given in \cite{Hei08ProFor}).  We seek,
therefore, confluent subsystems.  We conjecture, but as yet have no
proof, that minimal reduction is strongly normalizing.  

Similar structures to our annotated sequents arise as strategies
Coquand's game theoretical treatment of classical
arithmetic~\cite{Coq95SemEvi}.  Coquand gives a way to play a strategy
containing cuts, which amounts to a non-associative composition on
proofs, and it would be interesting to compare this with the
nonconfluent properties of Herbrand nets.

We look also to extend our system beyond prenex normal form, first to
encompass a treatment of the propositional connectives.  The paper
\cite{Mck09ExpNet} gives a multiplicative treatment of classical proof
nets which improves on \cite{Rob03ProNetCla} by replacing contraction
(binary, defined on all formulae) by expansion (n-ary, defined only on
positive formulae).  Contraction on negative \emph{atoms} (needed for
completeness) is handled by the same basic binding structure used here
to model quantification.

\noindent\textbf{Acknowledgements} The author thanks Willem Heijltjes
for many stimulating and helpful exchanges, and thanks Michel Parigot,
Lutz Strassburger, Kai Br\"unnler, Roman Kuznets and Stefan Hetzl for
usefule comments.  

\nocite{Girard91NewCon,Coq95SemEvi,BelSco94piCal} \bibliography{nets}
\appendix

\section{Properties of subnets of Herbrand nets}
\label{sec:prop-herbr-nets}
The proofs contained in the appendix are very minor variations on the
proofs of similar properties for \MLLminus proof nets,
as presented in \cite{212898}.  They are presented here for the sake
of completeness.

\subsection{Existence of kingdom and empire}
\begin{definition}
  Let $F$ be an ACC forest, $t$ a node of $F$, and $\sigma$ a
  switching of $F$.  Remove from $F_\sigma$ the edge from $t$ to its
  parent, if it has one.  $F(t, \sigma)$
  is the connected component of this graph containing $t$.
\end{definition}

\begin{proposition}
Let $e = \bigcap_\sigma F(X, \sigma)$, where $\sigma$ ranges over 
all switchings of $F$.  Let $e(X)$ be the intersection of $e$ with the
nodes of $F$.  $e(X)$ is a subnet of $F$, and $X$ is a root 
of $e(X)$.  
\end{proposition}
\begin{proof}
  We must first see that $e(X)$ is a substructure of $F$ -- that is,
  it must be closed under $\vartriangleleft$.  This is easy to see
  when passing from an unswitched node to its unique child.  Suppose
  now that $Z$ is a switched node in $e(X)$, and that one of its
  $\vartriangleleft$-successors $Y$ is not in $e(X)$.  Then there is a
  switching $\sigma$ such that $Z \in F(\sigma, X)$ and $Y \notin
  F(\sigma, X)$.  Thus there is a path $p$ $X$ to $Z$ in $F_\sigma$,
  and a path $p'$ from the parent $W$ of $X$ to $Y$, also in
  $F_\sigma$.  By changing the switching $\sigma$ to a switching
  $\sigma'$, where $Z$ chooses $Y$ and $W$ chooses $X$ (if $W$ is
  switched) and leaving all other switches unchanged, we obtain a
  cyclic swiyching graph $F_\sigma'$.  Hence $e(X)$ is a substructure.

  We next observe that $e(X)$ is an ACC forest: let $\sigma$ be a switching
of the nodes in $e(X)$, and let $\sigma'$ be an extension of that
  switching to $F$.   The graph $e(X)_\sigma$ is acyclic since $e(X)$ is a
substructure of $F$.  To see that $e(X)_\sigma$ is connected, observe
that it is the restriction of $F(X, \sigma')$, a connected graph, to 
$e = \bigcap_\sigma F(X, \sigma)$.

Suppose now that $X$ is not a root of $e(X)$.  Then there is a $Y$ in
$e(X)$ such that $Y \leftarrow X$.  Choose a switching $\sigma_X$ of
$F$ such that whenever $Z$ is a switched node with $Y \leq Z
\leq X$, we choose a switching $W$ for $Z$ such that $W$ is the
predecessor of  $X$.

  Because of these choices, the unique path from $X$ to $Y$ in
  $F_{\sigma_X}$ uses the edge from $X$ to its parent, and because of
  this does not provide a path from $X$ to $Y$ in $F(X, \sigma_X)$.
  If $Y$ is in $e(X)$, then there is some other path from $X$ to $Y$
  in $F_{\sigma_X}$, but this contradicts the fact that $F$ is correct
  (acyclicity of $F_{\sigma_X}$).
\end{proof}

\begin{proposition}
The subnet $e(X)$ is the largest subnet of $F$ having $X$ as a root.
\end{proposition}
\begin{proof}
  Suppose otherwise.  Let $G$ be a substructure of $F$, with $X$ as a
  root, which is larger than $e(X)$.  Then there is a node $Z$ of $G$,
  and a switching $\sigma$, such that $Z \notin F(\sigma, X)$. But
  then there is no path from $X$ to $Z$ in $G_\sigma$, and so $G$ is
  not an ACC forest.
\end{proof}

The following technical lemma will be
crucial:
\begin{lemma}
  Let $F$ be an Herbrand net, and let $s$ and $t$ be distinct nodes of
  $F$,  such that $t \in e(s)$.  Let
  $s'$ be the parent of $s$ and $t'$ the parent of $t$.  Then
\[ s' \in e(t) \text{ \ iff \ } t' \notin k(s')\]
\label{lem:empire-kingdom-nesting}
\end{lemma}
\begin{proof}
We have that 
\[ G_1 = e(t) \cap k(s') \qquad G_2 = e(t) \cup k(s')\]
are nets (since $G_1$ is nonempty).  If $s' \in e(t), t' \in k(s')$
then $G_1$ has $s'$ as a root and does not contain $t'$, and so is a 
subnet with $s'$ as a root smaller than $k(s')$ -- contradiction.

Similarly, if $t' \notin e(s), \ s' \notin k(t')$ then 
$G_2$ has $t$ as a root and contains $s'$, in contradiction of
the definition of empire.
\end{proof}

This allows us to show that the relation $\ll$ is a partial order on
the nodes of a structure.
\begin{lemma}
\label{lem:ll-is-an-order}
  Let $F$ be an Herbrand net, and let $X$, $Y$ be nodes of $F$ such
  that $X\ll Y$ and $Y \ll X$.  Then $X = Y$.
\end{lemma}
\begin{proof} 
  Let $X$ be labelled with $t$ and $Y$ with $s$.  Suppose that $X$ and
  $Y$ are not the same node.  We have that $k(X)= k(X) \cap k(Y) =
  k(Y)$, by minimality of the kingdom.
\begin{enumerate}
\item If $X$ is an $\alpha$ node, or expansion node, 
then removing $X$ from $k(Y)$ yields a smaller subnet with $Y$ as
  a root, contradicting minimality of $k(Y)$.
\item If $X$ is an $\varepsilon$ node with child $X'$, then its
  kingdom is equal to $k(X') \cup \finset{X}$, and so $Y \in k(X')$.
This contradicts the previous lemma, which says that $Y \notin e(X')$.
Similarly for $\bowtie$ nodes.
\end{enumerate} 
\end{proof}

$e(t)$ and $e(s)$ are clearly disjoint. Suppose
  that a node $ u \in F$ is a member of $e(t)$ but that $u'$, the
  parent of $u$, is not in $e(t)$.  By Lemma \ref{lem:empire-kingdom-nesting},
  $t\bowtie s$ is a member of $k(u')$, contradicting $\ll$-minimality of $t \bowtie s$.
  By connectedness of $F$, we have that $F = e(t)\cup e(s) \cup
  \finset{t \bowtie s}$. 

\subsection{Calculating the kingdom}

We know that the kingdom of a node always exists, but the definition
of the kingdom of $X$ as the intersection of all subnets having $X$ as
a root is unwieldy for calculations.  In this section we reconstruct
arguments from \cite{212898} showing that the kingdom of a node may be
calculated in time at most quadratic in the size (number of nodes) of
a net.

We will work on the dependency graph of a net.  We first see
how to calculate the empire of a node:
\begin{lemma}
  Let $F$ be an ACC forest and $X$ a node of $F$.  The graph $e$, defined
  above, is the smallest subgraph $E$ of $Dep(F)$ of $F$ closed under
  the following:
\begin{enumerate}
\item $X \in E$
\item \label{item:0}(Dependency) If a node $X$ is in $E$ then 
all vertices $Y$ with $X\leftarrow Y$ or $X \dependson Y$ are in $E$.
\item \label{item:2}An  unswitched node $t$ of $F$ is in $E$ if and only
  if there is a $\vartriangleleft$-predecessor $s$ of t also in $E$.
\item \label{item:3}An expansion node $t$ of $F$ is in $E$ if and only
  if \emph{all} its $\vartriangleleft$-predecessors are in $E$.
\end{enumerate}
\end{lemma}
\begin{proof}
  We have already seen that $e$ is closed under dependency: 
  it is easy to deduce that $e$ is also closed under items 
  \ref{item:2} and \ref{item:3}, and thus $E \subset e$.

  To see that $e \subseteq E$, we construct a switching $\sigma$ for
  $F$ such that $E= F(\sigma,X)$.  Choose the switching $\sigma$ ---
  the \emph{principal switching} for $X$ --- such that, whenever a switched
  node $Y$ has a choice of switching $Z$ which is not in $E$, we pick
  that switching.

  It follows from the properties of $E$ that each root of $E$ is
  either $t$, a root of $F$ or an expansion node $X$ one of whose
  children is not in $E$.  Now suppose that a node $W$ is in $e$
  but not in $E$.  Since $F(\sigma, X)$ is connected and $E$ is
  closed under dependency, the path connecting $X$ and $W$ in $F_\sigma$
  must exit $E$ at one of its roots.  But this is impossible by the choice
of a principal switching.  Thus $e \subseteq E$.
\end{proof}

\begin{corollary}
The empire $e(X)$ of a node $X$ of an ACC forest may be calculated in
a number of steps linear in the number of nodes in $F$.
\end{corollary}

We may now use Lemma~\ref{lem:empire-kingdom-nesting} to give 
an alternative characterization of the kingdom $k(X)$ of a node:

\begin{lemma}
Let $K$ be the smallest subset of vertices of $Dep(F)$ containing $X$ 
and closed as follows
\begin{enumerate}
\item $K$ is closed under dependency.
\item If $Z$ is a successor of $Y \in K$, and $Y, Z \neq X$, then
$Z \in k(X)$ if and only if $X \notin E$. 
\end{enumerate}
The non-tautology vertices of $K$ are precisely the nodes of $k(X)$.
\end{lemma}
  
\begin{corollary}
The kingdom $k(X)$ of a node $X$ may be calculated in 
time quadratic in the number of nodes of $F$.  
\end{corollary}

\end{document}